\documentclass[12pt]{article}
\renewcommand{\include}{\input}
\usepackage{amsmath,amsthm,verbatim,amssymb,amsfonts,amscd,diagrams,graphics}
\theoremstyle{plain}
\newtheorem{theorem}{Theorem}[section]
\newtheorem{corollary}[theorem]{Corollary}
\newtheorem{lemma}[theorem]{Lemma}
\newtheorem{proposition}[theorem]{Proposition}

\theoremstyle{definition}

\theoremstyle{remark}
\newtheorem{remark}[theorem]{Remark}

\newarrow{ul}---->
\newarrow{Backwards}<----

\newcommand{\td}[1]{\tilde{#1}}
\newcommand{\into}{\hookrightarrow}
\newcommand{\Z}{\mathbb{Z}}
\newcommand{\Q}{\mathbb{Q}}
\newcommand{\R}{\mathbb{R}}
\newcommand{\bd}{\partial}

\newcommand{\la}{\leftarrow}

\newcommand{\mv}{\mathversion{bold}}

\newcommand{\mc}[1]{\mathcal{#1}}
\newcommand{\bb}[1]{\mathbb{#1}}
\newcommand{\dlim}{\varinjlim}
\newcommand{\vg}{\varGamma}

\newcommand{\mf}{\mathfrak}

\begin{document}

\title{Intersection Alexander polynomials}
\author{Greg Friedman\\Yale University\\Dept. of Mathematics\\10 Hillhouse Ave\\PO Box 208283\\New Haven, CT 06520\\friedman@math.yale.edu\\Tel. 203-432-6473  Fax:  203-432-7316}
\date{April 29, 2002}
\maketitle

\begin{abstract}
By considering a (not necessarily locally-flat) PL knot as the singular locus of a PL stratified pseudomanifold, we can use intersection homology theory to define \emph{intersection Alexander polynomials}, a generalization of the classical Alexander polynomial invariants for smooth or PL locally-flat knots. We show that the intersection Alexander polynomials satisfy certain duality and normalization conditions analogous to those of ordinary Alexander polynomials, and we explore the relationships between the intersection Alexander polynomials and certain generalizations of the classical Alexander polynomials that are defined for non-locally-flat knots. We also investigate the relations between the intersection Alexander polynomials of a knot and the intersection and classical Alexander polynomials of the link knots around the singular strata.  To facilitate some of these investigations, we introduce spectral sequences for the computation of the intersection homology of certain stratified bundles. 

\end{abstract}



\tableofcontents
\section{Introduction}

In \cite{CS}, Cappell and Shaneson adapted the theory of intersection
homology
to study stratified manifolds with embedded subpseudomanifolds as the
singular
loci. In particular, if we consider the case of a knot $K$ given by a
PL-sphere pair
$S^{n-2}\subset S^n$, we can think of $S^n$ as a stratified manifold with
singular locus
$K\cong S^{n-2}$. If we choose a  local coefficient system $\vg$
defined
on $S^n-K$ with stalks $\Gamma=\Q[\Z]=\Q[t,t^{-1}]$ and action of the
fundamental group given so that $\alpha\in\pi_1(S^n-K)$ acts on
$\Gamma$ by multiplication by $t^{\ell_K(\alpha)}$, where $\ell_K(\alpha)$ is
the linking number of $\alpha$ with the knot $K$ (see \cite[Section 4.3]{GBF1}), then
following \cite{GM2}
and \cite{CS}, we can define the intersection homology groups
$IH^{\bar p}_i(S^n;\vg)$. These will be torsion $\Gamma$-modules, and we can study
their associated polynomials $I\lambda_i^{\bar p} (t)$, which, following Cappell and Shaneson, we
will call the
\emph{intersection Alexander polynomials}. As we shall see, these polynomials
possess interesting properties of their own as well as some relations with the
ordinary Alexander polynomials of non-locally-flat knots as studied in the author's dissertation (see \cite{GBF} and \cite{GBF1}).

The structure of this paper is as follows: 

Section \ref{S: poly alg} consists of some algebraic preliminaries that will be of use, while 
Sections \ref{S: IH con} and \ref{S: superperversities} contain some introductory
material concerning our approach to intersection homology theory. 

In Section \ref{S: IH dual}, we apply the superduality results of Cappell and Shaneson
\cite{CS} to show that the intersection Alexander polynomials possess a
duality analogous to that for traditional Alexander polynomials:

\begin{theorem}[Theorem \ref{T: superduality}]\label{T: dual*}
Let $K\cong S^{n-2}\subset S^n$ be a knot, not necessarily locally-flat, and let $\bar p$
and $\bar q$ be a superdual perversity and superperversity as defined in \cite{CS}, i.e.  $\bar p(k)+\bar
q(k)=k-1$ for all $k\geq 2$. Then $I\lambda^{\bar p}_i(t)\sim I\lambda^{\bar q}_{n-1-i}(t^{-1})$, where $\sim$ denotes similarity in $\Gamma$.
\end{theorem}

Note that (again following Cappell and Shaneson) we allow intersection homology
modules defined by the Deligne process (\cite{GM2},\cite{Bo}) to have
\emph{superperversities}, i.e. perversities $\bar q$ such that $\bar
q(2)=1$. In Section \ref{S: superperversities}, we show that these
superperverse intersection homology modules do \emph{not} necessarily
agree with those obtained from the geometric
(simplicial or singular) intersection homology theories. By contrast, it
is known that the theories do agree for traditional perversities,
$\bar p$, which satisfy $\bar p(2)=0$.
Theorem \ref{T: dual*} allows us to concentrate on the traditional perversity
polynomials for the remainder of the paper, since the corresponding
results for the superperverse polynomials can be obtained by duality.  Thus for the remainder of this introduction, the perversity $\bar p$ refers to a traditional perveristy with $\bar p(2)=0$.

In Section \ref{S: IH norm}, we obtain a normalization condition on the $I\lambda_i^{\bar p}$.
Suppose that $\gamma$ is an element of $\Gamma$. Recall that there is an
element, say $\bar \gamma$, in the similarity class of $\Gamma$ which is primitive
in
$\Lambda=\Z[\Z]=\Z[t,t^{-1}]$, i.e. the coefficients are relatively prime
(though not necessarily pairwise so), and this element is unique up to
similarity class in $\Lambda$ (see, e.g., \cite{L66} or \cite{GBF1}). We will say that $\gamma$ is a polynomial of
Alexander type if $\bar\gamma(1)=\pm 1$. Another classical property of the Alexander polynomials of locally-flat knots $S^{n-2}\subset S^n$ is that they are of Alexander type for $0<i<n-1$ (\cite{L66}). In fact, this is also true of the Alexander polynomials of knots which are not locally flat (\cite{GBF1}, \cite{GBF}), where in this case the polynomials are again defined to be those associated to the modules $H_i(S^n-K;\vg)$. We  show that the same property holds for intersection Alexander polynomials:
\begin{theorem}[Theorem \ref{T: normal}]
For any PL-knot $K\cong S^{n-2}\subset S^n$, not necessarily locally-flat, and traditional
perversity $\bar p$,
$I\lambda_i^{\bar p}$ is of Alexander type for $i>0$, $I\lambda_0^{\bar p}\sim t-1$,
and $I\lambda_i^{\bar p}\sim 1$ for $0\neq i\geq n-1$.
\end{theorem}
\noindent This is proven
by a double induction on the dimension of the knot and the codimension of
the strata, using the results of \cite{GBF1} on the ``ordinary'' Alexander polynomials of
knots as the ``base step''.

In Section \ref{S: ih and ord}, we compare the intersection Alexander
polynomials with the ``ordinary''
Alexander
polynomials in several interesting cases. For locally-flat knots, we see, not surprisingly, that they
are
identical, i.e. $I\lambda_i^{\bar p}\sim \lambda_i$. 
For a knot
with a point
singularity, we obtain a more complicated relationship between the intersection Alexander
polynomials and the ordinary Alexander polynomials. 
It will be useful to introduce some notation from \cite{GBF} and \cite{GBF1}. Let $D$ be the open regular neighborhood of the singular point of the embedding. Then the complement of $D$ in the pair $(S^n, S^{n-2})$ is a locally flat disk knot bounded by a locally-flat sphere knot, which is the link knot around the singular point. Futhermore, the knot complement $S^n-K$ is homotopy equivalent to the complement of the induced disk knot (see \cite[\S 3]{GBF1} for details). We can then define Alexander polynomials $\lambda_i$, $\nu_i$, and $\mu_i$ associated to the homology modules of the disk knot complement $C$, the boundary sphere knot complement $X$, and the pair $(C,X)$. Note that the  $\lambda_i$ also represent the Alexander polynomials of the sphere knot $K$ and the $\nu_i$ are just the usual Alexander polynomials of the locally-flat link knot. Furthermore, these polynomials factor into terms which they share, i.e. we can write $\nu_i\sim a_ib_i$, $\lambda_i\sim b_ic_i$ and $\nu_i\sim c_ia_{i-1}$. (Analogous polynomials can be defined for any PL knots; see \cite{GBF1}, \cite{GBF}, and Section \ref{S: man triv} below.) In this language, we obtain the following formula:

\begin{proposition}[Proposition \ref{P: ih point}]\label{P: point*}
With the notation as above, the intersection Alexander polynomial of a knot $K$ embedded with a single point singularity is given by 
\begin{equation*}
I\lambda_i^{\bar p}(t) \sim
\begin{cases}
\lambda_i(t), & i<n-1-\bar p(n)\\
c_i(t), & i=n-1-\bar p(n)\\
\mu_i(t), & i> n-1-\bar p(n).
\end{cases}
\end{equation*}
\end{proposition}
\noindent Notice how this provides  a nice example of the sort of filtering of ordinary homology theories that we often see in intersection homology as we run through a range of perversities.

More generally, if the singularity of the embedding consists of only one connected singular stratum, $\Sigma$, of any dimension, then this will be a closed manifold. In this case, the
open regular neighborhood $N(\Sigma)$ will be a block bundle pair  with blocks
of the form $D^i\times c(S^k, S^{k-2})$, the product of a disk with the open
cone on a locally-flat knot pair (the \emph{link pair}). If we make the further assumption
that
$N(\Sigma)\cong
\Sigma\times c(S^k, S^{k-2})$,  a product, then we obtain similar explicit,
but more complicated, filtering formulae:
\begin{equation*}
I\lambda_i^{\bar p}(t)\sim \mathfrak a_{i-1}^{\geq k-\bar
p(k+1)}(t) \mathfrak b_i^{<k-\bar p(k+1)}(t)c_i(t),
\end{equation*} 
where $ \mathfrak a_{i-1}^{\geq k-\bar
p(k+1)}(t)$ is a polynomial which divides $a_{i-1}(t)$ and $\mf b_i^{<k-\bar
p(k+1)}(t)$ is a polynomial which divides $b_i(t)$. See Section
\ref{S: man triv} for the exact definitions of these polynomials.
 
In Section \ref{S: SS,primes}, we continue to assume that that the singular set $\Sigma$ is a manifold, but we no longer assume that its neighborhood can be written as a product.  
We then obtain some relationships between the prime divisors of the intersection Alexander polynomials and those of both the ordinary Alexander polynomials of the knot and
to the Alexander polynomials of the locally-flat knotted link pair. At first, we assume
$N(\Sigma)$ can be given the
structure of a fiber bundle pair with fiber given by the cone on the link knot
pair (so that the neighborhood will in fact be tubular). Then, as a tool to obtain our results,
we develop the hypercohomology spectral sequence  of a Leray
intersection homology sheaf associated to the bundle neighborhood:

\begin{proposition}[Proposition \ref{P: SS}]
Let $(E,B,F, \pi)$ be a fiber bundle with base space $B$ a manifold,
total space $E$,
paracompact stratified fiber $F$, and
projection $\pi$ such that for sufficiently small open $U\subset B$,
$\pi^{-1}(U)\cong
U\times F$, where the the stratification is given by $F_i\times U$, $F_i$
the strata
of $F$. Then, for any fixed
perversity $\bar p$ which we omit from the notation, there is a spectral
sequence abutting
to the sheaf intersection cohomology $IH^i_c(E;\vg)$ with $E_2$ term
\begin{equation*}
E_2^{p,q}=H_c^p(B; \mc{IH}^q_c(F;\vg|F)),
\end{equation*}
where $\mc{IH}^i_c(F;\vg|F)$ is a local coefficient system (sheaf) with
stalks $IH^i_c
(F;\vg|F)$.
\end{proposition}

Using this proposition, we prove the following theorems at first under the assumption that the singularity $\Sigma$ has such a bundle neighborhood. However, by invoking some further results from \cite{GBF} and \cite{GBF3}, we later argue that the theorems  holds for any
knots in which the singular set of the embedding consists of only one singular stratum. We state the theorems here in their full generality; the versions with the fiber bundle assumption occur in the text as Theorems \ref{T: at one} and \ref{T: divisors}:

\begin{theorem}\label{T: at one**} 
Let $K$ be a non-locally flat knot with singularity
$\Sigma=\Sigma_{n-k-1}$ a
manifold.  Let $\xi_j$ be
the Alexander polynomials of the locally-flat link knot $\ell$ which is given by the link pair of $\Sigma$ in $S^n$. Let $\bar
p$ be
a traditional perversity. Then, for
$0<i<n-1$ and for any
prime $\gamma\in \Gamma$,
$\gamma|I\lambda_i^{\bar p}$ only if $\gamma|c_i$ or $\gamma|\xi_s$ for
some $s$ such that $0\leq i-s \leq n-k$ and $0< s < k-1$.
In other words, the prime factors of
$c_i$ and $\xi_s$, $s$ in the allowable range, are the only possible
prime factors of $I\lambda_i^{\bar p}$.
\end{theorem} 
\begin{theorem}\label{T: divisors**} 
With the same hypotheses, suppose
$\gamma$ is
a prime element
of $\Gamma$ which does not divide
$\lambda_i(t)$. Suppose $\gamma|\xi_s$ only if $s<k-\bar p(k+1)$. Then
$\gamma\nmid I\lambda_i^{\bar
p}$.
\end{theorem} 
\noindent The polynomial $c_i$ mentioned in these theorems is a factor of the polynomial corresponding to the torsion $\Gamma$-module $H_i(S^n-K;\vg)$. If $\Sigma$ is a point singularity, this is the same $c_i$ as that in Proposition \ref{P: point*}. See \cite{GBF1} and Section
\ref{S: man triv} below for more details.

 As a corollary to these theorems, we can determine some cases in which the
intersection Alexander polynomial agrees with one of the ordinary Alexander
polynomials:

\begin{corollary}[Corollary \ref{S: degeneracies}]
For a knot $K\subset S^n$ with a  manifold singularity of dimension $n-k-1$:
\begin{enumerate}
\item If $i<k-\bar p(k+1)$, then $I\lambda_i^{\bar p}\sim \lambda_i(t)$.
\item If $\bar p(k+1)\leq 1$ or if $H_i(S^k-\ell;\vg)=0$ for $i\geq j$ and $\bar
p(k+1)\leq
k-j$, then $I\lambda_i^{\bar p}\sim \lambda_i(t)$ for all $i$.
\item If $i\geq n-\bar p(k+1)+1$, then $I\lambda_i^{\bar p}\sim \mu_i(t)$.
\end{enumerate}
\end{corollary}
\noindent Once again, the polynomial $\mu_i(t)$ is analgous to that in Proposition \ref{P: point*} (see \cite{GBF1} and below).

Finally, in Section \ref{S: general knots}, we develop some relations between the
intersection Alexander polynomials of a knot, its ordinary Alexander polynomials,
and both the intersection and ordinary Alexander polynomials of
its link knots. We show the following:
\begin{theorem}[Theorem \ref{T: divisors2}]
Let $\xi_{iks}$ denote the $s$th intersection
Alexander polynomial of the link $L_{i,k}$ of the $k$th connected
component $X_{i,k}$ of the $i$th stratum $X_i=\Sigma_i-\Sigma_{i-1}$
of a knot.
A prime element $\gamma\in \Gamma$ divides the intersection
Alexander polynomial $I\lambda_j^{\bar p}$ only
if $\gamma|\lambda_j$ or $\gamma|\xi_{iks}$ for some some set of indices $i$, $k$,
and $s$
such that $0\leq j-s \leq i-1$ and $0\leq s < n-i-2$. Furthermore, $\gamma\nmid
I\lambda_j^{\bar p}$ if, for all $i,k$, $\gamma|\xi_{iks}$ only if $s<n-i-1-\bar
p(n-i)$.
\end{theorem}

\begin{theorem}[Theorem \ref{T: divisors3}]
Let $\zeta_{iks}$ be the $s$th ordinary Alexander polynomial of the link
knot pair
$L_{i,k}$. A prime element $\gamma\in \Gamma$ divides the $s$th
intersection
Alexander polynomial $I\lambda_j^{\bar p}$ only
if $\gamma|\lambda_j$ or $\gamma|\zeta_{iks}$ for some set of indices $i$, $k$,
and $s$,
such that $0\leq j-s \leq i-1$ and $0\leq s < n-i-2$.
\end{theorem}

We also obtain some results on the maximal powers to which prime divisors of the intersection
Alexander polynomials can occur:
\begin{theorem}[Theorem \ref{T: max power}]
Let $\gamma_{ipq}$ be the maximal power to which the prime $\gamma$ occurs as
a divisor of
the polynomial $e_{ipq}$ of 
$H_p(\Sigma_{i}-\Sigma_{i-1}; \mc{IH}_q(L;\vg))$ (see Section \ref{S: general knots} for a description of the local coefficient system  $\mc{IH}_q(L;\vg)$ whose fibers are the intersection homology modules of the link knots of the stratum $\Sigma_i-\Sigma_{i-1}$).
In other words, $\gamma^{\gamma_{ipq}}|e_{ipq}$, but $\gamma^{\gamma_{ipq}+1}\nmid
e_{ipq}$.
Let
$\gamma_l$ denote the maximal power to which $\gamma$ occurs in the Alexander polynomial
$\lambda_l$ of the knot $K$.  
The prime $\gamma\in \Gamma$ cannot occur in the polynomial $I\lambda_j^{\bar p}$ to a
power greater than
\begin{equation*}
\gamma_j+\sum_{i=0}^{n-2}\left((\sum_{\overset{p+q=j}{
q=0, q<n-i-1-\bar p(n-i)}}\gamma_{ipq})+(\sum_{p+q=j-1}\gamma_{ipq})\right).
\end{equation*}
\end{theorem}

This work originally appeared in the author's dissertation \cite{GBF} as part of a general program to study polynomial invariants of  PL-knots which are not locally-flat (see also \cite{GBF1}). I thank my advisor, Sylvain Cappell, for all of his generous and invaluable guidance.

\section{Polynomial algebra}\label{S: poly alg}

In this section, we provide some basic results on what we call \emph{polynomial algebra}. The polynomials in question are those associated with torsion modules over the principle ideal domain of rational Laurent polynomials $\Gamma:=\Q[t,t^{-1}]$. In particular, given such a module, we can take as the associated polynomial the determinant of a square presentation matrix or, equivalently, the product of its torsion coefficients. So if the $\Gamma$-torsion module $M$ has the form $\oplus_i \Gamma/(p_i)$, $p_i\in \Gamma$, then the associated polynomial is $\prod_i p_i$. Our main interest is in the relations that  occur between such polynomials associated to torsion modules in exact sequences. A detailed study is given in \cite{GBF} and \cite{GBF1}. Here, we simpy summarize the results that we will need and prove one additional lemma.

Let $\Gamma=\Q[\Z]=\Q[t,t^{-1}]$ be the ring of Laurent polynomials with 
rational coefficients. In other words, the elements of $\Gamma$ are
polynomials $\sum_{i\in \Z} a_it^i$ such that each $a_i\in \Q$ and
$a_i=0$ for
all but a finite number of $i$. $\Gamma$ is a principal ideal domain
\cite[\S 1.6]{L66}. Unless otherwise specified, we will generally
not distinguish between elements of $\Gamma$ and their similarity classes
up to unit. Let $\Lambda=\Z[\Z]=\Z[t,t^{-1}]$, the ring of Laurent polynomials
with integer coefficients. Then
$\Gamma=\Lambda\otimes_{\Z}\Q$. We call a polynomial in $\Lambda$
\emph{primitive} if its set of non-zero coefficients have no common
divisor
except for $\pm 1$. 
Any element of $\Gamma$ has an associate
in $\Gamma$ which is a primitive polynomial in $\Lambda$: Any
element $at^i\in \Gamma$ is a unit and, in particular then, any $a\in \Q$.
So given an element of $\Gamma$, we can first
clear denominators and then divide out any common divisors without
affecting similarity (associate) class in $\Gamma$. We will often choose
to represent an element of $\Gamma$ (technically, its associate class) by
such a primitive element of $\Lambda$. 

\begin{lemma}\label{L: prime splitting}
Let $r$ and $s$ be powers of distinct (non-associate) prime elements of
$\Gamma$. Then the
only $\Gamma$-module morphism 
$f:\Gamma/(r)\to \Gamma/(s)$ is the $0$ map.
\end{lemma}

\begin{corollary}\label{C: sequence splitting}
Let $M_i(p)$ be the $p$-primary direct summand of the torsion $\Gamma$-module $M_i$. Given an exact sequence 
\begin{equation}\label{E: ex}
\begin{CD}
0 @>d_0>> M_1 @>d_1>> M_2 @>d_2>> \cdots @>d_{n-1}>> M_n @>d_n>> 0 ,
\end{CD}
\end{equation} then for any prime $p\in \Gamma$,  the sequence 
\begin{equation}\label{E: exact2}
\begin{CD}
0 @>e_0>> M_1(p) @>e_1>> M_2(p) @>e_2>> \cdots @>e_{n-1}>> M_n(p) @>e_n>> 0 
\end{CD}
\end{equation}
is exact, where the maps $e_i$ are the restrictions of the the maps $d_i$ to the direct
summands $M_i(p)$.
\end{corollary}

Note that this lemma together with its corollary 
allows us to write the exact sequence
\eqref{E: ex} as the direct sum of exact sequences of the form \eqref{E: exact2}.

This lemma and its corollary can then be used to prove the following proposition.

\begin{proposition}\label{P:alt. poly.}
Suppose we have an exact sequence of finitely generated torsion
$\Gamma$-modules
\begin{equation}\label{E: exact1}
\begin{CD}
0 @>d_0>> M_1 @>d_1>> M_2 @>d_2>> \cdots @>d_{n-1}>> M_n @>d_n>> 0,
\end{CD}
\end{equation}
and suppose that $\Delta_i$ is the determinant of a square presentation
matrix of
$M_i$ (which we will refer to as the polynomial associated to the module). Then,  taking
$\Delta_{n+1}=1$ if $n$ is odd, the
alternating product $\prod_{i=1}^{\lceil n/2 \rceil}
\frac{\Delta_{2i-1}}{\Delta_{2i}}\in \Q(t)$ is equal to a unit of $\Gamma$,
and, in
particular, with a consistent choice of normalization within associate
classes for the elementary
divisors of the $M_i$ (in the language of \cite{HU}), this product is equal
to
$1$. 
\end{proposition}

\begin{corollary}\label{C: kern} 
With the notation and assumptions as above, each
$\Delta_i=\delta_{i}\delta_{i+1}$, where
$\delta_{i+1}|\Delta_{i+1}$ and $\delta_i|\Delta_{i-1}$. Furthermore, if we represent the $\Delta_i$ by the elements in their similarity classes in $\Gamma$  which are primitive in $\Lambda$, the $\delta_i$ will also be primitive in $\Lambda$.  
\end{corollary}

This corollary will be used often in what follows.

For convenience, we introduce the following notation. Suppose $\Delta_i\in \Gamma$. We will refer to an 
\emph{exact sequence of polynomials}, denoted by
\begin{equation*}
\begin{CD}
@>>> \Delta_{i-1} @>>>\Delta_i @>>> \Delta_{i+1} @>>>,
\end{CD}
\end{equation*}
to mean a sequence of polynomials such that each $\Delta_i\sim \delta_i \delta_{i+1}$, $\delta_i\in
\Gamma$. As we have seen, such a sequence arises  in the case of an exact sequence of
torsion
$\Gamma$-modules, $M_i$, and, in that case, the factorization of the polynomials is determined by the
maps of
the modules. In fact, each $\delta_i$ is the polynomial of the module ker$(M_i\to M_{i+1})$. 
 
Observe that knowledge of two thirds of the terms of an exact sequence of polynomials (for
example, all $\Delta_{3i}$ and $\Delta_{3i+1}$, $i\in\Z$) and the common factors of those
terms (the $\delta_{3i+1}$), allows us to deduce the missing third of the sequence
($\Delta_{3i+2}= \delta_{3i+2}\delta_{3i+3}=
\frac{\Delta_{3i+1}}{\delta_{3i+1}}\cdot\frac{\Delta_{3i+3}}{\delta_{3i+4}}$).

Note also that for any bounded exact sequence of polynomials (or even a half-bounded 
sequence),
the collections $\{\Delta_i\}$ and
$\{\delta_i\}$ carry the same information. That is, suppose that one (or both) end(s) of the
polynomial sequence is an infinite number of $1$'s (by analogy to extending any bounded or 
half-bounded exact module
sequence to an infinite number of $0$ modules).  Clearly, the $\Delta_i$ can be reconstructed
from the
$\delta_i$ by $\Delta_i\sim \delta_i \delta_{i+1}$. On the other hand, if $\Delta_0$ is the
first nontrivial term in the polynomial sequence, then $\delta_0\sim 1$, $\delta_1\sim
\Delta_0$, and
$\delta_i\sim \Delta_{i-1}/\delta_{i-1}$ for all $i>1$. Similar considerations hold for a
sequence which is bounded on the other end. Therefore, we will often study properties of
the polynomials $\Delta_i$ in an exact sequence by studying the $\delta_i$ instead. We will
refer to the $\delta_i$ as the \emph{subpolynomials} of the sequence and to the process of
determining the subpolynomials from the polynomials as ``dividing in from the outside of
the sequence''.

We can also use this polynomial algebra to say something about the relationship between the prime
factors of the polynomial of a module and the prime factors of the polynomials of its submodules
and quotient modules:
 Suppose that $M$ is a torsion $\Gamma$-module with submodule
$N$. Associated to the short exact sequence
\begin{equation*}
\begin{CD}
0 @>>> N @>>> M @>>> M/N @>>> 0,
\end{CD}
\end{equation*}
we have a short exact polynomial sequence of the form
\begin{equation*}
\begin{CD}
0 @>>> f  @>>> h @>>> g @>>> 0,
\end{CD}
\end{equation*}
where $f, g,h\in \Gamma$, $f$ is the polynomial assocatied to $N$, $h$ is the
polynomial associated to $M$, and
$g$ is the polynomial associated to $M/N$. Further, from the properties of
exact
polynomial sequences, we know that we must have $h=fg$.
It is immediate, therefore, that if a prime
$\gamma\in \Gamma$ divides $f$ or $g$ then it divides $h$. Conversely, if
it divides
$h$ then it must divide $f$ or $g$.
We can then drawing the following conclusion: 
\begin{lemma} Suppose that $A$ is a
subquotient of the torsion $\Gamma$-module $M$ (i.e. a quotient module of a submodule of
$M$). Then a
prime $\gamma\in \Gamma$ can divide the polynomial associated to $A$ only if
it divides the polynomial associated to $M$.
\end{lemma}
\begin{proof}
For suppose $A=N/P$, where
$P\subset N\subset M$. If $\gamma$ divides the polynomial of $A$, then by the
above arguments it must divide the polynomial of $N$. But then similarly
$\gamma$ must divide the polynomial of $M$.
\end{proof}

In Section \ref{S: SS,primes}, we will use this fact to study the polynomial algebra of spectral
sequences.

\section{Intersection homology notations and conventions}\label{S: IH con}

We now fix our notation and conventions for dealing with intersection homology.

Let us recall the definition of a stratified pair of paracompact Hausdorff spaces
$(Y,X)$ as given in \cite{CS}. Let $c(Z)$ denote the open cone on the space $Z$,
and let $c(\emptyset)$ be a point. Then the stratification of $(Y,X)$ is defined
by a filtration
\begin{equation*}
Y=Y_n\supset Y_{n-1} \supset Y_{n-2}\supset \cdots \supset Y_0\supset Y_{-1}=\emptyset
\end{equation*}
such that for each point $y\in Y_i-Y_{i-1}$, there exists a distinguished neighborhood 
$N$, a compact Hausdorff pair $(G,F)$, a filtration 
\begin{equation*}
G=G_{n-i-1}\supset  \cdots \supset G_0\supset G_{-1}=\emptyset,
\end{equation*}
and a homeomorphism
\begin{equation*}
\phi: \R^i\times c(G,F)\to (N,N\cap X)
\end{equation*}
that takes $\R^i\times c(G_{j-1},G_{j-1}\cap F)$ onto $(Y_{i+j},Y_{i+j}\cap X)$.  For $(Y,X)$ a
compact PL pair, such a stratification exists with each $\phi$ a PL map and the filtration refining
the filtration by $k$-skeletons (see \cite{Bo}).

Now suppose that we have a PL knot $K$, i.e. a PL embedding $K: S^{n-2}\into S^n$ (though by the standard abuse of notation we also sometimes use $K$ to stand for the image of the embedding). In this case, we have
$Y=S^n$,  and $K$ represents the PL knotted subspace  $S^{n-2}$, which we take as
$Y_{n-1}=Y_{n-2}$ in the filtration. Thus $K$ forms the subspace which is usually
referred to as the singular locus. We will use $\Sigma_k$ to represent the lower
dimensional subspaces of the filtration. In particular, $\Sigma_{n-4}$, which we will
often abbreviate as simply $\Sigma$, will contain all of the points at which the
embedding of $K$ is not locally-flat. (Note that our notation differs from the usual
use of $\Sigma$ in this context.) Thus our stratification has the form
\begin{equation*} 
S^n\supset K\supset \Sigma_{n-4}\supset \Sigma_{n-3}\supset \cdots
\supset \Sigma_0\supset \Sigma_{-1}=\emptyset.
\end{equation*} 
Of course, we can have
$\Sigma_k=\Sigma_{k-1}$. Recall that by \cite{GM2} and \cite{CS}, the intersection
homology modules with a given coefficient system on $S^n-K$ are independent of the choice of
further stratification. We use $IH_i^{\bar p}(S^n;\vg)$ to stand for the intersection homology
modules with perversity $\bar p$ and with coefficient system $\vg$ defined on the complement  $S^n-K$ by stalks $\Gamma$ and action of the fundamental group given by multiplication by $t^{\ell (\alpha, K)}$, where $\ell (\alpha, K)$ is the linking number of $\alpha\in\pi_1(S^n-K)$ with $K$. For simplicity, we will also use $\vg$ to refer to
the restriction of the coefficient system to a subset.

As in \cite{CS}, we will allow
perversities, $\bar p$, with $\bar p(2)$ equal to either $0$ or $1$ (see the following
section for a discussion of the case where $\bar p(2)=1$, which we shall call
\emph{superperversities}). Recall,
however, that a perversity function must be a function from the integers
$\geq
2$ to the non-negative integers satisfying $\bar p(m)\leq \bar p(m+1)\leq
\bar p(m)+1$. We will usually employ the homology notation $IH^{\bar p}_i(S^n;\vg)$
with boundary maps decreasing dimension (as opposed to the cohomology
notation which is also commonly employed). Where sheaves can be be avoided,
we will think of these modules as being defined using finite PL-chains as
in \cite{GM1} or finite singular chains as in \cite{Ki}. As noted in
\cite{Ki}, for geometrically defined intersection homology (see Section \ref{S: superperversities},
below),
these
theories agree with that
of Goresky and MacPherson \cite{GM2} on compact spaces but
not in general unless the sheaf cohomology and hypercohomology 
are taken with compact supports. Hence, our homology theories will always be
those 
with compact supports, though we avoid referring to this in the notation
except
where confusion may arise. Note
also that the versions of the theory developed in \cite{GM1} and \cite{Ki} do
not take advantage of a local coefficient system defined on the complement of
the singular set, but
the definitions there
can be modified to do so easily by taking advantage of the fact, noted in
\cite{GM2} where local coefficients are first introduced to the theory, that
the allowability conditions prevent the simplices of any simplicial or singular chain
and the
simplices of its boundary from lying entirely within the singular set. 

Since it will be used often in the sequel, we state here for convenience the formula for
the intersection homology of a cone. This formula, proved in \cite{Ki}, holds for
geometrically defined intersection homology (e.g. singular intersection homology) for any
perversity. It holds for sheaf intersection homology for traditional perversities ($\bar
p(2)=0$). See the following section (Section \ref{S: superperversities}) for a discussion
of nontraditional perversities and the resulting difference between geometric and sheaf
intersection homology in those cases. Now, suppose that $X$ is a stratified space of
dimension $n-1$
with filtration $\{\Sigma_i\}$ and that the open cone
$c(X)$ is the stratified space filtered by the cone point, $*$, and the collection
$\{c(\Sigma_i)\}$. If $\vg$ is a local coefficient system on $X-\Sigma$, then $\vg\times
\R$ is a local coefficient system on $c(X)-c(\Sigma)$. Conversely, if $\vg$ is a
local coefficient system on $c(X)-c(\Sigma)$, then $\vg\cong \vg|_X\times \R$. Thus, to
simplify notation, we will simply refer to the local coefficient system $\vg$. 
The cone formula is then:
 \begin{equation*}
IH_i^{\bar p} (c(X);\vg)\cong
\begin{cases}
0, & 0\neq i\geq n-1-\bar p(n)\\ 
IH_0^{\bar p} (X;\vg), & i=0 \text{ and } \bar p(n)\geq n-1 \\
IH_i^{\bar p} (X;\vg), & i< n-1-\bar
p(n).
\end{cases}
\end{equation*}
This theorem in \cite{Ki} does not include local coefficients, but the proof goes
through unaltered.

For references to the intersection homology theory, the reader is advised to
consult \cite{GM1}, \cite{GM2}, \cite{Ki}, and \cite{Bo}.

\section{Superperversities}\label{S: superperversities}

Some extra care must be taken when considering the perversities which satisfy
$\bar p(2)=1$. In \cite{CS}, Cappell and Shaneson define the intersection
(co)homology
modules with such perversities via the Deligne sheaf process. For
perversities with $\bar p(2)=0$, the Deligne sheaf complex is quasi-isomorphic to the
complex of sheaves determined geometrically by the presheaf $U\to
IC_{n-i}^{\bar{p},\infty}(U)$, the module of intersection chains with closed
support on $U$ (see \cite{GM2}, \cite{Bo}). This is the route by which one shows that
the intersection (co)homology as defined by the Deligne sheaves coincides with the
intersection (co)homology defined geometrically via allowability conditions on
geometric chains. For $\bar p(2)=1$, however, these theories do \emph{not}
coincide. To see this, recall that the Deligne sheaves, and all sheaves
quasi-isomorphic to them, are required to satisfy a set of axioms determined by the
stratification of the space and the choice of perversity \cite{Bo}. These axioms
are satisfied
by the Deligne sheaves for any perversity, by construction, and by the geometrically
defined intersection chain sheaves for $\bar p(2)=0$. However, the axioms are not
necessarily satisfied by the geometrically defined intersection chains when $\bar
p(2)=1$. 

Consider, for example, the stratified space given by the sphere $S^2$ with singular
locus consisting of a single point, $x$. Let $\mc{IC}^{\infty}_{2-\bullet}$ denote the sheaf
complex
of intersection chains with $\Z$ coefficients and perversity satisfying $\bar
p(2)=1$ (note that we employ the codimension as
index so that the
differentials will be maps of ascending dimension in keeping with the general practices of
sheaf theory; see \cite[\S II]{Bo}). Note that the proof that these sheaves are soft
in \cite[\S II.5]{Bo} holds for any perversity. Now, if these sheaves were to satisfy
the axioms, then, for $j\leq \bar p(2)=1$, there would be an isomorphism
$\mc{H}^j(\mc{IC}^{\infty}_{2-\cdot})_x\cong \mc
H^j(i_*(\mc{IC}^{\infty}_{2-\cdot}|_{S^n-x}))_x$, where $i:S^2-x\to S^2$ is the inclusion map
and we have used the fact that the $\mc{IC}$ are soft to replace the functor $Ri_*$ with $i_*$ in the usual adjunction axiom. We compute 
\begin{align*}
\mc{H}^1(\mc{IC}_{2-\bullet})_x &=\dlim_{x\in U}
H^1(\mc{IC}^{\infty}_{2-\bullet}(U))\\
&=\dlim_{x\in U} H_1(IC^{\infty}_{\bullet}(U)).
\end{align*}
We can take each $U$ to be a distinguished neighborhood of $x$, i.e. a disk about
$x$. Then this group is $0$ because all $2$-chains and $1$-cycles are
allowable, and
we know that the closed-support (Borel-Moore) homology of a disk
$H_1^{\infty}(D^2)=0$. On the other hand:

\begin{align*}
 \mc H^1(i_*(\mc{IC}_{2-\bullet}|_{S^n-x}))_x &= \dlim_{x\in U}
H^1(i_*(\mc{IC}_{2-\bullet}|_{S^n-x})(U))\\
&= \dlim_{x\in U} H^1(\mc{IC}_{2-\bullet}|_{S^n-x}(U-x))\\
&= \dlim_{x\in U} H_1^{\infty}(U-x),
\end{align*}
the last equality because the sheaf of intersection chains restricted to $S^2-x$
is simply the sheaf of closed chains on $S^2-x$, as can be verified locally, 
because
there are no allowability conditions on $S^2-x$. Once again taking each $U\cong D^2$,
the unit disk, 
each $H_1^{\infty}(U-x)$ is non-trivial, generated by the cycle whose
support is given in polar coordinates  by $(r,0)$, $0<r<1$, and which does not
bound. In fact, the directed sequence is constant generated by restrictions of this cycle. Hence, the necessary condition of the axiom can not be satisfied. Similar
counter-examples can be constructed for higher dimensions and for local
coefficients.

The heart of the difficulty is the following: the fact that
the geometric
intersection chain sheaf satisfies the axioms is generally proven by exhibiting,
in a
certain range of dimensions,
isomorphisms between the intersection homology of a  distinguished
neighborhood of a point on a given stratum and the intersection homology of its link,
modulo certain dimension shifts (see \cite[\S II]{Bo} for the
details). In particular, it is shown that
$IH^{\infty}_i(\R^{n-k}\times cL)$ equals
$IH_{i-(n-k+1)}(L)$ for $i\geq n-\bar p(k)$ and equals $0$ otherwise. This is
proven
by constructing a map $\tau_{\geq n-\bar
p(k)}IC^{\infty}_{\bullet-(n-k+1)}(L)\to IC^{\infty}_{\bullet}(\R^{n-k}\times cL)$ and examining
the induced map on homology. This chain map is the composition of a ``cone map''
followed by several ``suspension maps'', where the cone map 
\begin{equation}\label{E: cone map}
c:\tau_{\geq
k-\bar p(k)}IC^{\infty}_{\bullet -1}(L)\to IC^{\infty}_{\bullet}(cL)
\end{equation} 
is induced simply by taking the
open cones on allowable chains in the link. Simple dimensional and allowability
arguments show that for perversities satisfying $\bar p(2)=0$, the cone on such an
allowable chain, $\xi_{i-1}\in IC_{i-1}(L)$, is itself allowable in the following
situations:
\begin{align*}
\begin{cases}
i>k-\bar p(k),& \text{any $\xi$}\\
i=k-\bar p(k), & \text{only if $\xi$ is a cycle}\\
i<k-\bar p(k), & \text{no $\xi$}.
\end{cases}
\end{align*}
Hence, the cone map is well-defined, one goes on to show that it is a quasi-isomorphism,
and the arguments proceed. For $\bar p(2)=1$, however, we are stopped dead
in our tracks. These perversities allow the possibility of $i=k-\bar p(k)=1$,
but an allowable  vertex in the link does \emph{not} cone to an allowable chain in
$cL$ because the cone on a vertex will have a boundary vertex in the singular
locus (note that this problem does not arise for the traditional perversities
with $\bar p(2)=0$ because the above dimension and allowability conditions already
forbid us from coning a vertex). Thus the cone map \eqref{E: cone map} does not
necessarily exist. Of
course, one could devise a new cone map by changing the range of the truncation, but
following the program through with this new cone map will not ultimately yield a set
of axioms equivalent to those satisfied by the Deligne sheaves with these
perversities. 

To keep matters straight, we shall adopt the following definitions. We will use
``perversity'' or ``traditional perversity'' to refer to the perversities with
$\bar p(2)=0$, and we will refer to those satisfying $\bar p(2)=1$ as
\emph{superperversities}. All intersection homology
groups with superperversities will be those as defined by the Deligne sheaf.
Hence, in these cases, we take $IH_i=IH^{n-i}$ for a space of dimension $n$, where
the latter are hypercohomology modules of the Deligne sheaf (recall that we
assume compact supports unless otherwise noted).  In the case where $\bar p(2)=1$, it is not even
clear how to define a geometrical (simplicial or singular) version of intersection homology with local
coefficients (though of course there is a perfectly well-defined theory for constant
coefficients; see, e.g, \cite{Ki}, \cite{Ki2}). 
Note that even for superperversities, sheaf intersection
homology is a topological invariant by a straightforward adaption of the proofs in \cite{GM2} and
\cite{Bo} (or even, in the cases with which we will be concerned, by the superduality of Cappell and Shaneson (see
below)). \emph{Unless specified otherwise, we will assume below that all
perversities are traditional perversities with $\bar p(2)=0$.} This is not a
large restriction, as we may employ Theorem \ref{T: superduality} below to
calculate the superperverse intersection homology modules from those with
traditional perversities.  Also, recall again that simiplicial, singular, and sheaf intersection
homology with compact supports all  agree for traditional perversities,
and in this case, we are free to make use of the singular geometric
theory of King \cite{Ki}.

\section{Duality properties of intersection Alexander polynomials}\label{S: IH dual}

We first  prove a general theorem concerning the duality of intersection
Alexander polynomials which is analogous to that for the usual Alexander
polynomials: Given a locally-flat knot $S^{n-2}\subset S^n$, let $\lambda_i$ be the polynomial associated to the homology module $H_i(S^n-K;\vg)$, where $\vg$ is defined as above. Then for $0<i<n-1$, $\lambda_i(t)\sim \lambda_{n-i-1}(t^{-1})$, where $\sim$ indicated similarity up to associates in $\Gamma$ (\cite{L66}).

For non-locally flat knots, we let $I\lambda^{\bar p}_i(t)$ represent the analogous
intersection Alexander
polynomials, which are defined as the polynomials of the modules $IH_i^{\bar p}
(S^n;\vg)$. (That these intersection homology modules are torsion modules, and hence
have well-defined associated polynomials, is a
consequence of \cite[Proposition 2.4]{CS}.)

\begin{theorem}\label{T: superduality}
Let $K\cong S^{n-2}\subset S^n$ be a knot, not necessarily locally-flat, and let $\bar p$
and $\bar q$ be a superdual perversity and superperversity as defined in \cite{CS}, i.e.  $\bar p(k)+\bar
q(k)=k-1$ for all $k\geq 2$.Then $I\lambda^{\bar p}_i(t)\sim I\lambda^{\bar q}_{n-1-i}(t^{-1})$.
\end{theorem}
\begin{proof}
In \cite{CS}, Cappell and Shaneson define a sub-pseudomanifold $X$  of a sphere $S^n$ to be
of \emph{finite (homological) type} if $H_i(S^{n+2}-X;\Upsilon)$ is finite dimensional as
a $\Q$-vector space. They define a sub-pseudomanifold $X$ of a manifold $Y$ to be of
\emph{finite
local type} if the link of each component of any stratification is of finite type. (See
\cite[Section 1]{CS} for a detailed description of the local coefficient system
$\Upsilon$ (there denoted $\Lambda$, but changed here to avoid conflict with our other
$\Lambda=\Z[\Z]$)). For our purposes, $\Upsilon$ will always be the coefficient system $\vg$,
as above, and in this case, by \cite[Proposition 2.2]{CS}, the knot $K$ has finite type and finite local type.
Thus by \cite[Corollary 3.4]{CS}, for the knot $K$ and superdual perversities $\bar p$
and $\bar q$, 
\begin{equation*}
\overline{IH_i^{\bar p}(S^n;\vg)}\cong \text{Hom}(IH_{n-1-i}^{\bar q}(S^n;\vg),
\Q(t)/\Gamma),
\end{equation*}
where $\Q(t)$ is the field of fractions of $\Gamma$, i.e. the field of rational
functions,
and $\bar{A}$ is the
$\Gamma$-module obtained from the $\Gamma$-module $A$ by composing all module structures
with the involution $p(t)\to p(t^{-1})$. 

We claim that $\text{Hom}(\Gamma/(p), \Q(t)/\Gamma)\cong \Gamma/(p)$. 
Since $\Gamma$ is an integral domain, we have $\text{Hom}(\Gamma/(p),
\Q(t)/\Gamma)\cong$ $\text{Ext}(\Gamma/(p),\Gamma)$ by \cite[Proposition VII,2.3]{CE}. Then, from the short exact sequence (and
free resolution)
\begin{equation*}
\begin{CD}
0@>>> \Gamma @>p>> \Gamma @>>> \Gamma/(p) @>>> 0,
\end{CD}
\end{equation*}
the long exact sequence of left derived functors of $Hom(-,\Gamma)$ gives
\begin{equation*}
0\la \text{Ext}(\Gamma/(p),\Gamma)\la
\text{Hom}(\Gamma,\Gamma) \overset{p}{\la}\text{Hom}(\Gamma,\Gamma) \la 
\text{Hom}(\Gamma/(p),\Gamma) \la 0,
\end{equation*}
where the leftmost term is $0$ because Ext$(\Gamma,\Gamma)=0$, $\Gamma$ being a free
$\Gamma$-module. Since $\text{Hom}(\Gamma/(p),\Gamma)=0$, because any map from a
torsion module to a free module is $0$, and since Hom$(\Gamma,\Gamma)\cong \Gamma$, this
becomes the short exact sequence
\begin{equation*}
\begin{CD}
0@<<< \text{Ext}(\Gamma/(p),\Gamma)@<<< 
\Gamma @<p<< \Gamma  @<<< 0. 
\end{CD}
\end{equation*}
Therefore,  Hom$(\Gamma/(p), \Q(t)/\Gamma)\cong$Ext$(\Gamma/(p),\Gamma)\cong \Gamma/(p)$. 

By  \cite[Proposition 2.4]{CS}, the intersection homology modules $IH_i^{\bar
p}(S^n;\vg)$ are finite dimensional as $\Q$-vector spaces because $K$ has finite type and
finite local type. Therefore, they must be finitely generated torsion
$\Gamma$-modules. Thus, $IH_{n-1-i}^{\bar q}(S^n;\vg)\cong
\oplus_j \Gamma/(p^{\bar q}_{n-1-i,j})$, and $I\lambda_{n-1-i}^{\bar q}\cong \prod_j
p^{\bar q}_{n-1-i,j}$. By the results of the preceding paragraphs and the fact that
Hom$(\oplus A_j, B)\cong\oplus$Hom$(A_j,B)$ for finite direct sums, 
\begin{align*}
\overline{IH_i^{\bar p}(S^n;\vg)}&\cong \text{Hom}(IH_{n-1-i}^{\bar
q}(S^n;\vg),
\Q(t)/\Gamma)\\
&\cong IH_{n-1-i}^{\bar q}(S^n;\vg)\\
&\cong \oplus_j
\Gamma/(p^{\bar
q}_{n-1-i,j}).
\end{align*}
Thus 
\begin{equation*}
IH_i^{\bar p}(S^n;\vg)\cong \oplus_j \Gamma/(p^{\bar q}_{n-1-i,j}(t^{-1})),  
\end{equation*}
and $I\lambda_i^{\bar p} (t)\sim \prod_j p^{\bar
q}_{n-1-i,j}(t^{-1})=I\lambda_{n-1-i}^{\bar q} (t^{-1})$.
\end{proof}

\begin{corollary}
For a not necessarily locally-flat knot $S^{n-2}\subset S^n$, $\lambda^{\bar p}_{i}\sim 1$ for
$i\geq n$.
\end{corollary}
\begin{proof}
This is a consequence of the theorem and the fact that $IH_i^{\bar p}(S^n;\vg)=0$ for all $i<0$.
\end{proof}

This is also a convenient place to point out the following:
\begin{corollary}
For a not necessarily locally-flat knot $S^{n-2}\subset S^n$, $I\lambda^{\bar p}_{0}\sim t-1$ for a
traditional perversity $\bar p$.
\end{corollary}
\begin{proof}
Because $\bar p(2)=0$, we can employ the geometric theory, and the allowable zero- and one-chains are those which lie in $S^n-K$. Therefore, 
\begin{equation*}
IH_0^{\bar p}(S^n;\vg)\cong H_0(S^n-K;\vg)\cong \Gamma/(t-1),
\end{equation*}
the second isomorphism being clear from the identification of the homology of $S^n-K$
with local coefficient system $\vg$ with the rational homology (viewed as a
$\Gamma$-module) of the infinite cyclic
cover of $S^n-K$ (see \cite[\S 4.3.2]{GBF1} or \cite{GBF}).

\end{proof}

\section[Normalization properties of $I\protect\lambda^{\bar
p}_i$]{Normalization properties of {\mv $I\protect\lambda^{\bar
p}_i$}}\label{S: IH
norm}

Suppose that $p$ is an element of $\Gamma$. Recall that there is an
element, say $p'$, in the similarity class of $\Gamma$ which is primitive
in
$\Lambda=\Z[t,t^{-1}]$, i.e. the coefficients are relatively prime
(though not necessarily pairwise so), and this element is unique up to
similarity class in $\Lambda$ (see, e.g., \cite{L66} or \cite{GBF1}). We will say that $p$ is a polynomial of
Alexander type if $p'(1)=\pm 1$. Another classical property of the Alexander polynomials of locally-flat knots $S^{n-2}\subset S^n$ is that they are of Alexander type for $0<i<n-1$ (\cite{L66}). In fact, this is also true of the Alexander polynomials of knots which are not locally flat (\cite{GBF1}, \cite{GBF}), where in this case the polynomials are again defined to be those associated to the modules $H_i(S^n-K;\vg)$. We wish to show that the same property holds for intersection Alexander polynomials.

Observe, first of all,  that if the polynomial associated to a torsion $\Gamma$-module,
$M$,
is of Alexander type, then so are the polynomials associated to any
submodule or
quotient module of $M$, as follows immediately from the short exact
polynomial
sequence (see Section \ref{S: poly alg}) associated to the short exact module sequence which
represents
the
inclusion or quotient and from the fact that an integer polynomial factors over
$\Q$   
if and only if it factors over $\Z$. Of course, if a primitive polynomial
in
$\Lambda$ evaluates to $\pm 1$ when $t=1$, then this will be true of any
of
its factors in $\Lambda$, which will also be primitive. Thus, any factor of
a
polynomial of Alexander type is also of Alexander type. Note, in particular, that the element $1\in
\Gamma$ is of Alexander type.

\begin{theorem}\label{T: normal}
For any PL-knot $K\cong S^{n-2}\subset S^n$, not necessarily locally-flat, and traditional
perversity $\bar p$,
$I\lambda_i^{\bar p}$ is of Alexander type for $i>0$, $I\lambda_0^{\bar p}\sim t-1$,
and $I\lambda_i^{\bar p}\sim 1$ for $0\neq i\geq n-1$.
\end{theorem}
\begin{proof}
We will proceed by induction on the dimension $n$. We begin with a trivial
low-dimensional case. For $n=1$, we define the knot by the pair $(S^1,\emptyset)$ and
define the local coefficient system $\vg$ so that the generator of $\pi_1(S^1)$ acts
on
the stalk $\Gamma$ by multiplication by $t$. (This choice of $\vg$ is made to be consistent with the coefficient system that will appear on $S^1$ when it is considered as  the link of the top stratum for knots of higher dimension.) In this case it is clear that
\begin{equation*}
IH_i^{\bar p}(S^1;\vg)\cong H_i(S^1;\vg)\cong
\begin{cases}
\Gamma/(t-1),& i=0\\
0,& i>0.
\end{cases}
\end{equation*}
Since the polynomial associated to the zero module is $1$, up to similarity, the
theorem holds for $n=1$.

From this case, we will proceed by induction on $n$, so let us assume that the theorem
holds for knots in $S^{n-1}$ and show that it holds for knots in $S^n$. This
step will
also proceed by induction, this time on the codimension of the strata. First, however,
we will choose a useful stratification, as we are free to do by the topological
invariance of intersection homology. In fact, we will merely refine the ``natural''
stratification defined by the embedding, so there will be no difficulty with the definition
of the local coefficient system $\vg$. We continue to allow $K$, itself, to determine
the singular locus, $\Sigma_{n-2}$, but we now define $\Sigma_j$, $j<n-2$, to be the
$j$-skeleton of $K\subset S^n$ for some triangulation of $S^n$ for which $K$ is a
full subpolyhedron.

Let $U_j=S^n-\Sigma_{n-j}$, $2\leq j \leq n+1$. Thus $U_2\cong S^n-K$, $U_3\cong
S^n-K\cup\{$the open $n-2$-simplices of $K\}$, and so on up to $U_{n+1}\cong S^n$. The
induction will be over this codimension $j$. In other words, for each $j$ we show
$IH_i^{\bar p}(U_j;\vg)$ has associated polynomials which are of Alexander type
for $0<i<n-1$, trivial ($\sim 1$) for $i\geq n-1$, and similar to $t-1$ for $i=0$. The theorem
will then be proven for knots in $S^n$ once we have inducted up to $j=n+1$.

For $j=2$, $U_2\cong S^n-K$ and $IH_i^{\bar p}(S^n-K;\vg)\cong H_i(S^n-K;\vg)$. In this
case, the desired conclusions hold because they are true for the ordinary
(singular) homology modules of the knot complement (\cite[Thm. 4.3]{GBF1}).

Suppose now that the claim holds for $U_{j-1}$, $j\geq 3$. We will show that it is true
for $U_j$. Note that $U_j-U_{j-1}$ is the union of the (finite number of) open simplices
of $K$ of dimension $n-j+1$, say $\{e^{\alpha}_{n-j+1}\}$.  For each $
e^{\alpha}_{n-j+1}$, consider its neighborhood in $U_j$ defined by taking, in a derived
subdivision of the triangulation of $S^n$, the union of the open simplices whose closures
intersects $ e^{\alpha}_{n-j+1}$. Note that, because $ e^{\alpha}_{n-j+1}$ is the open
simplex, we do not include in the neighborhood those open simplices whose closures only
intersect $\bd \bar e^{\alpha}_{n-j+1}$. Let us call these neighborhoods
$N^{\alpha}_{n-j+1}$. Then $N=\cup_{\alpha} N^{\alpha}_{n-j+1}$ is a neighborhood
of $U_j-U_{j-1}$ in $U_j$, and $ N^{\alpha}_{n-j+1}\cap N^{\beta}_{n-j+1}=\emptyset$ for
$\alpha\neq \beta$. Furthermore, each $ N^{\alpha}_{n-j+1}$ is homeomorphic to
$\R^{n-j+1}\times cL^{\alpha}$, where $cL^{\alpha}$ is the open cone on the link determined by
the stratification (hence $L^{\alpha}\cong S^{j-2}$, but the coefficient system on the
link is determined by $\vg|L^{\alpha}$, which may depend on $\alpha$).  Note that $
N^{\alpha}_{n-j+1}-e^{\alpha}_{n-j+1}\cong \R^{n-j+1}\times (cL^{\alpha}-*)$ where $*$
represents the cone point.

With $N$ as above, let $N'=N\cap U_{j-1}=N-\cup_{\alpha}
e^{\alpha}_{n-j+1}$, and consider the Mayer-Vietoris sequence
\begin{equation}\label{E: MVU}
\to IH_i^{\bar p}(N';\vg) \to IH_i^{\bar p}(N;\vg)\oplus IH_i^{\bar
p}(U_{j-1};\vg) \to IH_i^{\bar p}(U_j;\vg)\to.
\end{equation}

Firstly, for $i=0$, $IH_0^{\bar p}(U_j;\vg)\cong H_0(U_j-K\cap U_j;\vg)\cong
\Gamma/(t-1)$,
the first isomorphism because we assume $\bar p(2)=0$ so that the $0$- and
$1$-intersection chains must lie outside of $K$ and the second because $U_j$ is
connected and the homology of the complement with coefficients in $\vg$ is equal to the
rational homology of the infinite cyclic cover as a $\Gamma$-module with trivial $t$ action, i.e. $\Q\cong \Gamma/(t-1)$.
Furthermore, the same holds for the zero degree intersection homology of each component
of $N'$ and ${N}$ with the inclusion of the generating points inducing
an isomorphism $ IH_0^{\bar p}(N';\vg) \overset{\cong}{\to} IH_0^{\bar
p}(N;\vg)$. In particular, the corresponding map of the Mayer-Vietoris sequence is
injective.

Next, we consider $i>0$. By the induction step, the polynomial associated to $IH_i^{\bar
p}(U_{j-1};\vg)$ is of Alexander type, $i>0$, and it is similar to $1$ for $0\neq i\geq n-1$.  By
the
K\"unneth theorem \cite{Ki},
\begin{align*}
IH_i^{\bar p}(N;\vg) &\cong \oplus_{\alpha} IH_i^{\bar p}(cL^{\alpha};\vg)\\
IH_i^{\bar p}(N;\vg) &\cong \oplus_{\alpha} IH_i^{\bar p}(cL^{\alpha}-*;\vg)\cong \oplus_{\alpha}IH_i^{\bar p}(L^{\alpha};\vg).
\end{align*}
But $ IH_i^{\bar p}(L^{\alpha};\vg)$ is the intersection homology of the link knot pair
which, by induction (as the dimension of the link sphere is $<n$), has associated
polynomial of Alexander type. Applying the formula for the intersection homology of a
cone (see \cite{Ki}), according to which $IH_i^{\bar p}(cL)$ is equal to either  $IH_i^{\bar
p}(L)$ or $0$, the polynomial associated to each $IH_i^{\bar
p}(cL^{\alpha};\vg)$ is of Alexander type. Since the product of polynomials of Alexander type is of
Alexander
type, the same is true of the polynomials associated to $ IH_i^{\bar p}(N;\vg)$ and
$ IH_i^{\bar p}(N';\vg)$. Examining the long exact polynomial sequence (see Section \ref{S: poly 
alg}) associated to the Mayer-Vietoris sequence \eqref{E: MVU}, we can conclude that the
polynomials associated to the $ IH_i^{\bar p}(U_j;\vg) $ are of Alexander type as each
is
the product of factors of the preceding and following terms, and factors of
polynomials of Alexander type are of Alexander type. (Note that we use here the
injectivity results of the last paragraph to see that the $t-1$ factors of the polynomial of
$IH_0^{\bar
p}(N';\vg)$ are not shared with the polynomial associated with $
IH_1^{\bar p}(U_j;\vg))$.

Finally, by the induction on the dimension of the knot, $ IH_i^{\bar
p}(N';\vg) \cong \oplus_{\alpha}IH_i^{\bar p}(L^{\alpha};\vg)$ has
polynomial similar to $1$ for $0\neq i\geq n-2$ (since the link is a sphere of dimension
$<n$), and
similarly for $IH_i^{\bar p}(N;\vg)\cong \oplus_{\alpha} IH_i^{\bar p}(cL^{\alpha};\vg)$
by again applying the formula for intersection homology of a cone. The same holds for $IH_i^{\bar
p}(U_{j-1};\vg)$, $0\neq i\geq n-1$, by the induction on the codimension $j$. Hence, the
polynomial associated to $IH_i^{\bar p}(U_{j};\vg)$ for $0\neq i\geq n-1$ is $1$.

This completes the induction step, and the proof follows.
\end{proof}

\begin{corollary}\label{C: IH high dim}
Suppose that $\bar p$ is a superperversity, and let $I\lambda_i^{\bar p}$ be the
polynomial associated to the sheaf intersection homology module $IH_i^{\bar
p}(S^n;\vg)$ of the knot $K\subset S^n$, $n\geq 2$. Then $I\lambda_i^{\bar p}$ is of Alexander
type for $0<i<n-1$, is similar to $1$ for $i=0$ or $i>n-1$, and is similar to $t-1$ for
$i=n-1$.
\end{corollary}
\begin{proof}
This follows immediately from the theorem and superduality (see Theorem \ref{T: superduality}).
\end{proof}

\section{Some relations with the ordinary Alexander
polynomials}\label{S: ih and ord}

We now turn to calculating intersection Alexander polynomials in some special cases in terms of other polynomial
invariants. We will consider a number of cases,
each more general than the preceding. We could start at the end and deduce some of 
the earlier
conclusions as special cases of the later ones, but it is more instructive to show the
development in order of increasing complexity. Recall that we assume all perversities satisfy $\bar p(2)=0$ unless otherwise specified.
 
\subsection{Intersection homology of locally-flat knots}

We will begin by showing that in the case of a locally-flat knot $K$,
$IH_{i}^{\bar p}(S^n;\vg)\cong H_i(S^n-K;\vg)$, the usual Alexander module of the knot.

\begin{proposition}\label{P: locally-flat}
Let $K$ be a locally-flat PL-knot $S^{n-2}\subset S^n$. Then, with the
notation as
above, $IH_i^{\bar p}(S^n;\vg)=H_i(S^n-K;\vg)$.
\end{proposition}
\begin{proof}
Let $N(K)$ be an open regular neighborhood of $K$, let $N'(K)$ be a ``smaller'' open
regular neighborhood of $K$ with $\overline{N'(K)}\subset N(K)$, and let $X(K)$ be the knot
exterior given by the open subset $S^n-\overline{N'(K)}$. Then, using the ``generalized
annulus theorem'' (see \cite[Proposition 1.5]{Sto}), $X(K)\cap N(K)\cong \bd
\overline{N(K)} \times \R$. 

We have the Mayer-Vietoris sequence
\begin{equation*}
\to IH_i^{\bar p} (X(K)\cap N(K);\vg) \overset{i_*}{\to}IH_i^{\bar p}
(X(K);\vg)\oplus
IH_i^{\bar p} (N(K);\vg) \to IH_i^{\bar p} (S^n;\vg) \to
\end{equation*}
because, as noted in \cite[p. \S 2]{Ki}, the standard singular homology proof carries
through for intersection homology defined by singular intersection chains. This
proof is unaltered for local coefficients. We will use this sequence to compute $
IH_i^{\bar p} (S^n;\vg)$.

We begin by observing that $IH_i^{\bar p} (X(K);\vg)\cong H_i(X(K);\vg)\cong
H_i(S^n-K;\vg)$ because the singular set $K$ does not intersect $X(K)$, and it is well
known that in this case intersection homology agrees with the standard homology
theories (recall that we are assuming compactly supported homology at all times). The
last isomorphism then follows from the homotopy equivalence of $X(K)$ and $S^n-K$.
Hence this term is the usual Alexander module of the knot.

Similarly, 
\begin{align*}
IH_i^{\bar p} (X(K)\cap N(K);\vg)&\cong H_i(X(K)\cap N(K);\vg)\\
&\cong H_i(\bd
\overline{N(K)} \times \R;\vg)\\
&\cong H_i(\bd \overline{N(K)};\vg). 
\end{align*}
From the standard obstruction theory argument in the proof of the existence of Seifert
surfaces (see\cite{L66}), we know that there is a trivialization $\bd \overline{N(K)}\to S^1$ so that
$H_i(\bd \overline{N(K)};\vg)\cong H_i(S^{n-2}\times S^1;\vg)$. This can be computed
using
the K\"{u}nneth theorem. Since $S^1$ represents a meridian, an easy computation gives
\begin{equation*}
H_i(S^1;\vg)\cong
\begin{cases}
\Gamma/(t-1)\cong \Q, & i=0\\
0,& i\neq 0.
\end{cases}
\end{equation*}
(Whenever we write $\Q$, we will mean it as a $\Gamma$-module with trivial $t$ action, i.e.  $\Q\cong \Gamma/(t-1)$.)
If $n\geq 4$, $S^{n-2}$ is simply-connected so that, if $\pi_1$ and $\pi_2$ are the
projections, the local coefficient system on
$S^{n-2}\times S^1$ is $\pi_1^{*}\Gamma\otimes \pi_2^{*}(\vg|{S^1})$, and
\begin{equation*}
H_i(S^{n-2};\Gamma)\cong
\begin{cases}
\Gamma, & i=0, n-2\\
0,& i\neq 0,n-2.
\end{cases}
\end{equation*}
If $n=3$, we can obtain the same equations by choosing for the first factor $S^{n-2}=S^1$ a
``preferred longitude'' (see \cite{Rolf}). 
Thus by the K\"{u}nneth theorem, we have
\begin{equation*}
IH_i^{\bar p} (X(K)\cap N(K);\vg)\cong
\begin{cases}
\Gamma\otimes_{\Gamma} \Q\cong \Q, & i=0, n-2\\
0,& i\neq 0, n-2.
\end{cases}
\end{equation*}

Lastly, we consider the term $IH_i^{\bar p} (N(K);\vg)$. Since the knot is locally
flat, $N(K)$ is homeomorphic to a disk bundle over the knot, and, by extending the trivialization of the
boundary $\bd \overline{ N(K)}$ to the interior of each disk, we have that $ N(K)\cong
S^{n-2}\times D^2$. Since $S^{n-2}$ is an unstratified  manifold, the K\"{u}nneth theorem applies as
proved in \cite[\S 2]{Ki}, since the argument there easily extends to include local
coefficient systems. We need to compute the terms of the K\"{u}nneth formula. Here,
$D^2\cong c(S^1)$, where $S^1$ represents a meridian of the knot and the cone point is
the singular set. In this case, we can use the formula for the intersection homology
of a cone, also in \cite[\S 2]{Ki}. Since $IH_i^{\bar p}(S^1;\vg)\cong H_i(S^1;\vg)$
as in the last paragraph, we can calculate
\begin{equation*}
IH_i^{\bar p} (D^2;\vg)\cong
\begin{cases}
\Gamma/(t-1)\cong \Q, & i=0\\
0,& i\neq 0.
\end{cases}
\end{equation*}
And, just as above, we have
\begin{equation*}
H_i(S^{n-2};\Gamma)\cong
\begin{cases}
\Gamma, & i=0, n-2\\
0,& i\neq 0,n-2.
\end{cases}
\end{equation*}
Therefore,
\begin{equation*}
IH_i^{\bar p} (N(K);\vg)\cong
\begin{cases}
\Gamma\otimes_{\Gamma} \Q\cong \Q, & i=0, n-2\\
0,& i\neq 0, n-2.
\end{cases}
\end{equation*}

From these calculations, we see that the proposition must be true for
$i\neq 0,1,n-1,n-2$. But, for $i=0$, the map $i_*: IH_0^{\bar p} (X(K)\cap
N(K);\vg) \to$ $ IH_0^{\bar p} (X(K);\vg)$ corresponds to the map $H_0(X(K)\cap N(K);\vg) \to H_0(X(K);\vg)$ induced by the inclusion of
a
point. Hence, this map is an injection $\Q\to\Q$ by the usual arguments in
dimension $0$.  Thus, the proposition will be proven if we can show that
$i_*: IH_{n-2}^{\bar p} (X(K)\cap N(K);\vg) \to IH_{n-2}^{\bar p}
(X(K);\vg)$ is an injection. Since this Mayer-Vietoris map is induced by
inclusion, the calculations above and the naturality of the K\"{u}nneth
theorem provide that $i_*$ is equivalent to the map $\text{id}_*\otimes j_*:
H_{n-2}(S^{n-2};\vg)\otimes H_0(S^1;\vg)\to H_{n-2}(S^{n-2};\vg)\otimes
H_0(D^2;\vg)$, where $j$ is the inclusion $j:S^1\to D^2=c(S^1)$. By the
computation of the intersection homology of a cone in \cite[\S 4]{Ki2},
$j_*$ is a surjection, but any surjection $\Q \to \Q$ is also an
isomorphism. Thus $i_*$ is an injection.
\end{proof}

\subsection{Point singularities}
We next consider the case where the knot $K$ has a single point singularity. In other
words, $\Sigma=\Sigma_0=*$. It will be useful to introduce some notation from \cite{GBF} and \cite{GBF1}. Let $D$ be the open regular neighborhood of the singular point of the embedding. Then the complement of $D$ in the pair $(S^n, S^{n-2})$ is a locally flat disk knot bounded by a locally-flat sphere knot, which is the link knot around the singular point. Futhermore, the knot complement $S^n-K$ is homotopy equivalent to the complement of the induced disk knot (see \cite[\S 3]{GBF1} for details). We can then define Alexander polynomials $\lambda_i$, $\nu_i$, and $\mu_i$ associated to the homology modules of the disk knot complement $C$, the boundary sphere knot complement $X$, and the pair $(C,X)$. Note that the  $\lambda_i$ also represent the Alexander polynomials of the sphere knot $K$ and the $\nu_i$ are just the usual Alexander polynomials of the locally-flat link knot. 

Let $\gamma_i\in \Gamma$. Recall that, in Section \ref{S: poly alg}, we defined  an
\emph{exact sequence of polynomials}, denoted by
\begin{equation*}
\begin{CD}
@>>> \gamma_{i-1} @>>>\gamma_i @>>> \gamma_{i+1} @>>>,
\end{CD}
\end{equation*} 
to mean a sequence of polynomials such that each $\gamma_i\sim \delta_i \delta_{i+1}$, $\delta_i\in
\Gamma$. Recall also that such a sequence is determined by an exact sequence of
torsion
$\Gamma$-modules, and, in that case, the factorization of the polynomials is determined by the maps of
the modules (in fact, each $\delta_i$ will be the polynomial of the kernel module of a map in the
exact module sequence).

For knots,  we therefore have the exact sequence
\begin{equation*}
\begin{CD}
@>>> \nu_{i}(t)@>>> \lambda_i(t) @>>> \mu_i(t) @>>> \nu_{i-1}(t) @>>>.
\end{CD}
\end{equation*}
To emphasize the various shared factors, we sometimes rewrite this as 
\begin{equation*}
\to a_{i}(t)b_i(t)  \to b_i(t)c_i(t) \to c_i(t)a_{i-1}(t) \to
a_{i-1}(t)b_{i-1}(t) \to,
\end{equation*}
and we refer to the polynomials $a_i$, $b_i$, and $c_i$ as \emph{Alexander subpolynomials} or just \emph{subpolynomials}.

Recall also that knowledge of  two thirds of the terms of an exact sequence of polynomials
(for example,   
all $\gamma_{3i}$ and $\gamma_{3i+1}$, $i\in\Z$) and the common
factors of those terms (the $\delta_{3i+1}$),  allows us to deduce the missing third of the sequence   
($\Delta_{3i+2}=
\delta_{3i+2}\delta_{3i+3}=\frac{\Delta_{3i+1}}{\delta_{3i+1}}\cdot\frac{\Delta_{3i+3}}{\delta_{3i+4}}$). Therefore, the Alexander polynomials and subpolynomials carry the same information in a sequence that is only finitely non-trivial.

\begin{proposition}\label{P: ih point}
With the notation as above, the intersection Alexander polynomial of a knot, $K$, embedded with a single point singularity is given by 
\begin{equation*}
I\lambda_i^{\bar p}(t) \sim
\begin{cases}
\lambda_i(t), & i<n-1-\bar p(n)\\
c_i(t), & i=n-1-\bar p(n)\\
\mu_i(t), & i> n-1-\bar p(n).
\end{cases}
\end{equation*}
\end{proposition}
\begin{proof}
Let $N(\Sigma)$ be an open regular neighborhood of $\Sigma$, which will be an open
PL-ball, $D^n$. Let $N'(\Sigma)$ be another open regular neighborhood of $\Sigma$
with $\overline{N'(\Sigma)}\subset N(\Sigma)$, and let $X(\Sigma)$ be the knot exterior given by
the open subset $S^n-\overline{N'(\Sigma)}$, which will also be a PL-ball, $D^n$. As
in the
last proposition, we will compute the intersection homology via the Mayer-Vietoris
sequence
\begin{equation*}
\to IH_i^{\bar p} (X(\Sigma)\cap N(\Sigma);\vg) \overset{i_*}{\to}IH_i^{\bar
p}
(X(\Sigma);\vg)\oplus IH_i^{\bar p} (N(\Sigma);\vg)\to IH_i^{\bar p}
(S^n;\vg) \to.
\end{equation*}
 
The subspace $X(\Sigma)$ is an open ball, and the pair $(X(\Sigma), X(\Sigma)\cap K)$ is a
locally-flat (open) disk knot. The arguments of Proposition \ref{P: locally-flat} carry over with only minor
alteration 
to show that $IH_i^{\bar p} (X(\Sigma)-K\cap
X(\Sigma);\vg)\cong H_i(X(\Sigma);\vg)\cong
H_i(S^n-K;\vg)$ 
(see also the proof of Proposition \ref{P: man triv1} below).

If $L(\Sigma)$ is the link of the point $\Sigma$, which will be a sphere $S^{n-1}$,
then $ X(\Sigma)\cap N(\Sigma)\cong L(\Sigma)\times \R$ by the generalized annulus
property. Therefore, using the K\"unneth Theorem for intersection homology (\cite[Lemma 3]{Ki}), $
IH_i^{\bar p} (X(\Sigma)\cap
N(\Sigma);\vg)\cong IH_i^{\bar p} (L(\Sigma)\times \R;\vg)\cong IH_i^{\bar p}
(L(\Sigma);\vg)$. But the pair $( L(\Sigma), L(\Sigma)\cap K)$ is a locally-flat
sphere knot pair, and $\vg|_{L(\Sigma)}$ gives the usual local coefficient system for a
sphere knot because the linking number of $K$ with a loop in $L$ is the same
as if we consider the loop to be in $S^n$. Thus, by Proposition \ref{P: locally-flat},
$IH_i^{\bar p} (X(\Sigma)\cap N(\Sigma);\vg)\cong H_i(L(\Sigma)-K\cap L(\Sigma);\vg)$,
the Alexander module of the link knot.

Lastly, $N(\Sigma)\cong c(L(\Sigma))$, and we can use \cite[Proposition 5]{Ki}  to compute
{\small\begin{equation*}
IH_i^{\bar p} (cL(K);\vg)\cong
\begin{cases}
0, & 0\neq i\geq n-1-\bar p(n)\\
IH_0^{\bar p} (L(\Sigma);\vg)\cong \Gamma/(t-1), & i=0, \bar
p(n)\geq n-1 \\
IH_i^{\bar p} (L(\Sigma);\vg)\cong H_i (L(\Sigma)-K;\vg)
, & i<
n-1-\bar p(n).\\
\end{cases}
\end{equation*}}

Since the map of the Mayer-Vietoris sequence $i_*: H_i^{\bar p} (X(\Sigma)\cap
N(\Sigma);\vg) \to IH_i^{\bar p} (N(\Sigma);\vg)$ is induced by inclusion and the
above
calculations show that this is equivalent to the map induced by the inclusion $i_*:
H_i^{\bar p} (L(\Sigma);\vg) \to IH_i^{\bar p} (cL(\Sigma);\vg)$, Proposition 5 of
\cite{Ki} allows us to conclude that $i_*$ is an isomorphism for $i=0$ or $i<n-1-\bar
p(n)$.
Thus for $i<n-1-\bar p(n)$, $i_*: IH_i^{\bar p} (X(\Sigma)\cap N(\Sigma);\vg) \to
H_i^{\bar p} (X(\Sigma);\vg)\oplus IH_i^{\bar p} (N(\Sigma);\vg)$ is an injection.
Using this fact and the calculations of the preceding paragraphs, we can break the
Mayer-Vietoris sequence in this range into the short exact sequences
\begin{equation*}
0\to H_i (L(\Sigma)- K;\vg) \overset{i_*}{\to} H_i (L(\Sigma)-K;\vg) \oplus
H_i (S^n-K;\vg) \to IH_i^{\bar p} (S^n;\vg) \to 0.
\end{equation*}
Therefore, $IH_i^{\bar p} (S^n;\vg)\cong H_i (S^n-K;\vg)$ and
$I\lambda^{\bar p}_i(t)\sim \lambda_i(t)$ for $i<n-1-\bar p(n)$.

For $0\neq i\geq n-1-\bar p(n)$, we have 
\begin{equation}\label{E: MV IH point}
\to IH_i^{\bar p} (X(\Sigma)\cap N(\Sigma);\vg) \overset{i_*}{\to}IH_i^{\bar
p} (X(\Sigma);\vg) \to IH_i^{\bar p} (S^n;\vg) \to,
\end{equation}
where the first map is induced by inclusion. Thus, by the calculations above and the
homotopy equivalence properties for ordinary homology, the first map is equivalent to
the map $i_*: H_i(L(\Sigma)-L(\Sigma)\cap K;\vg) \to H_i(\overline{X(\Sigma)}-
\overline{X(\Sigma)}\cap K;\vg)$, which is the inclusion map of the usual exact
sequence of
a pair for a knot with a point singularity (see above or section \cite{GBF1}). We obtain the exact
sequence of polynomials
\begin{equation*}
\begin{CD}
@>>> \nu_{i}(t) @>>> \lambda_{i}(t) @>>> I\lambda_{i}^{\bar p} (t) @>>>
\nu_{i-1}(t) @>>>  
\end{CD}
\end{equation*}
from the exact sequence of modules \eqref{E: MV IH point}. Since the map $i_*$ there is the same as
the map of the usual homology sequence of the pair of the knot, the splitting of $\nu_i(t)$ and
$\lambda_i(t)$ into factors in the polynomial sequence is determined in the same manner, and we obtain
the sequence
\begin{equation*}
\begin{CD}
@>>> a_{i}(t)b_i(t) @>>> b_i(t)c_i(t) @>>> I\lambda_{i}^{\bar p} (t) @>>>
a_{i-1}(t)b_{i-1}(t) @>>>.
\end{CD}
\end{equation*}
Thus, for $i> n-1-\bar p(n)$, the polynomial sequence determines that $
I\lambda_{i}^{\bar p}(t)\sim c_i(t)a_{i-1}(t)\sim \mu_i(t)$.

For $i= n-1-\bar p(n)$, the exact sequence can be truncated to the exact sequence  
\begin{equation*}
\begin{CD}
@>>> \nu_{n-1-\bar p(n)}(t) @>>> \lambda_{ n-1-\bar p(n)}(t) @>>> I\lambda_{n-1-\bar 
p(n)}^{\bar p} (t) @>>> 0 
\end{CD}
\end{equation*}
due to the injectivity of $i_*$ in the Mayer-Vietoris sequence for $i= n-2-\bar p(n)$. So
calculating $I\lambda_{ n-1-\bar p(n)}^{\bar p} (t) $ from the exact polynomial sequence,
using the known factorizations of $ \lambda_{n-1-\bar p(n)}(t)$ and $ \nu_{n-1-\bar
p(n)}(t)$, gives $I\lambda_{ n-1-\bar p(n)}^{\bar p} (t) \sim c_{ n-1-\bar p(n)} (t)$.

This completes the proof of the proposition. 
\end{proof}

The result for knots with multiple point singularities is not quite as nice, but it will occur as a special case of the computations in the following section.

\subsection{Manifold singularities with trivial neighborhoods}\label{S: man triv}

We turn next to the cases where $\Sigma=\Sigma_{n-k-1}$ is a manifold, so that the
filtration of $S^n$ is
\begin{equation*}
S^n\supset K\supset \Sigma_{n-k-1}.
\end{equation*}
Furthermore, we assume that $N(\Sigma)$, the open regular neighborhood of $\Sigma$, has the
structure of a product. In other words, we assume there exists a sphere
knot $\ell \cong S^{k-2}\subset S^k$ so
that the space pair of the regular neighborhood of $\Sigma$, $(N(\Sigma),
N(\Sigma)\cap K)$, is homeomorphic to the product space $\Sigma\times c(S^k, \ell)$.

For
example, this will
allow us to compute the intersection Alexander polynomials of frame twist-spun knots (see \cite{GBF1} and \cite{GBF}). 
Note that if $\pi_1$ and $\pi_2$ are the projections of $\Sigma\times (S^k-\ell)$ to its
factors, then $\vg|\Sigma\times (S^k-\ell)\cong
\pi_1^*\vg|_{\Sigma}\otimes\pi_2^*\vg|_{S^k-\ell}$. For convenience of notation, however, we will
simply refer to each restricted coefficient system as $\vg$.

We first recall that in \cite{GBF} and \cite{GBF1} we extended the definitions of the polynomials $\lambda_i$, $\mu_i$, $\nu_i$, $a_i$, $b_i$, and $c_i$ as follows: Let $\bar N(\Sigma)$ denote the closed regular neighborhood of the singular set $\Sigma$, and let $D=S^n-\bar N(\Sigma)$. Then we define $\lambda_i$ as the polynomial of the module $H_i(D-D\cap K;\vg)$,  $\nu_i$ as the polynomial of the module $H_i(\bd \bar N(\Sigma)-(\bd \bar N(\Sigma)\cap K);\vg)$, and $\mu_i$ as the corresponding polynomial of the homology of the pair $ (D-D\cap K,\bd \bar N(\Sigma)\cap K)$ with coefficients $\vg$. Since  $D-D\cap K$ is homotopy equivalent to $S^n-K$, $\lambda_i$ is, in fact, just what we would expect to be the ordinary Alexander polynomial of the knot. However, $\nu_i$ is the polynomial of a locally-flat knotted manifold pair which does not necessarily consist of spheres and is not even necessarily connected. It is shown in \cite{GBF} and \cite{GBF1} that these modules are all torsion $\Gamma$-modules so that these definitions make sense. We can then define the subpolynomials $a_i$, $b_i$, and $c_i$ as above using the polynomial splitting in the exact sequence $\to \nu_i\to\lambda_i\to \mu_i\to$. In case $\Sigma$ is a point singularity, these definitions reduce to those discussed in the last section. 

Let $A_i$ denote the kernel of the inclusion map in long exact homology sequence of a knot, i.e. the module
whose polynomial we have referred to as $a_i(t)$. Let   $\mathfrak A_i^{\geq
k-\bar p(k+1)}$ denote the submodule of  $H_i(\bd \bar N(\Sigma)-(\bd \bar N(\Sigma)\cap K);\vg)$ given by 
\begin{align*}
\mathfrak A_i^{\geq
k-\bar p(k+1)}&\cong A_i\cap\left( \underset{\underset{0\neq s \geq k-\bar
p(k+1)}{i=r+s}}{\oplus} \left[H_r(\Sigma;\vg)\otimes H_s
(S^k-\ell;\vg)\right]\right.\\
&\qquad\qquad\qquad \left.\oplus \underset{\underset{0\neq s \geq k-\bar
p(k+1)}{i-1=r+s}}{\oplus}
\left[H_r(\Sigma;\vg)*H_s (S^k-\ell;\vg)\right]\right),
\end{align*}
where we have identified the latter terms as submodules of $H_i(\bd \bar N(\Sigma)-(\bd \bar N(\Sigma)\cap K);\vg)$ using the K\"unneth theorem and the triviality assumptions concerning the neighborhood of $\Sigma$. 
Let $\mathfrak a_i^{\geq k-\bar p(k+1)}(t)$ denote the polynomial of
$\mathfrak A^{\geq k-\bar p(k+1)}_i$, and let $\mathfrak
b_i^{\geq k-\bar
p(k+1)}(t)$  denote the polynomial of 
\begin{equation*}\underset{\underset{0\neq s \geq k-\bar
p(k+1)}{i=r+s}}{\oplus}
\left[H_r(\Sigma;\vg)\otimes H_s (S^k-\ell;\vg)\right]\oplus
\underset{\underset{0\neq s \geq k-\bar p(k+1)}{i-1=r+s}}{\oplus}
\left[H_r(\Sigma;\vg)*H_s (S^k-\ell;\vg)\right]
\end{equation*}
divided by $\mathfrak a_i^{\geq k-\bar
p(k+1)}(t)$. (See below for a formula for  $\mathfrak b_i^{\geq k-\bar p(k+1)}(t)$).
We will show below that this is a polynomial and that $ \mathfrak b_i^{\geq k-\bar p(k+1)}(t)|b_i(t)$. Call the quotient $\mathfrak  b_i^{< k-\bar p(k+1)}(t)$.
Then we can prove the following:

\begin{proposition}\label{P: man triv2}

Given a knot $K\subset S^n$ as above, with $\Sigma=\Sigma_{n-k-1}$ a manifold
and $(N(\Sigma), N(\Sigma)\cap K)\cong\Sigma\times c(S^k, \ell)$,
then 
\begin{equation*}
I\lambda_i^{\bar p}(t)\sim \mathfrak a_{i-1}^{\geq k-\bar p(k+1)}(t) \mathfrak b_i^{<
k-\bar p(k+1)}(t)c_i(t).
\end{equation*}
\end{proposition}

\begin{proof}

With the notation of Proposition \ref{P: ih point} adapted in the obvious manner, we can once again form the Mayer-Vietoris
sequence
\begin{equation}\label{E: man triv 2 MV}
\to IH_i^{\bar p} (X(\Sigma)\cap N(\Sigma);\vg) \overset{i_*}{\to}IH_i^{\bar
p}
(X(\Sigma);\vg)\oplus IH_i^{\bar p} (N(\Sigma);\vg) \to IH_i^{\bar p}
(S^n;\vg) \to.
\end{equation}

As before, we first claim that  $IH_i^{\bar p} (X(\Sigma);\vg)\cong
H_i(S^n-K;\vg)$. Let $n(K)\cong N(K)\cap X(\Sigma)$ and $x(K)\cong X(K) \cap
X(\Sigma)$. Then we have a Mayer-Vietoris sequence 
{\small
\begin{equation*}
\to IH_i^{\bar p} (x(K)\cap n(K);\vg) \overset{i_*}{\to}IH_i^{\bar p}
(x(K);\vg)\oplus IH_i^{\bar p} (n(K);\vg) \to IH_i^{\bar p} (X(\Sigma);\vg) \to.
\end{equation*}}

But $x(K)$ does not intersect the singular set and is homotopy equivalent
to the knot complement $S^n-K$, so $IH_i^{\bar p} (x(K);\vg)\cong
H_i(S^n-K;\vg)$. To prove the claim, it is thus sufficient to show that
the map $i_*: IH_i^{\bar p} (x(K)\cap n(K);\vg) \to IH_i^{\bar p}
(n(K);\vg)$ induced by inclusion is an isomorphism. Now, since the
embedding of $K$ is locally-flat away from $\Sigma$, $n(K)$ is homeomorphic to a bundle of
disks $D^2\cong c(S^1)$, and once again there exists a trivialization of
this bundle as given in the proof of the existence of a Seifert surface
(\cite[Prop. 4.2]{GBF1}). Thus, if $\kappa=K\cap n(K)$, $n(K)\cong
\kappa\times c(S^1)$ and $ x(K)\cap n(K)\cong \kappa\times S^1\times \R$.
Since $\kappa$ is an unstratified  manifold, we can employ the K\"{u}nneth theorem for intersection homology, as
proven in \cite[Theorem 4]{Ki}, observing that the usual functorial
naturality holds since the theorem is proven by verifying the existence of
an ``Eilenberg-Zilber''-type chain map which induces the appropriate
homology isomorphisms. We obtain a diagram with split exact rows (in
which all coefficients are the suitable restrictions of $\vg$):

{\scriptsize
\begin{equation*}
\begin{diagram}
0&\rTo & \underset{i=r+s}{\oplus} H_r(\kappa)\otimes IH_s^{\bar
p}(S^1\times \R) &\rTo & IH_i^{\bar p} (x(K)\cap n(K)) &\rTo
&\underset{i-1=r+s}{\oplus}  H_r(\kappa)* IH_s^{\bar p}(S^1\times
\R) &\rTo&0\\
& & \dTo^{\oplus(\text{id}_*\otimes j_*)} & & \dTo^{i_*} &
&\dTo^{\oplus(\text{id}_*\otimes j_*)} & &\\
0&\rTo & \underset{i=r+s}{\oplus}  H_r(\kappa)\otimes IH_s^{\bar
p}(c(S^1)) &\rTo & IH_i^{\bar
p} (n(K)) &\rTo &\underset{i-1=r+s}{\oplus}  H_r(\kappa)*
IH_s^{\bar p}(c(S^1)) &\rTo&0.
\end{diagram}
\end{equation*}}

\noindent But we saw in the proof of Proposition \ref{P: locally-flat} that $j_*$ is an
isomorphism. Therefore,
the left and right vertical maps are isomorphisms, and so $i_*$ is an isomorphism by the five
lemma. Note that the proof of the claim does not rely on any of the assumptions concerning the number
or properties of the strata of the knot and therefore holds for any knot. In other words, we have
shown that $IH_*^{\bar p}(X(\Sigma);\vg)\cong H_*(S^n-K;\vg)$ for any knot. 

Similarly, due to the hypotheses of the theorem, $N(\Sigma)\cong\Sigma\times
c(S^k)$ and $N(\Sigma)\cap X(\Sigma)\cong\Sigma\times S^k\times \R$, and,
because $\Sigma$ is a manifold, we can again use the K\"unneth theorem to obtain
\begin{align*}
IH_i^{\bar p} (N(\Sigma);\vg) &\cong \underset{i=r+s}{\oplus} \left[ 
H_r(\Sigma;\vg)\otimes IH_s^{\bar p}(c(S^k);\vg)\right]\\
&\qquad \oplus
\underset{i-1=r+s}{\oplus} \left[ H_r(\Sigma;\vg)* IH_s^{\bar
p}(c(S^k);\vg)\right]\\
IH_i^{\bar p} (N(\Sigma)\cap X(\Sigma);\vg) &\cong \underset{i=r+s}{\oplus} \left[ 
H_r(\Sigma;vg)\otimes IH_s^{\bar p}(S^k\times\R;\vg)\right]\\
&\qquad\oplus
\underset{i-1=r+s}{\oplus} \left[ H_r(\Sigma;\vg)* IH_s^{\bar p}(S^k\times
\R;\vg)\right].
\end{align*}

Since the link knot pair $\ell\subset S^k$ is locally-flat, $IH_s^{\bar
p}(S^k\times \R;\vg)\cong H_s(S^k-\ell;\vg)$. For $IH_s^{\bar p}(c(S^k);\vg)$, we can use the cone formula of \cite{Ki} once again to obtain
\begin{equation*}
IH_s^{\bar p}(c(S^k);\vg)\cong
\begin{cases}
0, & 0\neq s\geq k-\bar p(k+1)\\
IH_s^{\bar p}(S^k;\vg)\cong H_s(S^k-\ell;\vg),& s=0, s< k-\bar p(k+1).
\end{cases}
\end{equation*}

Note also that the embedding of $K$ is locally flat in $N(\Sigma)\cap X(\Sigma)$ and therefore    $IH_i^{\bar p} (N(\Sigma)\cap X(\Sigma);\vg)\cong H_i( N(\Sigma)\cap X(\Sigma)-K;\vg)$. But this space is homemorphic to the product $[\bd \bar N(\Sigma)-(\bd \bar N(\Sigma)\cap K)]\times \R$ by the generalized annulus property, so using homotopy equivalence, $IH_i^{\bar p} (N(\Sigma)\cap X(\Sigma);\vg)$ is isomorphic
to the ordinary homology module of the link complement $\bd \overline{N(\Sigma)}-K$ with
corresponding polynomial $\nu_i(t)$.

 Meanwhile, if $\xi_{il}$ are the Alexander invariants of the
knot $\ell\subset
S^k$ and $H_i(\Sigma;\vg)=\Gamma^{\mathfrak B_i}\oplus \oplus_l \Gamma/{\zeta_l}$, then the polynomial
of $ IH_i^{\bar p} (N(\Sigma);\vg)$ is 
\begin{equation*}
\prod_{\underset{s=0, s< k-\bar p(k+1)}{i=r+s}}\left[\prod_l
\xi_{sl}^{\mathfrak B_r}\cdot\prod_{j,l}d(\zeta_{rj},\xi_{sl})\right]\cdot \prod_{\underset{s=0, s<
k-\bar p(k+1)}{i-1=r+s}} \left[\prod_{j,l}d(\zeta_{rj},\xi_{sl})\right],
\end{equation*}
where $d(\cdot,\cdot)$
indicated
the greatest common divisor in $\Gamma$. So, from the Mayer-Vietoris sequence \eqref{E: man triv 2 MV}, we obtain the long exact
polynomial sequence
{\footnotesize
\begin{equation*}
\to\nu_i \to \lambda_i \cdot \prod_{\underset{s=0, s< k-\bar
p(k+1)}{i=r+s}}\left[\prod_l \xi_{sl}^{\mathfrak B_r} 
\cdot\prod_{j,l}d(\zeta_{rj},\xi_{sl})\right]\cdot \prod_{\underset{s=0, s<
k-\bar p(k+1)}{i-1=r+s}} \left[\prod_{j,l}d(\zeta_{rj},\xi_{sl})\right]
\to I\lambda_i^{\bar p} \to.
\end{equation*}}

In order to calculate the intersection Alexander polynomial, we now need only determine the
polynomial of the kernel of $i_*$ in the Mayer-Vietoris sequence \eqref{E: man triv 2 MV}.  The map to
the first summand,
\begin{equation*}
i_*:IH_i^{\bar p}(X(K)\cap N(K);\vg) \to H_i(S^n-K;\vg)
\end{equation*}
is induced by the inclusion of $X(\Sigma)\cap N(\Sigma)$ into $X(\Sigma)$, and on homology this
induces the map $H_i(\Sigma\times (S^k-\ell);\vg)\to H_i(X(K);\vg)$ which, by homotopy equivalences,
is the standard map from the homology of the link exterior of $\Sigma$ to the homology of the exterior
of $K$. This is isomorphic to the inclusion map in the long exact homology sequence of the knot. 

The map to the second summand is the middle vertical map in the following
diagram induced by the naturality of the K\"unneth theorem (in
which all coefficients are the suitable restrictions of $\vg$):

{\scriptsize
\begin{equation}\label{D: Kunneth}
\begin{diagram}
0& \rTo & \underset{i=r+s}{\oplus} H_r(\Sigma)\otimes IH_s^{\bar
p}(S^k\times\R) &\rTo & IH_i^{\bar p} (N(\Sigma)\cap X(\Sigma)) &\rTo &
\underset{i-1=r+s}{\oplus}  H_r(\Sigma)* IH_s^{\bar p}(S^k\times \R) &
\rTo&0 \\
&&\dTo^{\oplus (\text{id}_*\otimes j_*)} &&\dTo^{i_*}&&\dTo^{\oplus
(\text{id}_* * j_*)} &&\\
0& \rTo & \underset{i=r+s}{\oplus}  H_r(\Sigma)\otimes IH_s^{\bar
p}(c(S^k)) &\rTo & IH_i^{\bar p}
(N(\Sigma)) &\rTo&\underset{i-1=r+s}{\oplus}  H_r(\Sigma)* IH_s^{\bar
p}(c(S^k)) & \rTo&0 . \end{diagram}
\end{equation}}

\noindent The leftmost and rightmost vertical maps are all induced by
inclusions and are the direct
sums of $\text{id}_{\Sigma*}\otimes j_*$ or $\text{id}_{\Sigma*} * j_*$, where $j_*$ is the homology
map induced by the
inclusion $S^k\times \R\to c(S^k)$. By the calculation of intersection homology of a cone in \cite[\S
4]{Ki2}, the maps $j_*$ are surjections on homology. Thus so are the maps $\text{id}\otimes j_*$ by
the right exactness of the tensor product functor.  Thus \emph{cok}$(\oplus \text{id}_*\otimes
j_*)=0$, and the serpent lemma yields a short exact sequence
\begin{equation*}
\begin{CD}
0@>>>\text{ker}[\oplus (\text{id}_*\otimes j_*)] @>>>\text{ker}[i_*] @>>>\text{ker}[\oplus (\text{id}_* * j_*)] @>>> 0.
\end{CD}
\end{equation*}
This sequence is split because the splitting map $ IH_i^{\bar p} (N(\Sigma)\cap X(\Sigma);\vg)\to
\underset{i=r+s}{\oplus} H_r(\Sigma;\vg)\otimes IH_s^{\bar p}(S^k\times\R;\vg)$ of the K\"unneth short
exact sequence restricts to give a splitting of the sequence of kernels. Hence, the kernel of $i_*$
is
the direct sum of the kernels of the left- and righthand maps of diagram \eqref{D: Kunneth}.
Furthermore, each of these is the direct sum of the kernels of the maps $\text{id}_*\otimes j_*$ or
$\text{id}_* * j_*$. But once again, the calculations of \cite{Ki} and \cite{Ki2} tell us that $j_*$
is
an isomorphism for $s=0$ or $s<k-\bar p(k+1)$ and it is the zero map for $0\neq s \geq k-\bar p(k+1)$.
From this, we conclude that the kernel of the map $i_*$ to the summand $H_i(S^n-K;\vg)$ is
\begin{equation*}
\underset{\underset{0\neq s \geq k-\bar p(k+1)}{i=r+s}}{\oplus}
\left[H_r(\Sigma;\vg)\otimes H_s
(S^k-\ell;\vg)\right]\oplus \underset{\underset{0\neq s \geq k-\bar p(k+1)}{i-1=r+s}}{\oplus}
\left[H_r(\Sigma;\vg)*H_s (S^k-\ell;\vg)\right].
\end{equation*}

So now, we let $A_i$ denote the
kernel of the inclusion map in long exact homology sequence of the knot. It is a submodule of
$H_i(\Sigma\times (S^k-\ell);\vg)$ with polynomial $a_i(t)$. We let 
\begin{align*}
\mathfrak A_i^{\geq k-\bar p(k+1)}&\cong A_i\cap
\left(\underset{\underset{0\neq s
\geq k-\bar p(k+1)}{i=r+s}}{\oplus} \left(H_r(\Sigma;\vg)\otimes H_s
(S^k-\ell;\vg)\right)\right.\\
&\qquad\qquad\qquad \left.\oplus \underset{\underset{0\neq s \geq k-\bar
p(k+1)}{i-1=r+s}}{\oplus} \left(H_r(\Sigma;\vg)*H_s
(S^k-\ell;\vg)\right)\right),
\end{align*}
the kernel of $i_*$ in the Mayer-Vietoris sequence \eqref{E: man triv 2 MV}.
Let
$\mathfrak a_i^{\geq k-\bar p(k+1)}(t)$ denote the polynomial of $\mathfrak A^{\geq k-\bar
p(k+1)}_i$, and let $\mathfrak b_i^{\geq k-\bar p(k+1)}(t)$ denote the polynomial of
\begin{equation*}
\underset{\underset{0\neq s \geq k-\bar p(k+1)}{i=r+s}}{\oplus}
\left(H_r(\Sigma;\vg)\otimes H_s
(S^k-\ell;\vg)\right)\oplus \underset{\underset{0\neq s \geq k-\bar p(k+1)}{i-1=r+s}}{\oplus} \left(H_r(\Sigma;\vg)*H_s (S^k-\ell;\vg)\right)
\end{equation*}
divided by $\mathfrak a_i^{\geq k-\bar
p(k+1)}(t)$. In other words,
\begin{equation*}
\mathfrak b_i^{\geq k-\bar p(k+1)}=\frac{\prod_{\underset{0\neq s \geq k-\bar
p(k+1)}{i=r+s}}\left[\prod_l\xi_{sl}^{\mathfrak B_r} \prod_{jl} d(\zeta_{rj},\xi_{sl})\right]\cdot \prod_{\underset{0\neq s \geq k-\bar p(k+1)}{i-1=r+s}}\left[\prod_{jl} d(\zeta_{rj},\xi_{sl})\right]}{ \mathfrak a_i^{\geq k-\bar p(k+1)}}.
\end{equation*} 
Then $\mathfrak a_i^{\geq k-\bar p(k+1)}(t)|a_i(t)$ and, from the calculations of
the last paragraph, this is the factor of the exact sequence of
polynomials which is shared by $\nu_i(t)$ and $I\lambda_{i+1}^{\bar
p}(t)$. 

We also claim that $\mathfrak b_i^{\geq k-\bar p(k+1)}(t)|b_i(t)$.
To see
this, note that by its definition, $\mathfrak b_i^{\geq k-\bar p(k+1)}(t)$ is the
polynomial of the module 
\begin{align*}
\left[\underset{\underset{0\neq s
\geq k-\bar p(k+1)}{i=r+s}}{\oplus} \left(H_r(\Sigma;\vg)\otimes H_s
(S^k-\ell;\vg)\right)\oplus \underset{\underset{0\neq s \geq k-\bar
p(k+1)}{i-1=r+s}}{\oplus} \left(H_r(\Sigma;\vg)*H_s
(S^k-\ell;\vg)\right)\right]/\mathfrak A_i^{\geq
k-\bar p(k+1)}. 
\end{align*}
But this module is isomorphic to
\begin{multline*}
IH_i^{\bar p}(X(\Sigma)\cap N(\Sigma);\vg)/\left(A_i +
\left[\underset{\underset{s=0, s
<k-\bar p(k+1)}{i=r+s}}{\oplus} \left(H_r(\Sigma;\vg)\otimes H_s
(S^k-\ell;\vg)\right)\right.\right.\\
\qquad\qquad\qquad \left.\left.\oplus \underset{\underset{s=0, s < k-\bar
p(k+1)}{i-1=r+s}}{\oplus} \left(H_r(\Sigma;\vg)*H_s
(S^k-\ell;\vg)\right)\right]\right), 
\end{multline*}
which is a quotient module of $IH_i^{\bar p}(X(\Sigma)\cap
N(\Sigma);\vg)/A_i$, whose polynomial is $b_i(t)$, because this intersection
homology module is the module of the link complement, as previously
noted. But from the polynomial sequence associated to a short exact
sequence of modules, it is clear that the polynomial of a quotient of
a module must divide the polynomial of the module.  Thus $\mathfrak b_i^{\geq k-\bar
p(k+1)}(t)|b_i(t)$.  

Denote by $\mathfrak b_i^{< k-\bar 
p(k+1)}(t)$ the quotient $b_i(t)/\mathfrak b_i^{\geq k-\bar
p(k+1)}(t)$.  
Then the polynomial shared by $I\lambda_{i}^{\bar
p}(t)$ and the ``middle'' term of the exact polynomial sequence will be
\begin{align*}
&\lambda_i\left[\prod_{\underset{s=0, s< k-\bar
p(k+1) }{i=r+s}}\left[\prod_l\xi_{sl}^{\mathfrak B_r} \prod_{jl}
d(\zeta_{rj},\xi_{sl})\right] \prod_{\underset{s=0, s< k-\bar p(k+1) }{i-1=r+s}}\left[\prod_{jl} d(\zeta_{rj}\xi_{sl})\right] \right] \div\frac{\nu_i}{\mathfrak a_i^{\geq k-\bar p(k+1)}}\\ 
&=
\frac{ b_i c_i \cdot{\displaystyle\prod_{\underset{s=0, s< k-\bar
p(k+1) }{i=r+s}}}\left[\prod_l\xi_{sl}^{\mathfrak B_r}\cdot \prod_{jl}
d(\zeta_{rj},\xi_{sl})\right]\cdot {\displaystyle \prod_{\underset{s=0, s<
k-\bar
p(k+1) }{i-1=r+s}}}\left[\prod_{jl} d(\zeta_{rj},\xi_{sl})\right]}{\mathfrak   
b_i^{\geq k-\bar p(k+1)}{\displaystyle \prod_{\underset{s=0, s< k-\bar
p(k+1) )}{i=r+s}}}\left[\prod_l\xi_{sl}^{\mathfrak B_r}\cdot\prod_{jl}
d(\zeta_{rj},\xi_{sl})\right]\cdot {\displaystyle\prod_{\underset{s=0, s<
k-\bar p(k+1) }{i-1=r+s}}}\left[\prod_{jl} d(\zeta_{rj},\xi_{sl})\right]}
\\&=\frac{b_i c_i}{ \mathfrak b_i^{\geq k-\bar p(k+1)}}
\\&= \mathfrak b_i^{< k-\bar p(k+1)}c_i.
\end{align*}

Thus we conclude that 
\begin{equation*}
I\lambda_i^{\bar p}(t)\sim \mathfrak a_{i-1}^{\geq k-\bar p(k+1)}(t) \mathfrak b_i^{< k-\bar p(k+1)}(t)c_i(t).
\end{equation*}
\end{proof}

If we assume that the local system of coefficients restricted to
$\Sigma$ is simple, i.e. we can identify $\Sigma$ with $\Sigma\times *\subset \Sigma \times (S^k-\ell)$ in such a way that the action of the fundamental group on $\vg|_{\Sigma}$ is
trivial, then the results of the previous proposition simplify slighly. This is the situation which occurs, for example, for frame-spun knots (see \cite[\S 4.3]{GBF1}). We simply state the results of this special case of the previous proposition.

Let $\xi_i(t)$ denote the $i$th Alexander
polynomial of the link knot $\ell$. With the above assumptions, and identifying $\Sigma$ with $\Sigma\times *\subset \Sigma\times (S^k-\ell)$, $H_i(\Sigma;\vg)$ is the free $\Gamma$ module, $\Gamma^{\beta_i}$, where $\beta_i$ is the $i$th
Betti
number of $\Sigma$. Then, using the K\"unneth theorem to calculate the
Alexander module of the link pair complement of $\Sigma$ (i.e. $\bd\overline{N(\Sigma)}-K$), we have
$\nu_i(t)=a_i(t)b_i(t)=\prod_{i=r+s}\xi_i(t)^{\beta_r}$ as the polynomial of the
module $H_i(\Sigma\times (S^k-\ell);\vg)\cong \oplus_{i=r+s} H_r(\Sigma;\Gamma)\otimes
H_s(S^k-\ell;\vg|S^k)$. Let $A_i$ denote the kernel of the inclusion map of the link
complement in the long exact
homology sequence of a knot, i.e. the module whose polynomial we have referred to as
$a_i(t)$. This, of course, is a submodule of $H_i(\Sigma\times (S^k-\ell);\vg)$. Let
$A^{\geq k-\bar p(k+1)}_i\cong A_i\cap \underset{\underset{0\neq s \geq k-\bar
p(k+1)}{i=r+s}}{\oplus} H_r(\Sigma;\Gamma)\otimes H_s(S^k-\ell;\vg|_{S^k-\ell})$. Let $a_i^{\geq
k-\bar
p(k+1)}(t)$ denote the polynomial of $A^{\geq k-\bar p(k+1)}_i$, and let $b_i^{\geq
k-\bar p(k+1)}(t)$ denote the polynomial of $\underset{\underset{0\neq s \geq k-\bar
p(k+1)}{i=r+s}}{\oplus} H_r(M;\vg)\otimes H_s(S^k-\ell;\vg)$ divided by $a_i^{\geq
k-\bar p(k+1)}(t)$. In other words,
\begin{equation*}
b_i^{\geq k-\bar p(k+1)}(t)=\frac{\prod_{\underset{0\neq s \geq k-\bar p(k+1)}{i=r+s}}\xi_s(t)^{\beta_r}}{ a_i^{\geq k-\bar p(k+1)}(t)}.
\end{equation*}
Then $ b_i^{\geq k-\bar p(k+1)}(t)|b_i(t)$, and we call the quotient $
b_i^{< k-\bar p(k+1)}(t)$. Then the preceding proposition specializes to the following:

\begin{proposition}\label{P: man triv1}
Given a knot $K\subset S^n$ as above, with
$\Sigma=\Sigma_{n-k-1}$ a manifold and
$(N(\Sigma), N(\Sigma)\cap K)\cong\Sigma\times c(S^k, \ell)$, suppose also that,
taking $\Sigma\cong \Sigma\times *$ for $*\in S^k-\ell$, $\vg|_\Sigma$ is a simple
system of local coefficients (so that the action of each element of $\pi_1(\Sigma)$ on
$\Gamma$ is trivial).  Then
\begin{equation*}
I\lambda_i^{\bar p}(t)\sim a_{i-1}^{\geq k-\bar p(k+1)}(t) b_i^{< k-\bar
p(k+1)}(t)c_i(t).
\end{equation*}
\end{proposition}

In this case we have the following corollaries:

\begin{corollary}\label{C: stable range}
For a knot, $K$, as in Proposition \ref{P: man triv1}:
\begin{enumerate}
\item For $i<k-\bar p(k+1)$, $I\lambda_i^{\bar p}(t)\sim \lambda_i(t)$.
\item If $\bar p(k+1)\leq 1$, then  $I\lambda_i^{\bar p}(t)\sim \lambda_i(t)$ for all
$i$. More generally, if $H_i(S^k-\ell;\vg)=0$ for $i\geq j$ and $\bar p(k+1)\leq k-j$, then
$I\lambda_i^{\bar p}(t)\sim \lambda_i(t)$ for all
$i$.
\item If $i\geq n-\bar p(k+1)$, then $I\lambda_i^{\bar p}(t)\sim \mu_i(t)$.
\end{enumerate}
\end{corollary}
\begin{proof}
\begin{enumerate}
\item From the definitions, it is apparent that in this case $A_{i-1}^{\geq k-\bar
p(k+1)}=0$ because, with the assumptions of the proposition, $i-1=r+s$ and $s\geq k-\bar p(k+1)$
together
imply that $r<0$ (in fact that $r<-1$, but we shall need this extra room
shortly).  Therefore, $A_{i-1}^{\geq k-\bar p(k+1)}=A_{i-1}\cap 0$. Hence, again
from the definitions, $a_{i-1}^{\geq k-\bar p(k+1)}\sim 1$. Similary, these arguments
hold for  $A_{i}^{\geq k-\bar p(k+1)}$ and   $a_{i}^{\geq k-\bar p(k+1)}$ (here
again using the implication $r<0$) so that, by definition, we also obtain
$b_{i}^{\geq k-\bar p(k+1)}\sim 1$ and hence $b_{i}^{< k-\bar p(k+1)}\sim
b_i(t)$. Therefore, by the proposition,  $I\lambda_i^{\bar p}(t)\sim
b_i(t)c_i(t)\sim \lambda_i(t)$.

\item If  $\bar p(k+1)\leq 1$, then $k-\bar p(k+1)\geq k-1$, thus by the definitions
and
the fact that $H_s(S^k-\ell;\vg)=0$ for $s\geq k-1$, $A_i^{\geq k-\bar p(k+1)}=0$ for
all $i$. Similiarly, with $\bar p(k+1)\leq k-j$, we get  $k-\bar
p(k+1)\geq j$ so that $A_i^{\geq k-\bar p(k+1)}=0$ for all $i$ if $H_i(S^k-\ell;\vg)=0$ for
$i\geq j$. The rest of the argument now follows as for the previous item.

\item If $i\geq n-\bar p(k+1)$, $i=r+s$ or $i-1=r+s$, and $s<k-\bar p(k+1)$, then
$r>n-k-1$. Therefore, $A_i^{\geq k-\bar p(k+1)}=A_i$ and  $A_{i-1}^{\geq k-\bar
p(k+1)}=A_{i-1}$. It follows then that $a_i^{\geq k-\bar p(k+1)}\sim a_i$, $a_{i-1}^{\geq
k-\bar p(k+1)}\sim a_{i-1}$,  $b_i^{\geq k-\bar p(k+1)}\sim 1$, and $b_{i-1}^{\geq    
k-\bar p(k+1)}\sim 1$. Therefore, by the proposition,  $I\lambda_i^{\bar p}(t)\sim
a_{i-1}(t)c_i(t)\sim\mu_i(t)$.
 \end{enumerate}
\end{proof}

In the general case, where the coefficient bundle is not simple, the corollaries generalize as follows:

\begin{corollary}
For a knot, $K$, as in Proposition \ref{P: man triv2}:                      
\begin{enumerate}
\item For $i<k-\bar p(k+1)$, $I\lambda_i^{\bar p}(t)\sim \lambda_i(t)$.   
\item If $\bar p(k+1)\leq 1$, then  $I\lambda_i^{\bar p}(t)\sim \lambda_i(t)$ for all
$i$.  More generally, if $H_i(S^k-\ell;\vg)=0$ for $i\geq j$ and $\bar p(k+1)\leq k-j$, then
$I\lambda_i^{\bar p}(t)\sim \lambda_i(t)$ for all
$i$.
\item If $i\geq n-\bar p(k+1)+1$, then $I\lambda_i^{\bar p}(t)\sim \mu_i(t)$.
\end{enumerate}
\end{corollary}
\begin{proof}
The proof is essentially the same as that of Corollary \ref{C: stable range}. For the first 
two items, it is
easy to check that, in the definition of $\mathfrak A_i^{k-\bar p(k+1)}$, the
torsion product terms are also $0$ in the ranges for which we checked above that  the tensor
product terms are $0$. For the last item, we need to vary the range slightly to account for
the fact that the torsion product terms of $\mathfrak A_i^{k-\bar p(k+1)}$ have total degree
$i-1$ and not $i$.   
\end{proof}

\section{Spectral sequences and theorems on prime
components for knots with manifold singularities}\label{S: SS,primes} 

In the following sections, we will use spectral sequences to derive some
results concerning what prime elements in $\Gamma$ may arise as factors of
the intersection Alexander polynomials of knots with one singular
stratum. Initially, for simplicity, we will assume that the neighborhood of
this stratum can be given the structure of a fiber bundle. Then the
computations can proceed by showing how the Leray-Serre spectral sequence can be
used to compute the  intersection homology of a fiber bundle
with an unstratified manifold as the base space and a stratified
pseudomanifold as the fiber.

As we will discuss at the end of this section, the assumption of the existence of a bundle neighborhood is unnecessary. In fact, as proven in \cite{GBF3}, there always exists a spectral sequence for computing the intersection homology of the regular neighborhood of the bottom stratum of a stratified PL-pseudomanifold, and the $E^2$ terms consist of the homology of the bottom stratum with coefficients in a bundle whose stalks are given by the intersection homology of the cone on the link of the stratum. This generalizes  the $E^2$ terms of spectral sequences we calculate here. Thus in  each of the following results which depends on the hypothesis of a fiber bundle neighborhood, this particular hypothesis can be ignored. However, we leave it in for now in order to provide a more coherent framework and also as an excuse to develop a sheaf theoretic Leray-Serre spectral sequence for intersection homology which applies to more general base spaces than PL-pseudomanifolds and also to  sheaf theoretic superperverse intersection homology. We will not need to treat such general base spaces in our applications (our manifolds will always be triangulable), but we begin with the sheaf theoretic approach
mainly to establish the existence of the spectral sequence in this
generality and to illustrate its application.

\subsection{The sheaf theoretic spectral sequence for the intersection
homology of a fiber bundle}

We will construct a sheaf theoretic spectral sequence for  the intersection homology of a fiber bundle with unstratified manifold base space and paracompact stratified fiber. We begin with a lemma.

\begin{lemma}\label{L: restrict} 
Let $A$ be a paracompact subspace of a
paracompact space $X$ and $\Phi$ a family of paracompactifying supports on
$X$
(for example the collection of closed sets of $X$). Let $\mathcal N$ be a
collection of open subspaces of $X$ containing $A$ and directed downward by
inclusion. Assume that, for each $K\in \Phi|X-A$, there is an $N\in \mathcal N$
with $N\subset X-K$. Suppose $A$ is $\Phi$-taut. Lastly, suppose that $\mc L^*$
is a bounded differential sheaf. Then there is an isomorphism induced by
restriction:
\begin{equation}
\vartheta: \dlim_{N\in \mc N} \bb  H^*_{\Phi\cap N}(N;\mc L^*|N)\to \bb H_{\Phi\cap A}^*(A;\mc L^*|A), 
\end{equation}
where $\bb H$ represents hypercohomology (see \cite[p. 213]{GelMan}).
\end{lemma}
\begin{proof}
By \cite[II.10.6]{BR2}, there is an isomorphism for each $\mc L^i$,
\begin{equation}
\theta: \dlim_{N\in \mc N}  H^*_{\Phi\cap N}(N;\mc L^i|N)\to  H_{\Phi\cap A}^*(A;\mc L^i|A), 
\end{equation}
induced by restriction (here $H^*$ denotes sheaf cohomology).  Let $\mc J^{**}$
be the Cartan-Eilenberg resolution of $\mc L^*$ given by the sheaves
$\mc J^{p,q}=\mc C^p(X;\mc L^q)$ (see \cite{BR2} for the definition of these flabby
sheaves). Then, letting $\mc J^*$ be
the single complex associated to the double complex, $\bb H^*_{\Phi}(X;\mc
L^*)=H^*(\Gamma_{\Phi}(X; \mc J^*))$, and, more generally, for any
left-exact functor $F$
such that the $\mc J^{**}$ are $F$-acyclic, $H^i(F(\mc J^*))$ is the $i$th
right derived
functor of $F$ on $\mc L^*$. See, e.g.,
\cite[\S III.7]{GelMan}.
Then, 
\begin{align*}
\dlim_{N\in \mc N} \bb H^*_{\Phi\cap N}(N;\mc L^*|N)&=\dlim  H^*(\Gamma_{\Phi\cap
N}(N;\mc
J^*|N))\\
&=H^*(\dlim \Gamma_{\Phi\cap N}(N;\mc J^*|N)).
\end{align*}
By spectral sequence theory, this is the abutment of a spectral sequence with $E_2$ term
\begin{align*}
E_2^{pq}&=H^p_{II}(H_I^q(\dlim \Gamma_{\Phi\cap N}(N;\mc J^{**}|N)))\\
&= H^p_{II}(\dlim H_I^q (\Gamma_{\Phi\cap N}(N;\mc J^{**}|N)))\\
&= H^p_{II}(\dlim H_{\Phi\cap N}^q (N;\mc L^*|N)),
\end{align*}
where the last identity holds because the $\mc J^{*,q}|N$ form a flabby resolution of the $\mc
L^q|N$ because restriction is an exact functor and the restriction of a flabby sheaf to an
open set is flabby.

On the other hand, if $\mc J_A^{**}$ is a Cartan-Eilenberg resolution of $\mc L|A$ and
$\mc J_A^*$ is the associated single complex, then
\begin{align*}
\bb H_{\Phi\cap A}^*(A;\mc L^*|A) &= H^*(\Gamma_{\Phi\cap A}(A;\mc J_A^*)),
\end{align*}
which is the abutment of a spectral sequence with $E_2$ term
\begin{align*}
E_2^{pq}&=H^p_{II}(H_I^q(\Gamma_{\Phi\cap A}(A;\mc J_A^{**}|N)))\\
&= H^p_{II}(H_{\Phi\cap A}^q (A;\mc L^*|A)).
\end{align*}

Now if $r:\Gamma_{\Phi}(\mc L^*)\to \Gamma_{\Phi\cap A}(\mc L^*|A)$ is the restriction of
sections, by \cite[IV.4.2]{BR2} $r$ induces the isomorphisms $\theta: 
 \dlim_{N\in \mc N}  H^*_{\Phi\cap N}(N;\mc L^i|N)\to  H_{\Phi\cap A}^*(A;\mc L^i|A)$. 
It is easy to check that $r$ is  a natural transformation
of functors (for example, it is obvious if we think of the sheaves as ``sheaf spaces'').
Therefore, we have a natural isomorphism of the $E_2$
terms of the spectral sequence, which proves the lemma.
\end{proof}

\begin{proposition}\label{P: SS}
Let $(E,B,F, \pi)$ be a fiber bundle with base space $B$ a manifold, total space $E$,
paracompact stratified fiber $F$, and
projection $\pi$ such that for sufficiently small open $U\subset B$, $\pi^{-1}(U)\cong
U\times F$, where the the stratification is given by $F_i\times U$, $F_i$ the strata
of $F$. Then, for any fixed
perversity, $\bar p$, which we omit from the notation, there is a spectral sequence abutting
to $IH^i_c(E;\vg)$ with $E_2$ term
\begin{equation*}
E_2^{p,q}=H_c^p(B; \mc{IH}^q_c(F;\vg|F)),
\end{equation*}
where $\mc{IH}^i_c(F;\vg|F)$ is a local coefficient system (sheaf) with stalks $IH^i_c
(F;\vg|F)$ and $c$ denotes the system of compact supports.
\end{proposition}
\begin{proof}

Let $\mc{IC}^*(E)$ be the sheaf of intersection chains on $E$ with appropriate local
coefficient system and perversity. $\mc{IC}^*(E)$ is soft \cite{Bo} and hence $\Phi$-soft
for any paracompactifying system of supports. Then $\pi_c \mc{IC}^*(E)$ is c-soft by
\cite[p. 493, Property h]{BR2}, taking there
$\Phi=\Psi=c$. So, applying \cite[2.1]{BR2} to the differential sheaf $\pi_c
\mc{IC}^*(E)$,
there is a spectral sequence which abuts to
$H^*(\Gamma_c(\pi_c\mc{IC}^*(E)))$ and which has
$E_2$ term $E_2^{p,q}=H_c^p(B; \mc H^q(\pi_c\mc{IC}^*(E)))$. By \cite[\S
IV.5]{BR2},
\begin{align*}
H^*(\Gamma_c(\pi_c\mc{IC}^*(E)))=H^*(\Gamma_{c(c)}(
\mc{IC}_*(E))=H^*(\Gamma_{c}( \mc{IC}^*(E))=IH^*_c(E). 
\end{align*}
 It remains to show that $\mc
H^q(\pi_c\mc{IC}^*(E))$ is locally constant and that its stalks,
$\mc
H^q(\pi_c\mc{IC}^*(E))_y$ for $y\in B$, are isomorphic to $IH^i_c (F;\vg|F)$.

As in the arguments in \cite[p. 213]{BR2}, $\pi_c\mc{IC}_*(E)$ is the sheaf generated by
the presheaf 
\begin{align*} 
U &\to \Gamma_{c\cap \pi^{-1}(U)}( \mc{IC}^*(E)|\pi^{-1}(U))\\
&=\Gamma_{c\cap \pi^{-1}(U)}( \mc{IC}^*(\pi^{-1}(U))), 
\end{align*} 
where the equality is due to the fact that the sheaf of
intersection chains on the open set $\pi^{-1}(U)$ is equal to the restriction to
$\pi^{-1}(U)$ of the sheaf of intersection chains on $E$ (for example, they each satisfy
the axioms for $\vg|\pi^{-1}(U)$).  Thus the derived homology sheaf $\mc
H^q(\pi_c\mc{IC}^*(E))$ is generated by the presheaf 
\begin{align*} 
U &\to H^q(\Gamma_{c\cap \pi^{-1}(U)}( \mc{IC}^*(\pi^{-1}(U)))) \\ &=IH^q_{c\cap
\pi^{-1}(U)}(\pi^{-1}(U)) 
\end{align*}
as $\mc{IC}_*(\pi^{-1}(U))$ is $c\cap \pi^{-1}(U)$-soft, $c\cap \pi^{-1}(U)$ being
paracompactifying. 

Now, the stalk $\mc H^q(\pi_c\mc{IC}^*(E))_y$ is thus by definition
\begin{equation*}
\dlim IH^q_{c\cap \pi^{-1}(U)}(\pi^{-1}(U))
=\dlim \bb{H}^q_{c\cap \pi^{-1}(U)}(\pi^{-1}(U); \mc{IC}^*(E)|\pi^{-1}(U)),
\end{equation*}
where the limit is taken over a descending series of open sets $U$ with $y\in U$. We can
now apply
Lemma \ref{L: restrict}: $\pi^{-1}(y)$ is $c$-taut because it is a closed set and
$c$ is paracompactifying \cite[p. 73]{BR2}.
The other condition on the $U$'s is obviously satisfied. 
Therefore, all of the conditions of Lemma
\ref{L: restrict} are satisfied, and  $\mc H^q(\pi_c\mc{IC}^*(E))_y$ 
is
isomorphic to $\bb H_{c\cap \pi^{-1}(y)}^*(\pi^{-1}(y);\mc
IC^*(E)|\pi^{-1}(y))$.

Take now the neighborhoods $U$ small enough that $\pi^{-1}(U)=U\times F$. Let
$f$ be the projection $U\times F\to F$, and let $\mc P^*(X; \mc{E})$ be the
Deligne sheaf (see \cite{Bo}) on the space $X$ with coefficient system, $\mc
E$, and perversity
$\bar p$. $\mc P^*$ is quasi-isomorphic to
$\mc{IC}^*$, and, because restriction is an exact functor, $\mc P^*|F$ is
quasi-isomorphic to $\mc{IC}^*|F$. By \cite[V.3.14]{Bo}, $\mc{P}^*(U\times
F;\vg|U\times F)=f^*\mc{P}^*(F;\vg|F)$. Therefore, if $i:F\to U\times F$ is the
inclusion homeomorphism which takes $F$ homeomorphically onto $\pi^{-1}(y)$,
then
\begin{align}\label{E: restrict}
\mc{P}^*(U\times F)|\pi^{-1}(y)&=i^*\mc{P}^*(U\times F)\\\notag
&=i^* f^*\mc{P}^*(F;\vg|F)\\\notag
&=(fi)^* \mc{P}^*(F;\vg|F)\\\notag
&=\text{id}^*\mc{P}^*(F;\vg|F)\\\notag
&=\mc{P}^*(F;\vg|F).\notag
\end{align}
Again by the Lemma \ref{L: restrict} and the fact that we can restrict the
issue to a small
neighborhood of $y\in B$, the stalk is therefore  
\begin{align*}
\bb{H}_{c\cap \pi^{-1}(y)}^*(\pi^{-1}(y);\mc IC^*(U\times F)|\pi^{-1}(y))&=
 \bb {H}_{c\cap F}^*(F;\mc {IC}^*(U\times F)|F)\\
&\cong\bb {H}_{c\cap F}^*(F;\mc P^*(U\times F)|F)\\
&=\bb {H}_{c\cap F}^*(F;\mc P^*(F))\\
&\cong \bb {H}_{c}^*(F;\mc IC^*(F))\\
&=IH_c^*(F;\vg|F),
\end{align*}
where we have used the fact that $F$ is closed in $u\times F$ to replace $c\cap F$ by $F$ and also the fact that quasi-isomorphisms induce hypercohomology isomorphisms for any system of supports.

That $\mc
H^q(\pi_c\mc{IC}^*(E))$ is locally constant follows, with the obvious
modifications, as in the arguments in \cite[p. 227-228]{BR2}
for the sheaf cohomology of a fiber bundle with coefficients in a single sheaf.

This completes the proof. 
\end{proof}

\begin{remark}\label{R: SS}
For superperversities, the above proof remains true for sheaf intersection homology if we
replace the differential sheaf
$\mc{IC}$ with a soft resolution of the Deligne sheaf. Hence we obtain an analogous spectral
sequence for superperverse sheaf intersection homology. 
\end{remark}

\subsection{Prime factors of the intersection Alexander polynomials}

Recall the following discussion from Section \ref{S: poly alg}:

Suppose that $M$ is a torsion $\Gamma$-module with submodule $N$. Associated
to the short exact sequence
\begin{equation*}
\begin{CD}
0 @>>> N @>>> M @>>> M/N @>>> 0,
\end{CD}
\end{equation*}
we have a short exact polynomial sequence (see
Section \ref{S: poly alg}) of the form 
\begin{equation*}
\begin{CD}
0 @>>> f  @>>> h @>>> g @>>> 0,
\end{CD}
\end{equation*}
where $f, g,h\in \Gamma$, $f$ is the polynomial associated to $N$, $h$ is the
polynomial associated to $M$, and
$g$ is the polynomial associated to $M/N$. Further, from the properties of
exact
polynomial sequences, we know that we must have $h=fg$.
It is immediate, therefore, that if a prime
$\gamma\in \Gamma$ divides $f$ or $g$, then it divides $h$. Conversely, if
it divides
$h$ then it must divide $f$ or $g$.
We can then drawing the following conclusion: Suppose that $A$ is a
subquotient of $M$ (i.e. a quotient module of a submodule of $M$). Then a
prime $\gamma\in \Gamma$ can divide the polynomial associated to $A$ only if
it divides the polynomial associated to $M$. For suppose $A=N/P$, where
$P\subset N\subset M$. If $\gamma$ divides the polynomial of $A$, then by the
above arguments it must divide the polynomial of $N$. But then similarly,
$\gamma$ must divide the polynomial of $M$.

We will use these elementary facts to identify what can be the prime factors
of the intersection Alexander polynomials of a knot whose singular set is a
manifold $\Sigma_{n-k-1}$ and with the property that the open regular
neighborhood of the singularity is a fiber bundle with fiber the
cone on the link knot $\ell$ given by $S^{k-2}\subset S^k$. For a given range
of dimensions, this will always be the case for a manifold singularity
(see \cite{Hi}). With these assumptions, we have the
following theorem:

\begin{theorem}\label{T: at one}
Let $K$ be a non-locally-flat knot with singularity $\Sigma=\Sigma_{n-k-1}$ a
manifold and such that $N(\Sigma)$ is a fiber bundle with base $\Sigma$
and  with fiber the
cone on the link knot $\ell$ given by $S^{k-2}\subset S^k$. Let $\xi_j$ be
the Alexander polynomials of the locally-flat link knot $\ell$. Let $\bar p$ be
a traditional perversity. Then, for
$0<i<n-1$ and for any
prime $\gamma\in \Gamma$,
$\gamma|I\lambda_i^{\bar p}$ only if $\gamma|c_i$ or $\gamma|\xi_s$ for
some $s$ such that $0\leq i-s \leq n-k$ and $0< s < k-1$. 
In other words, the prime factors of
$c_i$ and $\xi_s$, $s$ in the allowable range, are the only possible
prime factors of $I\lambda_i^{\bar p}$.
\end{theorem}
\begin{proof}
Once again, we will employ the Mayer-Vietoris sequence 
\begin{equation}\label{E: MVM}
\to IH_i^{\bar p} (X(\Sigma)\cap N(\Sigma);\vg) \overset{i_*}{\to}IH_i^{\bar
p}
(X(\Sigma);\vg)\oplus IH_i^{\bar p} (N(\Sigma);\vg) \to IH_i^{\bar p}
(S^n;\vg) \to
\end{equation}
with the notation as in Section \ref{S: ih and ord} (see Propositions
\ref{P: man triv1} and \ref{P: man triv2}). 

Also as before, the module $IH_i^{\bar p} (X(\Sigma);\vg)$ is isomorphic to
$H_i
(X(K);\vg)$, and  $IH_i^{\bar p} (X(\Sigma)\cap N(\Sigma);\vg)$ is isomorphic 
to the homology of the of the link complement of $\Sigma$, $\bd
\overline{N(\Sigma)}-K$. The arguments we applied in the proofs of
Propositions \ref{P: man triv1} and \ref{P: man triv2} to $X(\Sigma)$ apply
again here because the embedding of $K$ is locally flat in the complement of
$\Sigma$,
and the trivialization of the circle bundle over $K\cap X(\Sigma)$ restricts to
a trivialization of the circle bundle over $K\cap X(\Sigma)\cap N(\Sigma)$.
Therefore,  the same reasoning as applied in those proofs shows that
$IH_i^{\bar p} (X(\Sigma);\vg)\cong
H_i(X(\Sigma)-K;\vg)$ and  
$IH_i^{\bar p} (X(\Sigma)\cap N(\Sigma);\vg)\cong
H_i(X(\Sigma)\cap N(\Sigma)-K;\vg)$. But $X(\Sigma)-K$ is
homotopy equivalent to $S^n-K$ and $X(\Sigma)\cap N(\Sigma)-K$ is
homotopy equivalent, by the generalized annulus property (see
\cite[Proposition 1.5]{Sto}), to the link
complement $\bd \overline{N(\Sigma)}-K$. In particular,
$\lambda_i(t)$ is the polynomial associated to the module $IH_i^{\bar
p}(X(\Sigma);\vg)$ and $\nu_i(t)=a_i(t)b_i(t)$ is the polynomial associated to
$IH_i^{\bar
p} (X(\Sigma)\cap N(\Sigma);\vg)$.

We will see below that all of the terms of the
Mayer-Vietoris sequence \eqref{E: MVM} are torsion modules. Then from the
exact polynomial
sequence associated to the exact module sequence, we will know that
$I\lambda^{\bar p}_i$ is
the product of two polynomials, one dividing the polynomial associated to
$IH_i^{\bar p}
(X(\Sigma);\vg)\oplus IH_i^{\bar p} (N(\Sigma);\vg)$ and one dividing the
polynomial of $IH_{i-1}^{\bar p} (X(\Sigma)\cap N(\Sigma);\vg)$. Hence,
$\gamma|I\lambda^{\bar p}_i$ only if it divides one of these factors, and it
can divide the appropriate factor only if it divides the whole polynomial
associated to the respective module. Thus it suffices to determine which primes
divide the polynomials associated to the modules $IH_i^{\bar p}
(X(\Sigma);\vg)\oplus IH_i^{\bar p} (N(\Sigma);\vg)$ and $IH_{i-1}^{\bar p}
(X(\Sigma)\cap N(\Sigma);\vg)$.

The polynomial of $IH_i^{\bar p}
(X(\Sigma);\vg)\oplus IH_i^{\bar p} (N(\Sigma);\vg)$ is the product of the
polynomials of the summands, the former of which we have already identified
as $\lambda_i(t)$. We can actually do slightly better with this term. Since the map
$i_*$ of the Mayer-Vietoris sequence is induced by inclusion, we know that
its kernel must be
a submodule of the kernel of the restriction of $i_*$ to the first summand.
Thus the polynomial associated to the kernel of $i_*$ must divide $a_i$ (see the
discussion prior to Proposition \ref{P: man triv1}). Then from the exact
polynomial sequence
associated to the Mayer-Vietoris sequence, this implies that $b_i$ divides the
polynomial of $IH_i^{\bar p}(X(\Sigma);\vg)\oplus IH_i^{\bar p} (N(\Sigma);\vg)$. This
further implies that the polynomial factor which divides both $IH_i^{\bar
p}(X(\Sigma);\vg)\oplus IH_i^{\bar p} (N(\Sigma);\vg)$ and $IH_i^{\bar p}(S^n;\Gamma)$
is a factor of the product of the polynomial of $IH_i^{\bar p} (N(\Sigma);\vg)$ and $\lambda_i/b_i=c_i$. (Note,
however, that any prime factor of $b_i$ may yet occur in one of the other remaining
terms.) 

Therefore, it remains to identify the prime factors of the
polynomials associated to $IH_i^{\bar p} (N(\Sigma);\vg)$ and $IH_{i-1}^{\bar
p} (X(\Sigma)\cap N(\Sigma);\vg)$ in order to determine what other primes might divide
$I\lambda_i^{\bar p}$.

First, we show that each $IH_i^{\bar p} (N(\Sigma);\vg)$ is a
$\Gamma$-torsion module whose associated polynomial is divisible only by
primes that divide one of the $\xi_s$, $0\leq i-s \leq n-k-1$ and $0< s < k-1$.
For this we will employ the
spectral sequence of Proposition \ref{P: SS}. With our notation for
compactly supported singular or simplicial intersection homology and Borel's
\cite{Bo} notation for sheaf intersection cohomology (and dropping the
explicit perversity from each for simplicity), we have $IH_i=IH^{n-i}_c$,
where the latter is the sheaf intersection homology with compact supports and
we have assumed a space of dimension $n$. So, equivalently, we need to show
that $IH^{n-i}_c (N(\Sigma);\vg)$ is a $\Gamma$-torsion modules whose
associated polynomial has the desired properties. Since $N(\Sigma)$ is a
fiber bundle, we can employ Proposition \ref{P: SS} by which $IH^{n-i}_c
(N(\Sigma);\vg)$ is the abutment of a spectral sequence with $E_2$ terms
given by
\begin{equation*}
E_2^{p,q}=H_c^p(\Sigma; \mc{IH}^q_c(F;\vg|F)),
\end{equation*}
where $F$ is the cone on the link knot pair $(S^k,\ell)$.  Since $\Sigma$ is
compact, this is the same as $ H^p(\Sigma; \mc{IH}^q_c(F;\vg|F))$. We first
show that each of these is a $\Gamma$-torsion module whose 
associated polynomial is a product of prime factors of $\xi_{k+1-q}$.
 
Since $\mc{IH}^q_c(F;\vg|F)$ is a locally constant sheaf,  $
H_c^p(\Sigma; \mc{IH}^q_c(F;\vg|F))$ is isomorphic to the
classical singular cohomology with coefficients in a local system by \cite[p.
179-80]{BR2},
and, by
\cite[\S 3.4.i]{Ha}, this is equivalent to the simplicial cohomology with local
coefficients. Since the fiber of $\mc{IH}^q_c(F;\vg|F)$ is $IH^q_c(F;\vg|F)$
and the number of simplices of $\Sigma$ is finite in each dimension, the
classical simplicial cohomology theory with local coefficients \cite{Ste} tells
us that the $i$th dimensional cochain module $C^i(\Sigma;
\mc{IH}^q_c(F;\vg|F))$ is isomorphic to a direct sum of a finite number of
copies of $IH^q_c(F;\vg|F)$. Once again, $F$ is the cone on the link knot pair,
$(S^k,\ell)$, and therefore, as we have noted several times already, each of
the modules $IH^q_c(F;\vg|F)=IH_{k+1-q}(F;\vg|F)$ is isomorphic either to
zero or to the intersection homology group in the same dimension of the link
pair $(S^k,\ell)$,
which is isomorphic to the usual Alexander module of $\ell$. Therefore, 
$C^i(\Sigma; \mc{IH}^q_c(F;\vg|F))$ is the direct sum of a finite number of
torsion $\Gamma$-modules whose associated polynomials are trivial or $\xi_{k+1-q}$,
and thus
its
polynomial is $1$ or a power of $\xi_{k+1-q}$. Since the cohomology modules $
H_c^p(\Sigma; \mc{IH}^q_c(F;\vg|F))$ are quotients of submodules of the
cochain modules, we see that the prime factors of the polynomials associated
to $
H_c^p(\Sigma;
\mc{IH}^q_c(F;\vg|F))$ must divide $\xi_{k+1-q}$ according to the discussion
preceding the theorem.

Now we turn the crank of the spectral sequence. Each of the $E^{p,q}_r$ terms
is
the quotient of a submodule of the $E^{p,q}_{r-1}$ term and hence, by
induction,
each has the property that the prime factors of its associated polynomial must
divide $\xi_{k+1-q}$.
Since this is a bounded first
quadrant spectral sequence, each term converges in a finite number of steps, and
$E_{\infty}^{p,q}$ is a $\Gamma$-torsion module such that the prime factors
of its associated polynomial must divide $\xi_{k+1-q}$.

By spectral sequence theory (see \cite{Mc}) and Proposition \ref{P: SS}, 
\begin{equation*}
E_{\infty}^{p,q}\cong F^p IH^{p+q}_c (N(\Sigma);\vg)/ F^{p+1} IH^{p+q}_c (N(\Sigma);\vg)
\end{equation*}
where the modules  $ F^p IH^{p+q}_c (N(\Sigma);\vg)$ form an ascending
bounded
module filtration of
$ IH^{p+q}_c
(N(\Sigma);\vg)$.  For simplicity, following McCleary \cite{Mc}, let us set $A=IH^{*}_c
(N(\Sigma);\vg)$ as a graded module which is filtered by $F^p A$ and set
$E_0^p(A)=F^pA/F^{p+1}A$.  Then, for some $N$, we have
\begin{equation*}
{0}\subset F^N A \subset F^{N-1}A \subset \cdots \subset F^1 A\subset F^0 A\subset
F^{-1}A=A.
\end{equation*}
This yields the series of short exact sequences
\begin{align}\label{E: filter}
\begin{CD}
0 &@>>>& F^N A &@>\cong>>&E_0^N(A)& @>>>& 0\\ 
0 &@>>>& F^N A &@>>>& F^{N-1}A &@>>>&E_0^{N-1}(A) &@>>>& 0\\ 
&&&&&&\vdots\\ 
0 &@>>>& F^k A &@>>>& F^{k-1} A &@>>>&E_0^{k-1}(A) &@>>>& 0\\ 
&&&&&&\vdots\\ 
0 &@>>> &F^1 A &@>>>& F^{0}A &@>>>&E_0^{0}(A) &@>>>& 0\\ 
0 &@>>>& F^0 A &@>>>& A &@>>>&E_0^{-1}(A) &@>>>& 0.\\ 
\end{CD}
\end{align}

Let us see what happens at the $j$th grade of these graded modules. For
clarity, we will
indicate the grade with a superscript following the argument. For any $p$,
\begin{align*}
E_0^p(A)^j&=(F^pA/F^{p+1}A)^j\\
&=F^pA^j/F^{p+1}A^j\\
&=F^pA^{p+j-p}/F^{p+1}A^{p+j-p}\\
&=E_{\infty}^{p,j-p}.
\end{align*}
We know that each of the prime factors of the polynomial of this module must be
a prime factor
of $\xi_{k+1-(j-p)}$. 
Further, from the construction of the spectral sequence and consideration of dimensions of spaces, we know that the
$E_{\infty}^{p,j-p}$ are non-trivial only if $0\leq p\leq n-k-1$ and $0\leq
j-p\leq k+1$. Hence, as $p$ varies, the only prime factors under consideration
are those of $\xi_{k+1-(j-p)}$ in this range, i.e they are the only possible
prime factors of the $E_0^p(A)^j$,
collectively in $p$ (but within the fixed grade $j$).

Now, by induction down the above list of short exact sequences
\eqref{E: filter} and
their
corresponding polynomial sequences, we can conclude that $F^NA^j$, and
subsequently $F^{N-1}A^j$, $F^{N-2}A^j$,$\dots$, $F^0A^j$, and $A^j$, have the
property of being torsion modules whose polynomials are products of
polynomials whose prime factors are all factors of one of the $\xi_{k+1-s}$,
where $s$ must be chosen in the range $0\leq j-s\leq n-k-1$ and
$0\leq s\leq k+1$.
Since $IH^{j}_c
(N(\Sigma);\vg)$ is the submodule of $A$ corresponding to the $j$th grade, it
too has this property.

Lastly, to draw our conclusions concerning 
$IH_{i}(N(\Sigma);\vg)$, we need only apply the above discussion to
$IH^{n-i}_c (N(\Sigma);\vg)$. Then the relevant factors are those of
$\xi_{k+1-s}$ for $0\leq n-i-s\leq n-k-1$ and
$0\leq s\leq k+1$. Reindexing, these are the polynomials $\xi_{s}$, $0\leq
i-s\leq n-k-1$ and $0\leq s\leq k+1$. We can then strengthen this slightly
by recalling that the $\xi_{s}$ must be
similar to $1$ for $s\geq k-1$, as these are the Alexander polynomials of a
knot $S^{k-2}\subset S^k$. Furthermore, since $\xi_0\sim t-1$ and
$I\lambda_i(1)\neq 0$, we can conclude that the only prime factors of
$IH_{i}(N(\Sigma);\vg)$ which can also divide $I\lambda_i^{\bar p}$ are those
which divide at least one of the $\xi_{s}$, $0\leq
i-s\leq n-k-1$ and $0< s< k-1$.

Now, turning to the term $IH_{i-1}^{\bar p} (X(\Sigma)\cap N(\Sigma);\vg)$, the
arguments are the same, as $X(\Sigma)\cap N(\Sigma)$ will also be a bundle with base
$\Sigma$ and fiber $(S^k,\ell)\times \R$. Thus the intersection homology
groups of the fiber will be the ordinary Alexander modules, and the spectral
sequence
argument will be identical except that we must replace each occurrence of $i$ by
$i-1$. Thus the prime factors which can divide the polynomial of $IH_{i-1}^{\bar p}
(X(\Sigma)\cap N(\Sigma);\vg)$ are those which divide at least one of the $\xi_{s}$,
$0\leq i-1-s\leq n-k-1$ and $0< s< k-1$. Putting these results together yields the
conclusion of the theorem.
\end{proof}

Although it may seem that we have imposed rather weak restrictions on where we
allowed ourselves to look for possible prime factors in the proof, these
results are perhaps the best that one can hope for in large generality. For
example, if the perversities are sufficiently large, then most of the
intersection homology modules of the fibers $IH_i^{\bar p}(c(S^k);\Gamma)$ will
be trivial, and the map $i_*$ will closely approximate the inclusion map
$IH_{i-1}^{\bar p} (X(\Sigma)\cap N(\Sigma);\vg)\to IH_{i-1}^{\bar p}
(X(\Sigma);\vg)$. In this case, which prime factors of the $\xi_s$ divide the
polynomial of the kernel of $i_*$ will depend entirely on the homological
properties
of this inclusion map. If $i_*$ is trivial, then all of the prime factors of
the $\xi_s$ which appear in polynomials of $IH_{i-1}^{\bar p} (X(\Sigma)\cap
N(\Sigma);\vg)$ will appear in $I\lambda^{\bar p}_i(t)$. But which of these
occur in general of course relies heavily on the geometry of the bundle and the ensuing
spectral sequence. For a trivial bundle, they might all occur (see Proposition
\ref{P: man triv1}). Nonetheless, there are some more specific conclusions that
can be drawn in certain situations:

\begin{theorem}\label{T: divisors}
We continue to assume the hypotheses of Theorem \ref{T: at one}. Suppose
$\gamma$ is
a prime element
of $\Gamma$ which does not divide
$\lambda_i(t)$. Suppose $\gamma|\xi_s$ only if $s<k-\bar p(k+1)$. Then
$\gamma\nmid I\lambda_i^{\bar
p}$.
\end{theorem}
\begin{proof}
We will show that the $\gamma$-primary summand of $IH_i^{\bar p} (X(\Sigma)\cap
N(\Sigma);\vg)$ maps isomorphically to the $\gamma$-primary summand of
$IH_i^{\bar p} (N(\Sigma);\vg)$ under the map of the Mayer-Vietoris sequence
\eqref{E: MVM} and similarly for dimension $i-1$. Then the result will follow
from the polynomial sequence associated to the Mayer-Vietoris sequence. In
fact, we can split the long
exact module sequence into the direct sum of an exact sequences involving
the $\gamma$-primary summands of the modules and a sequence involving the other
summands (see the remarks following Corollary \ref{C: sequence splitting} in
Section \ref{S: poly alg}), and we can study the associated
exact polynomial sequences of each module exact sequence. For the latter
sequence
(from which all $\gamma$-primary terms have been removed), the associated
polynomials have no $\gamma$ factors. On the other hand, if the map $IH_i^{\bar
p} (X(\Sigma)\cap N(\Sigma);\vg)\to IH_i^{\bar p} (N(\Sigma);\vg)$ is an
isomorphism in the sequence of $\gamma$-primary summands and $H_i(X(K);\vg)$
has no $\gamma$-primary summand, then the $\gamma$-primary summand of
$IH_i^{\bar p}(S^n;\vg)$ must be $0$ with associated polynomial $1$. Therefore
the total polynomial associated to $IH_i^{\bar p}(S^n;\vg)$, which is the
product of the polynomials of its summands, can have no $\gamma$ factors.

The remainder of the proof will entail a spectral sequence argument in which
we keep special track of only the $\gamma$-primary terms. In particular, notice
that just as in the case of an exact sequence of torsion modules, we can
``split off'' the direct summand corresponding to any given $p$-primary summand
in a spectral sequence of torsion modules. For example, consider the
$p$-primary
summands of each torsion module $E_2^{p,q}$ of a spectral sequence
corresponding to a bounded filtration. All non-trivial maps out of or into each
term must go to or come from a like summand of another term (see Section
\ref{S: poly alg}). When we turn the crank to go to the $E_3$ stage, the
$p$-primary summand of the homology of each $E_2$ term is determined entirely
by
``$p$-primary summands of the $E_2$ layer'' as a $\td p$-primary module cannot
arise as the quotient (or quotient of a submodule of) a $p$-primary module if
$\td p\nsim p$. This follows  by considering short exact
polynomial
sequences. Hence to determine the $p$-primary summands at the $E_3$ stage, we
are free to ignore any non-$p$-primary summands at the $E_2$ stage. Repeating
this argument at each stage, $E_r$, we see that we can ``peel off' the
$p$-primary direct summands of the entire spectral sequence, all the way up to
$E_{\infty}$. Similarly, for the filtration of the abutment, $A$, of the
spectral sequence, it can be seen from diagram \eqref{E: filter} that
the $p$-primary components of $A$ depend only on the $p$-primary components of
the $E_{\infty}\cong E_0(A)$ modulo the usual extension problems. The point of
this argument is that if we care only about some particular $p$-primary
component of the abutment, we can ignore all of the non-$p$-primary terms along
the way, and, in particular, if a map of spectral sequences induces an
isomorphism on $p$-primary components in the $E_2$ term then it will induce 
isomorphisms of the  $p$-primary terms of their abutments.

Consider now the map $i_*: IH_i^{\bar p} (X(\Sigma)\cap N(\Sigma);\vg) \to
IH_i^{\bar p} (N(\Sigma);\vg)$ of the Mayer-Vietoris sequence.  It will be
perhaps more convenient to study the equivalent sheaf theoretic cohomology form
of this sequence. In general, consider a stratified pseudomanifold, $X$, and
let $\mc{IC}^*_U$ denote the extension by zero of the restriction of
$\mc{IC}^*(X)$ to the open set $U$. If we take open sets $U_1$ and $U_2$ with
$U_1\cap U_2=U$ and $U_1 \cup U_2=X$, then we have an exact sequence of
differential sheaves
\begin{equation*}
\begin{CD}
0 @>>> \mc{IC}^*_U @>>> \mc{IC}^*_{U_1}\oplus  \mc{IC}^*_{U_2} @>>>
\mc{IC}^*(X) @>>>0
\end{CD}
\end{equation*}
induced by inclusions, as shown in \cite[\S II.13]{BR2} for a single
sheaf. Since these
sheaves are each soft and hence c-soft (a property preserved by the restriction and
extension by zero \cite[II.9.13]{BR2}), we obtain an exact sequence of chain
modules
\begin{equation*}
\begin{CD}
0 @>>> \Gamma_{c}(\mc{IC}^*_U) @>>>\Gamma_{c}( \mc{IC}^*_{U_1})\oplus
\Gamma_{c}(\mc{IC}^*_{U_2}) @>>> \Gamma_{c}(\mc{IC}^*(X) )@>>>0.
\end{CD}
\end{equation*}
We can now apply homology and use the obvious identification $\Gamma_{c}( \mc{IC}^i_{U})=\Gamma_{c|U}(\mc{IC}^i|U)$ to obtain the long exact hypercohomology sequence

{\small
\begin{equation*}
\to \bb{H}_c^i(\mc{IC}^*(U)) \to \bb{H}_c^i(\mc{IC}^*(U_1))\oplus  \bb{H}_c^i
(\mc{IC}^*(U_2)) \to 
\bb{H}_c^i(\mc{IC}^*(X))\to \bb{H}_c^{i+1}(\mc{IC}^*(U))\to
\end{equation*}}
or, equivalently,
\begin{equation*}
\to IH_c^i(U) \to IH_c^i(U_1) \oplus IH_c^i(U_2)  \to IH_c^i(X)  \to
IH_c^{i+1}(U) \to.
\end{equation*}

In our current situation, with the proper choices of open subsets, this gives us the sequence
\begin{equation*}
\to IH^i_c (X(\Sigma)\cap N(\Sigma);\vg) \overset{i_*}{\to}IH^i_c
(X(\Sigma);\vg)\oplus IH^i_c (N(\Sigma);\vg) \to IH^i_c
(S^n;\vg) \to,
\end{equation*}
which is equivalent to our original Mayer-Vietoris sequence. Note that the
increasing indices of the sequence are offset by the fact that
$IH_i=IH_c^{n-i}$.

With this notation, the map $i_*: IH_i^{\bar p} (X(\Sigma)\cap N(\Sigma);\vg) \to
IH_i^{\bar p}
(N(\Sigma);\vg)$ becomes 
\begin{equation*}
i_*: IH^{n-i}_c (X(\Sigma)\cap N(\Sigma);\vg) \to
IH^{n-i}_c(N(\Sigma);\vg)
\end{equation*}
 induced by the inclusion of sheaves $j:\mc{IC}^*_{
X(\Sigma)\cap N(\Sigma)}\into \mc{IC}^*_{N(\Sigma)}$. This induces a map
$\pi_c(j):\pi_c \mc{IC}^*_{ X(\Sigma)\cap N(\Sigma)}\to
\pi_c\mc{IC}^*_{N(\Sigma)}$ and hence a map of the spectral sequences which can
be used to compute $IH_i^{\bar p} (X(\Sigma)\cap N(\Sigma);\vg)$ and
$IH_i^{\bar p} (N(\Sigma);\vg)$ by Proposition \ref{P: SS} (note that $
X(\Sigma)\cap N(\Sigma)$ is a sub-fiber bundle of $N(\Sigma)$, and it is easy
to check that the proposition applies). We will show that, given the hypotheses
of the theorem, there is an isomorphism of the $\gamma$-primary summands of the
$E_2$ terms of the spectral sequences in such a range as to induce an
isomorphism of the $\gamma$-primary terms of the abutments in dimension $n-i$.
This will complete the proof of the theorem.

We will use the letters $E$, $F$, etc. when referring to the sequence for $N(\Sigma)$ in order to maintain consistency with the above. For $ X(\Sigma)\cap
N(\Sigma)$, we will use $\bar E$, $\bar F$, etc.
As we computed above, 
\begin{equation*}
E_2^{p,q}=H_c^p(\Sigma; \mc{IH}^q_c(F;\vg|F)),
\end{equation*}
where $F$ is the cone on the link knot pair $(S^k,\ell)$. By the formula for
the intersection homology of a cone, these will all be trivial modules for
$0\neq k+1-q\geq k-\bar p(k+1)$.  By assumption, no non-trivial 
$\gamma$-primary terms
can occur in $\bar E_2^{p,q}=H_c^p(\Sigma; \mc{IH}^q_c(S^k;\vg))$ in this
range (i.e.
for $0\neq k+1-q\geq k-\bar p(k+1)$), due to the restrictions on
$IH_*(S^k;\vg)\cong H_*(S^k-\ell;\vg)$.  Thus, restricting to $\gamma$-primary
summands, it is necessary only to show that the maps of the spectral sequence
\begin{equation*}
\bar E_2^{p,q}=H^p(\Sigma; \mc{IH}^{q}_c(S^k\times \R;\vg))\to
H^p(\Sigma; \mc{IH}^{q}_c(c(S^k);\vg))=E_2^{p,q} 
\end{equation*}
induce isomorphisms on their $\gamma$-primary summands in the range $k+1-q<k-\bar p(k+1)$. In
fact, we shall show that these are isomorphisms of the entire modules. This will imply that $i_*$ induces an isomorphism of the $\gamma$-primary summands of the spectral sequences.

Now, once again, the map of the spectral sequence is induced by the map
$\pi_c(j):\pi_c \mc{IC}^*_{
X(\Sigma)\cap N(\Sigma)}\to \pi_c\mc{IC}^*_{N(\Sigma)}$ of c-soft sheaves. 
Thus, if $\mc H^*(L^*)$ represents the derived cohomology sheaf of the
differential sheaf $L^*$, the induced map on the $E_2$ terms is then
\begin{equation*}
H_{c}^p(\Sigma;\mc{H}^q(\pi_c \mc{IC}^*_{ X(\Sigma)\cap
N(\Sigma)}))\to H_{c}^p(\Sigma;\mc{H}^q(\pi_c \mc{IC}^*_{N(\Sigma)})),
\end{equation*}
which comes from the functorial application of the homology functor to the sheaf map 
\begin{equation}\label{E: sheaf map}
\mc H^q(\pi_c(j)):\mc{H}^q(\pi_c \mc{IC}^*_{ X(\Sigma)\cap
N(\Sigma)})\to \mc{H}^q(\pi_c \mc{IC}^*_{N(\Sigma)}).
\end{equation}
As in Proposition \ref{P: SS}, $\mc{H}^q(\pi_c \mc{IC}^*_{ X(\Sigma)\cap
N(\Sigma)})$ and $\mc{H}^q(\pi_c \mc{IC}^*_{N(\Sigma)})$ are
the
sheaves generated by the presheaves
\begin{align*}
U\subset \Sigma &\to H^*(\Gamma_{c\cap\pi^{-1}(U)}( \mc{IC}^*_{N(\Sigma)\cap X(\Sigma)}|\pi^{-1}(U)))\\
&=\bb H^*_{c\cap\pi^{-1}(U)}(\pi^{-1}(U); \mc{IC}_{N(\Sigma)\cap
X(\Sigma)}^*|\pi^{-1}(U))
\end{align*}
and
\begin{align*}
U\subset \Sigma &\to H^*(\Gamma_{c\cap\pi^{-1}(U)}( \mc{IC}^*|\pi^{-1}(U)))\\
&=\bb{H}^*_{c\cap\pi^{-1}(U)}(\pi^{-1}(U); \mc{IC}^*|\pi^{-1}(U)),
\end{align*}
respectively (note again that the extension by zero of the restriction of a $\Phi$-soft
sheave is
$\Phi$-soft \cite{BR2}). The inclusion $ \mc{IC}_{N(\Sigma)\cap
X(\Sigma)}^*|\pi^{-1}(U) \into \mc{IC}^*|\pi^{-1}(U)$ thus induces a map of presheaves
\begin{equation*}
\bb H^*_{c\cap\pi^{-1}(U)}(\pi^{-1}(U); \mc{IC}_{N(\Sigma)\cap
X(\Sigma)}^*|\pi^{-1}(U))\to
\bb{H}^*_{c\cap\pi^{-1}(U)}(\pi^{-1}(U); \mc{IC}^*|\pi^{-1}(U)),
\end{equation*}
for $U\subset \Sigma$, which in turn induces the sheaf map \eqref{E: sheaf map}.
We wish to determine what this map is at the stalk $y\in \Sigma$.

In order to determine this map, we will replace the sheaves $\mc{IC}^*$ with
the quasi-isomorphic Deligne
sheaves $\mc P^*$, recalling that the two give us isomorphic hypercohomology. Because
restrictions and extensions by zero are exact functors, the
corresponding maps 
\begin{equation*}
\bb H^*_{c\cap\pi^{-1}(U)}(\pi^{-1}(U); \mc{P}_{N(\Sigma)\cap
X(\Sigma)}^*|\pi^{-1}(U))\to
\bb{H}^*_{c\cap\pi^{-1}(U)}(\pi^{-1}(U); \mc{P}^*|\pi^{-1}(U))
\end{equation*}
are induced by the corresponding inclusions  $ \mc{P}_{N(\Sigma)\cap
X(\Sigma)}^*|\pi^{-1}(U) \into \mc{P}^*|\pi^{-1}(U)$.
By Lemma \ref{L: restrict}, we may obtain the commutative diagram:

{\small\begin{equation}\label{D: sheaf diagram}
\begin{diagram}
\dlim_{y\in U} \bb H^*_{c\cap\pi^{-1}(U)}(\pi^{-1}(U); \mc{P}_{N(\Sigma)\cap
X(\Sigma)}^*|\pi^{-1}(U))& \rTo& \dlim_{y\in U} \bb 
H^*_{c\cap\pi^{-1}(U)}(\pi^{-1}(U);
\mc{P}^*|\pi^{-1}(U))\\ 
\dTo_{\bar{r}_y^*}&&\dTo_{r_y^*}\\
\bb H^*_{c\cap\pi^{-1}(y)}(\pi^{-1}(y); \mc{P}_{N(\Sigma)\cap
X(\Sigma)}^*|\pi^{-1}(y))&\rTo&\bb H^*_{c\cap
\pi^{-1}(y)}(\pi^{-1}(y);\mc{P}^*|\pi^{-1}(y)),
\end{diagram}
\end{equation}}
in which the vertical maps are isomorphisms and the bottom map is induced by
the inclusion of the
restriction, i.e. by the inclusion of sections  $\Gamma_{c\cap\pi^{-1}(y)}(
\mc{P}_{N(\Sigma)\cap
X(\Sigma)}^*|\pi^{-1}(y))\to \Gamma_{c\cap
\pi^{-1}(y)}(\pi^{-1}(y);\mc{P}^*|\pi^{-1}(y))$. The commutativity is clear at the sheaf
level of sections and is maintained upon applying the hypercohomology and
direct limit functors.

Let $\bar{\pi}: N(\Sigma)\cap X(\Sigma)\to
\Sigma$ be the projection $\pi|N(\Sigma)\cap X(\Sigma)$.
Using the computations of  \eqref{E: restrict} in Proposition \ref{P: SS} and
the
fact that restrictions commute with extensions by zero, the bottom
map of diagram \eqref{D: sheaf diagram} is 
\begin{equation*}
\begin{CD}
\bb H^*_{c\cap\pi^{-1}(y)}(\pi^{-1}(y); \mc{P}^*(\pi^{-1}(y))_{\bar{\pi}^{-1}(y)})
@>>>\bb H^*_{c\cap \pi^{-1}(y)}(\pi^{-1}(y);\mc{P}^*(\pi^{-1}(y))).
\end{CD}
\end{equation*}  
Here $\mc P^*(\pi^{-1}(y))$ is the Deligne sheaf on $\pi^{-1}(y)$ and the extension by zero in
the first term is extension only to the rest of $\pi^{-1}(y)$. By reversing the arguments of
the preceding paragraph in this context, this is the map 
\begin{equation*}
\begin{CD}
\bb H^*_{c\cap\pi^{-1}(y)}(\pi^{-1}(y); \mc{IC}^*(\pi^{-1}(y))_{\bar{\pi}^{-1}(y)})
@>>>\bb H^*_{c\cap \pi^{-1}(y)}(\pi^{-1}(y);\mc{IC}^*(\pi^{-1}(y)))       
\end{CD}
\end{equation*}
induced by the inclusion $\mc{IC}^*(\pi^{-1}(y))_{\bar{\pi}^{-1}(y)}\into
\mc{IC}^*(\pi^{-1}(y))$.
By the $c\cap \pi^{-1}(y)$-softness of these sheaves, this is the map \begin{equation*}
\begin{CD}
H^*(\Gamma_{c\cap\pi^{-1}(y)}( \mc{IC}^*(\pi^{-1}(y))_{\bar{\pi}^{-1}(y)}))
@>>> H^*(\Gamma_{c\cap \pi^{-1}(y)}(\mc{IC}^*(\pi^{-1}(y))))
\end{CD}
\end{equation*}
induced by the inclusion $\mc{IC}^*(\pi^{-1}(y))_{\bar{\pi}^{-1}(y)}\to \mc{IC}^*(\pi^{-1}(y)).$

But,
\begin{align*}
\Gamma_{c\cap\pi^{-1}(y)}( \mc{IC}^*(\pi^{-1}(y))_{\bar{\pi}^{-1}(y)})
&=\Gamma_{c\cap\pi^{-1}(y)|\bar{\pi}^{-1}(y)}(\mc{IC}^*(\pi^{-1}(y))|\bar{\pi}^{-1}(y))\\
&=\Gamma_c(\mc{IC}^*(\bar{\pi}^{-1}(y))
\end{align*}
because the restriction of the intersection chain sheaf to an open subset is the
intersection chain
sheaf of the subset and by the identity $(c\cap\pi^{-1}(y))|\bar{\pi}^{-1}(y)=c$ on
$\bar{\pi}^{-1}(y)$, which is easily verified. Therefore, the relevant inclusion of sheaves
induces the inclusion of chains $\Gamma_c(\mc{IC}^*(\bar{\pi}^{-1}(y)))\into
\Gamma_{c\cap
\pi^{-1}(y)}(\mc{IC}^*(\pi^{-1}(y)))=\Gamma_{c}(\mc{IC}^*(\pi^{-1}(y)))$,
the last
equality
because $\pi^{-1}(y)$ is a close subspace of $\pi^{-1}(U)$ 
(or of $N(\Sigma)$ in general). 
But this  is the familiar inclusion which induces the map from the link intersection homology to that of its cone in
simplicial intersection homology. 
Furthermore, we know that the induced map is an intersection homology
isomorphism
on this summand in the dimension range under consideration. Hence, the locally
constant sheaves $\mc{H}^q(\pi_c \mc{IC}^*_{ X(\Sigma)\cap
N(\Sigma)})$ and $\mc{H}^q(\pi_c \mc{IC}^*_{N(\Sigma)})$ have identical stalk
components which are identified isomorphically by the map induced by the sheaf
inclusion. Thus the inclusion induces a sheaf isomorphism.

Returning then to map between the $E_2$ terms of the spectral sequence.  It is, once again,
the map 
\begin{equation*} 
H_{c}^p(\Sigma;\mc{H}^q(\pi_c \mc{IC}^*_{ X(\Sigma)\cap
N(\Sigma)}))\to H_{c}^p(\Sigma;\mc{H}^q(\pi_c \mc{IC}^*_{N(\Sigma)})) 
\end{equation*}
induced by the sheaf inclusion. But we have just calculated that the map of 
coefficient sheaves is an isomorphism. Therefore, if, as in the proof of the
last
theorem, we think of these modules as given by simplicial homology with local
coefficients, the inclusion map induces an isomorphism on the cochain modules.
This in turn induces an isomorphism of the cohomology modules and hence of the
relevant $E_2$ terms of the spectral sequences, as we were to show.

Analogous consideration apply for the intersection homology in dimension $i-1$, the
slight fluctuation in allowable ranges accounted for by the hypotheses of the theorem.
\end{proof}

As we mentioned at the beginning of this section, the assumptions in the previous two theorems that the singular set $\Sigma$ has a fiber bundle neighborhood are unnecessary. This is due to the  following theorem from \cite{GBF3} and \cite{GBF}:

\begin{theorem}\label{P: Neigh IH}
Let $X$ be a finite-dimensional stratified pseudomanifold with
locally 
finite triangulation and filtration $\emptyset=X_{-1}\subset X_0\subset\cdots\subset X^n=X$ such that $X_i=\emptyset$ for $i<k$. Let $N=N(X_k)$ be an
open
regular neighborhood of $X_k$, and let $L$ be the link of the stratum $X_k$ (if $X_k$ is not connected, then we can treat each component separately and each component will have its own link). Then, for any fixed perversity $\bar p$ and local coefficient system $\mc{G}$ defined on $X-X_{n-2}$,
there are homological-type spectral sequences $\bar E^r_{p,q}$ and
$E^r_{p,q}$ that abut (up to isomorphism)
to
$IH^{\bar p}_i(N-X_k;\mc{G})$ and $IH^{\bar p}_i(N;\mc{G})$ with respective
$E^2$
terms
\begin{align*}
\bar  E^2_{p,q}=H_p(X_k; \mc{IH}^{\bar p}_q(L;\mc{G}|L)) &&
E^2_{p,q}=H_p(X_k; \mc{IH}^{\bar p}_q(cL;\mc{G}|cL))
\end{align*}
($cL=$ the open cone on $L$),
where $\mc{IH}^{\bar p}_q(L;\mc{G}|L)$ and $\mc{IH}^{\bar p}_q(cL;\mc{G}|cL)$ are local
coefficient systems with respective stalks $IH^{\bar p}_q (L;\mc{G}|L)$ and $IH^{\bar p}_q
(cL;\mc{G}|cL)$ . Furthermore, the map
$i_*: IH^{\bar p}_i(N-X_k;\mc{G})\to
IH^{\bar p}_i(N;\mc{G})$ induced by inclusion induces a map of spectral sequences which
on the the $E^2$ terms is determined by the coefficient homomorphism
$\mc{IH}^{\bar p}_q(L;\mc{G}|L)\to \mc{IH}^{\bar p}_q(cL;\mc{G}|cL)$ given by the map on the
stalk intersection homology modules induced by the inclusion $L\into cL$.
\end{theorem}

For our study of a knot $S^{n-2}\subset S^n$ with only one singular
stratum, recall that that we are stratifying $S^n$ as $S^n\supset
S^{n-2}\supset \Sigma$, where $\Sigma$ is the singular set of the embedding. So far in this section, we have deduced several results
concerning the intersection Alexander
polynomials of the knot by studying the long exact Mayer-Vietoris sequence corresponding
to the
pair $(N(\Sigma), X(\Sigma))$, where $N(\Sigma)$ was an open regular neighborhood of
$\Sigma$ and $X(\Sigma)$ was the complement a closed regular neighborhood of
$\Sigma$ contained in $N(\Sigma)$. Assuming that
$N(\Sigma)$ had the structure of 
a fiber bundle, the arguments then involved the use of spectral sequences to compute
and compare the intersection homologies of $N(\Sigma)$ and $N(\Sigma)\cap
X(\Sigma)$. However, we can apply Theorem \ref{P: Neigh IH} to show that the assumption
that $N(\Sigma)$ be a bundle is unnecessary. In fact, the spectral sequences of Theorem \ref{P: Neigh IH} simply take the role of the sheaf theoretic spectral sequences we have been using. It only remains to prove the following
lemma and its corollary which address the compatibility of the sets mentioned in the theorem with those of the previous arguments in this section.

\begin{lemma}\label{L: IH collapse}
Suppose that $Z$ is a subpolyhedron of the filtered polyhedron $Y$ and that $\bar N(Z)$ is a
filtered regular
neighborhood of $Z$ (such a neighborhood always exists; see \cite[p. 26]{Sto}). Let $N(Z)$ be
the
interior of $\bar N(Z)$, let $\bar N'(Z)$ be a closed regular neighborhood of $Y$ in $N(Z)$ (and
hence also in $Y$),  and let $X(Z)=Y-\bar N'(Z)$. Then $IH_*^{\bar p}(Y-Z ,X(Z) ;\mc{G})=0$.
\end{lemma}
\begin{proof}

By its definition, $X(Z)$ is the complement in $Y$ of a closed regular
neighborhood of $Z$, $\bar N'(Z)$, which is contained in $N(Z)$. Suppose that we are
given a singular intersection chain $\sigma:\Delta^i\to Y-Z$ (covered by
local coefficients) which is a relative cycle. Since the image of $\sigma$ is
compact, we can find another closed regular neighborhood of $Z$, $\bar N''(Z)\subset N'(Z)$,
such that Im$(\sigma)\subset Z-\bar N''(Z)$. Hence, we have the inclusions $\bar
N''\subset
\bar N'\subset \bar N$ in which each closed set is contained in the interior of the
next. By \cite[Prop. 1.5]{Sto}, cl$(\bar N(Z)-\bar N''(Z))\cong_{PL}I\times$fr$(N''(Z))$, and
this product respects the stratification. Thus, there exists a stratum-preserving
deformation retraction $Y-\bar N''(Z)\to Y-N(Z)$ by retracting along the product
lines. Furthermore, the homotopy which gives the retraction can be applied to give a
relative homology of $\sigma$ with a chain lying in $Y-\bar N'(Z)$ by the standard prism
process (see \cite{GBF3} for proofs that this is allowable when using intersection homology). Note that while the sides of the prism may
contribute to the boundary, those chains all lie in $Y-\bar N'(Z)$ by the
assumption that $\sigma$ is a relative cycle and the fact that our homotopy
equivalence takes $Y-\bar N'(Z)$ into itself. Since we can apply this argument to any
cycle in $IC_i^{\bar p}(Y-Z, X(Z);\mc{G})$, we conclude that $IH_i^{\bar
p}(Y-Z, X(Z);\mc{G})=0$.
\end{proof}

\begin{corollary}
Consider the commutative diagram of Mayer-Vietoris sequences induced by
the inclusions $X(\Sigma)\into S^n-\Sigma$:

{\small\begin{equation*}
\begin{diagram}
\,&\rTo& IH_i^{\bar p} (X(\Sigma)\cap N(\Sigma);\vg) & \rTo & IH_i^{\bar
p}
(X(\Sigma);\vg)\oplus IH_i^{\bar p} (N(\Sigma);\vg) &\rTo & IH_i^{\bar p}
(S^n;\vg) &\rTo&\\
&& \dTo &  & \dTo & & \dTo &\\
\,&\rTo& IH_i^{\bar p} (N(\Sigma)-\Sigma;\vg) & \rTo & IH_i^{\bar p}
(S^n-\Sigma;\vg)\oplus IH_i^{\bar p} (N(\Sigma);\vg) &\rTo & IH_i^{\bar p} 
(S^n;\vg) &\rTo &. 
\end{diagram}
\end{equation*}}

\noindent This diagram is an isomorphism of exact sequences.
\end{corollary}
\begin{proof}
The commutativity is obvious, as is the fact that the middle vertical map does not mix
terms. It is
sufficient to prove that the map induced by inclusion $i_*: IH_i^{\bar p} (X(\Sigma);\vg)\to
IH_i^{\bar p} (S^n-\Sigma;\vg)$ is an isomorphism, and the result will follow by the five-lemma. From
the long exact sequence of the pair $(S^n-\Sigma, X(\Sigma))$, it suffices to show that $IH_i^{\bar
p}(S^n-\Sigma, X(\Sigma);\vg)=0$, but this is the content of the lemma.
\end{proof}

We can now repeat, up to isomorphism, any spectral sequence arguments used earlier in this
section in
the context of the Mayer-Vietoris sequence of the pair $(N(\Sigma); S^n-\Sigma)$, but
without assuming any bundle structures. In fact, the bundle assumptions were only
used to obtain the spectral sequences in the first place and to determine the maps 
between them. Of course those arguments using the
cohomological language of sheaves must now be dualized from cohomological to homological indexing, a straightforward chore which we leave to the reader. Furthermore, we should point out  that, according to Theorem \ref{P: Neigh IH},  the map of the spectral sequences is induced by the inclusion $N(\Sigma)-\Sigma \into N(\Sigma)$ and is determined on the $E^2$ terms entirely by 
by the coefficient homomorphism
$\mc{IH}^{\bar p}_q(L;\mc{G}|L)\to \mc{IH}^{\bar p}_q(cL;\mc{G}|cL)$ given by the map on the
stalk intersection homology modules induced by the inclusion $L\into cL$. Hence, this is either a canonical isomorphism or the $0$ map, according to the dimension ranges in the usual cone formula. But this agrees precisely with the computation of the map for sheaf intersection cohomology used in the proof of Theorem \ref{T: divisors}.

Note also that,
by the proof of the preceding lemma, $IH_i^{\bar p}(S^n-\Sigma;\vg)\cong IH_i^{\bar
p}( X(\Sigma);\vg)$, so there is no difficulty added to the computations involving
the polynomials of the knot exterior, i.e. both of these modules are isomorphic to
$H_i(S^n-\Sigma;\vg)$ according to our previous calculations. 

Thus we conclude:

\begin{theorem}
Theorems \ref{T: at one} and \ref{T: divisors}
concerning intersection Alexander polynomials of knots with one singular stratum remain
true
without the assumption of the existence of a fiber bundle neighborhood of the
singular stratum.
\end{theorem}

We conclude this section with a few corollaries:

\begin{corollary}\label{S: degeneracies}
For a knot $K\subset S^n$ with a  manifold singularity of dimension $n-k-1$:  
\begin{enumerate}
\item If $i<k-\bar p(k+1)$, then $I\lambda_i^{\bar p}\sim \lambda_i(t)$.
\item If $\bar p(k+1)\leq 1$ or if $H_i(S^k-\ell;\vg)=0$ for $i\geq j$ and $\bar
p(k+1)\leq
k-j$, then $I\lambda_i^{\bar p}\sim \lambda_i(t)$ for all $i$.
\item If $i\geq n-\bar p(k+1)+1$, then $I\lambda_i^{\bar p}\sim \mu_i(t)$.
\end{enumerate}
\end{corollary}
\begin{proof}
\begin{enumerate}
\item For $q<k-\bar p(k+1)$, the map $IH^{\bar p}_q(S^k;\vg)\to IH^{\bar
p}_q(cS^k;\vg)$
induced
by inclusion is an isomorphism \cite{Ki}. Hence the maps $\mc{IH}^{\bar
p}_q(S^k;\vg)\to \mc{IH}^{\bar p}_q(cS^k;\vg)$ and $H_p(\Sigma;\mc{IH}^{\bar  
p}_q(S^k;\vg))\to H_p(\Sigma;\mc{IH}^{\bar p}_q(cS^k;\vg))$ are isomorphisms
induced by inclusion for $q<k-\bar p(k+1)$. But these are the maps of all of the
terms of
the spectral sequences which abut to $IH_i^{\bar p}(N(\Sigma)-\Sigma;\vg)$
and $IH_i^{\bar p}(N(\Sigma);\vg)$. Hence, by spectral sequence theory, the
inclusion induces an isomorphism $IH_i^{\bar p}(N(\Sigma)-\Sigma;\vg) \to
IH_i^{\bar p}(N(\Sigma);\vg)$ for $i<k-\bar p(k+1)$. This implies the corollary by
the Mayer-Vietoris sequence \eqref{E: MVM}.

\item In these cases, the maps  $H_p(\Sigma;\mc{IH}^{\bar p}_q(S^k;\vg))\to
H_p(\Sigma;\mc{IH}^{\bar p}_q(cS^k;\vg))$ are all isomorphism, so we obtain an
isomorphism of   
spectral sequences. Thus  $IH_i^{\bar p}(N(\Sigma)-\Sigma;\vg) \to IH_i^{\bar
p}(N(\Sigma);\vg)$ is an isomorphism for all $i$, and we conclude as above.

\item For this case, we note that if $p+q\geq n-\bar p(k+1)$ then
$H_p(\Sigma;\mc{IH}^{\bar p}_q(cS^k;\vg))=0$. Therefore, $IH_i^{\bar
p}(N(\Sigma);\vg)=0$ for $i\geq n-\bar p(k+1)$, and, in this range, the kernel of
the map $IH_i^{\bar p}(N(\Sigma)-\Sigma;\vg) \to IH_i^{\bar
p}(N(\Sigma);\vg)\oplus IH_i^{\bar p}(S^n-\Sigma ;\vg)$ of the Mayer-Vietoris
sequence is the kernel of the map $H_i(N(\Sigma)-\Sigma;\vg) \to
H_i(S^n-\Sigma;\vg)$, using the fact that the intersection homology groups coincide
with the ordinary homology groups of the complement for locally-flat
embeddings (see the proofs in Section \ref{S: ih and ord}). Thus for $i\geq
n-\bar p(k+1)$ the polynomial sequence associated to the intersection homology
Mayer-Vietoris sequence is isomorphic to the long exact polynomial sequence associated
to the
long exact sequence of the knot pair in ordinary homology in two thirds of its  
terms. Since the corresponding exact polynomial sequences argee in two-thirds of
their terms and in the common factors of those terms, they
must agree  in the remaining terms (see Section \ref{S: poly alg}) and
$I\lambda^{\bar p}_i\sim \mu_i(t)$ for
$i\geq n-\bar p(k+1)+1$ .

\end{enumerate}
\end{proof}

\begin{corollary}
For a knot $K\subset S^n$ with manifold singularity of dimension $n-k-1$ and a
superperversity $\bar p$, then:
\begin{enumerate}
\item If  $i>n-1-\bar p(k+1)$, then $I\lambda_{i}^{\bar
p}(t)\sim\lambda_{n-i-1}(t^{-1})$.
\item  If $\bar p(k+1)\geq k-1$ or if $H_i(S^k-\ell;\vg)=0$ for $i\geq j$ and $\bar
        p(k+1)\geq j$, then $I\lambda_{i}^{\bar p}(t)\sim \lambda_{n-i-1}(t^{-1})$ for all $i$.
\item If $i\leq k-\bar p(k+1)-2$, then $I\lambda_{i}^{\bar p}(t)\sim
        \mu_{n-i-1}(t^{-1})$.
\end{enumerate}
\end{corollary}
\begin{proof}
These statements follow from the previous corollary and
superduality.
\end{proof}


\section[Relations between $I\protect\lambda_i$ and the polynomials of the
links for general
non-locally-flat knots]{Relations between {\mv$\protect I\lambda_i$}
and the polynomials of the
links for general
non-locally-flat knots}\label{S: general knots}

In this section, we develop some relations between the intersection
Alexander polynomials of a knot, its ordinary Alexander polynomials,
and the intersection and ordinary Alexander polynomials of
its link knots. For example, we
determine what the possible prime divisors
of the intersection Alexander polynomials can be 
in terms of the prime divisors of these other polynomials. 

Once again, we consider the PL-knot $K$, given by $S^{n-2}\subset S^n$, as
a stratified
pseudomanifold pair. The top skeleton of the filtration is
$K\cong S^{n-2}$,
and the lower skeleta are denoted $\Sigma_{i}$; if there is no stratum 
of dimension $k$, then we let $\Sigma_k=\Sigma_{k-1}$. We do not place any
unnatural limitation on the number of non-empty strata as we did in
previous sections. Also, for consistency of notation, let
$\Sigma_{n-2}=K\cong S^{n-2}$ and let $\Sigma_n=S^n$. We continue to let $\Sigma$, with
no
index, denote the set of points at which the embedding is
non-locally-flat, i.e. the union of all strata of codimension greater than
$2$. We will assume, initially, that $S^n$
is given
a fixed triangulation with each $\Sigma_i$ triangulated as a full subcomplex.
We continue to employ the
local coefficient system $\vg$ defined on $S^n-S^{n-2}$. 

Let $\bar N_0$ be a closed regular neighborhood of
$\Sigma_0$ (or empty if $\Sigma_0$ is), let $W_1=S^n-\bar N_0$, and let
$X_1=\Sigma_1\cap W_1$. Then inductively define $\bar N_i$ to be a closed
regular neighborhood of $X_i$ in $W_i$,
$W_i=S^n-\cup_{j<i}\bar N_{j}$, and $X_i=\Sigma_i\cap
W_i$.
(We assume that all regular neighborhoods are given by derived neighborhoods in sufficiently fine subdivisions so that, in particular, there are no connected components of $N_i\cap \Sigma_j$, $j>i$, which do not also intersect $\Sigma_i$.)
 For consistency, we also
set $W_0=S^n$ and $X_0=\Sigma_0$. 
Notice that $W_i$ is
equal to either $W_{i-1}$ or $W_{i-1}-\bar
N_{i-1}$, depending on whether or not $X_i$ is empty. Each $W_i$ is
open in $S^n$ and in $W_{i-1}$, and $X_i$ is a close subpolyhedron of 
$W_i$. Furthermore, the triangulation of $S^n$ induces triangulations on 
$W_i$ and $X_i$ for each $i$, and $W_i$ is a stratified pseudomanifold
filtered by the restrictions of the $\Sigma_k$. Any distinguished
neighborhood in $W_i$ of a point in $X_i$
is also a distinguished neighborhood in $S^n$, so the
links, $L_{ik}$, of the connected components of $X_i$ are the same as the
links of the connected components of $\Sigma_i$.
Also,
for each $i$, let $\bar N'_i$ be a closed regular neighborhood
of $\bar N_i$ in $W_i$, and let $N'_i$ be the interior of $\bar
N'_i$. Then $\bar N'$ is also a closed regular neighborhood of $X_i$ in
$W_i$, $\bar N_i\subset N'_i$, and $W_i=N'_i\cup (W_i-\bar N_i)$, which is
a union of open sets in $W_i$. Therefore, there are Mayer-Vietoris
sequences
\begin{equation*}
\to IH_j^{\bar p}(N'_i\cap(W_i-\bar N_i);\vg) \to IH_j^{\bar
p}(N'_i;\vg)\oplus  IH_j^{\bar p}(W_i-\bar N_i;\vg) \to  IH_j^{\bar
p}(W_i;\vg) \to.
\end{equation*}

From the definitions, $W_i=W_{i-1}-\bar N_{i-1}$, and $IH_j^{\bar p}(N'_i\cap(W_i-\bar N_i);\vg)  \cong IH_j^{\bar p}(N'_i-X_i;\vg)$
by Lemma \ref{L: IH collapse}. So
this Mayer-Vietoris
sequence is isomorphic to the following
long exact sequence:
\begin{equation}\label{E: W sequence}
\to IH_j^{\bar p}(N'_i-X_i;\vg)  \to IH_j^{\bar
p}(N'_i;\vg)\oplus  IH_j^{\bar p}(W_{i+1};\vg) \to  IH_j^{\bar    
p}(W_i;\vg) \to.
\end{equation}
The map from $IH_j^{\bar p}(N'_i-X_i;\vg) $ to the summand $IH_j^{\bar
p}(N'_i;\vg)$ is the map induced by the inclusion $N'_i-X_i\into N'_i$. 

\begin{lemma}\label{L: IH(W)}
For all $i$, $IH_*^{\bar p}(W_i-W_i\cap \Sigma;\vg)\cong IH_*^{\bar
p}(S^n-\bar N(\Sigma);\vg)\cong H_*(S^n-K;\vg)$.
\end{lemma}
\begin{proof}
The second isomorphism is established in the proof of Proposition
\ref{P: man triv1}. For the first isomorphism, we will show that 
$IH_*^{\bar p}(W_i-W_i\cap \Sigma;\vg)\cong IH_*^{\bar
p}(W_{i+1}-W_{i+1}\cap
\Sigma;\vg)$ for all $i$. This suffices because $W_0=S^n$ and $IH_*^{\bar
p}(S^n-\Sigma;\vg)\cong IH_*^{\bar
p}(S^n-\bar N(\Sigma);\vg)$ by Lemma \ref{L: IH collapse}.

That $IH_*^{\bar p}(W_i-W_i\cap \Sigma;\vg)\cong IH_*^{\bar
p}(W_{i+1}-W_{i+1}\cap
\Sigma;\vg)$ is established by showing that 
$IH_*^{\bar p}(W_i-W_i\cap \Sigma_i, W_{i+1}-W_{i+1}\cap
\Sigma;\vg)=0$, and the demonstration of this  is essentially
the proof of Lemma \ref{L: IH collapse}. Of course if $W_i=W_{i+1}$ the
proof is trivial, so assume otherwise. If $C$ is a relative
cycle representing an element of $IH_*^{\bar p}(W_i-W_i\cap \Sigma_i,
W_{i+1}-W_{i+1}\cap
\Sigma;\vg)$, then $C$ lies in $W_i-W_i\cap \Sigma_i\subset W_i$ and it
boundary lies in $W_{i+1}-W_{i+1}\cap
\Sigma\subset W_{i+1}=W_{i}-\bar N_{i}$. Now, just as in the proof of
Lemma \ref{L: IH collapse}, there is a stratum preserving homotopy which
takes $C$ into $ W_{i+1}$ while keeping $\bd C$ in $W_{i+1}$, and, by the
prism construction (see the proof of Lemma \ref{L: IH collapse}), this induces a
relative homology from $C$ to a chain in $ W_{i+1}$. However, since $\Sigma$ is a
skeleton and the homotopy is stratum preserving, the homotopy on $C$ lies
entirely in $W_i-W_i\cap \Sigma$ and provides a homology in
$W_i-W_i\cap \Sigma$
between
$C$ and a chain in $W_{i+1}-W_{i+1}\cap\Sigma$. Since $C$ was an arbitrary
relative
cycle, we have $IH_*^{\bar p}(W_i-W_i\cap \Sigma_i, W_{i+1}-W_{i+1}\cap
\Sigma;\vg)=0$.
\end{proof}

Since $W_n-W_n\cap \Sigma=W_n$, we can hope to begin with $IH_*^{\bar
p}(W_n;\vg)$, which by the preceding lemma is simply the Alexander module 
$H_*(S^n-K;\vg)$, and determine something about the composition of the
polynomials of the $IH_*^{\bar p}(W_i;\vg)$, $i<n$, by an induction
involving the long exact
sequences \eqref{E: W sequence}. Since $W_0=S^n$, we induct down to the
intersection Alexander module of the knot. The first few steps are trivial
because $W_i-W_i\cap \Sigma=W_i$ for $n\geq i\geq n-3$. After that, the
singular strata of the embedding begin to come in. 
By polynomial algebra, we know that the polynomial of 
$ IH_j^{\bar p}(W_{i};\vg)$ must divide the product of the polynomials
of $ IH_j^{\bar p}(W_{i+1};\vg)$, $ IH_j^{\bar
p}(N'_i;\vg)$, and $IH_{j-1}^{\bar p}(N'_i-X_i;\vg)$. We can know
something about the first by induction. To study the latter two, we note
that $X_i$ is the bottom stratum in $W_i$. So, as we discussed in the 
last section, by Theorem \ref{P: Neigh IH} each is the abutment of a spectral sequence whose 
$E^2_{p,q}$ term is the homology of $X_i$ with a local coefficient system
with fiber isomorphic to the intersection homology of the link, which is
itself a knot pair. More specifically, by Theorem \ref{P: Neigh IH},
for each connected component $X_{i,k}$ of $X_k$ with
link $L_{i,k}$ and regular neighborhood $\bar N_{i,k}$,
there are homological-type spectral sequences $\bar E^r_{p,q}$ and
$E^r_{p,q}$ which abut (up to isomorphism) to
$IH^{\bar p}_*(N_{i,k}-X_{i,k};\vg)$ and $IH^{\bar p}_*(N_{i,k};\vg)$
with
respective
$E^2$
terms
\begin{align*}
\bar  E^2_{p,q}=H_p(X_{i,k}; \mc{IH}_q(L_{i,k};\vg)) &&
E^2_{p,q}=H_p(X_{i,k}; \mc{IH}_q(cL;\vg)),
\end{align*}
where $\mc{IH}_q(L_{i,k};\vg)$ and $\mc{IH}_q(cL_{i,k};\vg)$
are local
coefficient systems with respective stalks $IH_q (L_{i,k};\vg)$ and
$IH_q
(cL_{i,k};\vg)$. Furthermore, the map
$i_*: IH_*(N_{i,k}-X_{i,k};\vg)\to
IH_*(N_{i,k};\vg)$ induced by inclusion induces a map of spectral sequence
which
on the the $E^2_{p,q}$ terms is determined by the coefficient homomorphism
$\mc{IH}_q(L_{i,k};\vg)\to \mc{IH}_q(cL_{i,k};\vg)$ given by the map
on the
fiber intersection homology modules induced by the inclusion $L_{i,k}\into
cL_{i,k}$. Once again, this latter intersection homology map is the
identity map for $q=0,
q<n-i-1-\bar p(n-i)$ and the
the $0$ map for $0\neq q\geq n-i-1-\bar p(n-i)$.

It will be useful to know that $\bar E^2_{p,q}$ and $E^2_{p,q}$ are
finitely generated $\Gamma$-modules. To establish this, we begin by
demonstrating that if $\mc{G}$ is a local coefficient system on $X_i$ and
$U_i=\Sigma_i-\Sigma_{i-1}$, then $H_*(X_i;\mc{G})\cong
H_*(\Sigma_i-\Sigma_{i-1}; \bar{\mc{G}})$, where $\bar {\mc{G}}$ is a
suitable local coefficient system such that $\bar{ 
\mc{G}}|X_i=\mc{G}$. This
can once again be proven inductively by the methods of Lemmas \ref{L: IH
collapse} and \ref{L: IH(W)}. Recall that, by definition,
$X_i=\Sigma_i\cap W_i$ and $W_i$ is equal to either $W_{i-1}$ or
$W_{i-1}-\bar N_{i-1}$. Consider $\Sigma_i\cap W_{i-1}$. If $W_i=W_{i-1}$,
then $\Sigma_i\cap W_{i-1}=X_i$, and clearly $H_*(X_i;\mc{G})=H_*(\Sigma_i\cap
W_{i-1};\mc{G})$. If $W_i=W_{i-1}-\bar N_{i-1}$, then $X_i=\Sigma_i\cap
W_{i-1}-(\Sigma_i\cap\bar N_{i-1})$.  But since $\Sigma_{i-1}\cap
W_{i-1}\subset \bar N_{i-1}$, then $\Sigma_i\cap W_{i-1}-(\Sigma_i\cap\bar
N_{i-1})= (U_i)\cap (W_{i-1}-\bar N_{i-1})=(U_i\cap W_{i-1})-(U_i\cap \bar
N_{i-1})$. Define $ \Sigma_{i,j}= U_i\cap W_{j}$, and note that
$X_i=\Sigma_{i,i}$. We now claim that $H_*(\Sigma_{i,i-1}, X_i;\bar{
\mc{G}})=0$ for suitable $\bar{\mc{G}}$.

We note once again as in the proof of  Lemmas  \ref{L: IH
collapse} and \ref{L: IH(W)}
that, for any compact set in $W_{i-1}-X_{i-1}$, there exists
a stratum-preserving homotopy which retracts this set into $W_{i-1}-\bar
N_{i-1}$. Since the homotopy is stratum-preserving, any compact set in
$U_i\cap W_{i-1}=\Sigma_{i,i-1}$ retracts within $\Sigma_{i,i-1}$ into
$U_i\cap
(W_{i-1}-\bar N_{i-1})$. So in particular, any compact set in
$\Sigma_{i,i-1}$ retracts into $ X_i$, since $\Sigma_{i,i-1}\subset
W_{i-1}-X_{i-1}$ and $U_i\cap (W_{i-1}-\bar N_{i-1})=X_i$.
Now, we can define $\bar
{\mc{G}}$ on $\Sigma_{i,i-1}$. 
Choose a basepoint in $X_i$ and any loop $\gamma$ in
$\Sigma_{i,i-1}$ representing an element of
$\pi_1(\Sigma_{i,i-1})$. Since the image of $\gamma$ is compact, we can
retract the loop into a loop in $X_i=\Sigma_{i,i}$ representing the same
element of
$\pi_1(\Sigma_{i,i-1})$. Thus the action of $\gamma$ on the fiber over the
basepoint is given by the action of the retracted $\gamma$. Hence, we have
determined a local coefficient system on $\Sigma_{i,i-1}$ which clearly
restricts to $\mc{G}$ on $X_i$. The proof that $H_*(\Sigma_{i,i-1},
\Sigma_{i.i};\bar{\mc{G}})=0$ now also proceeds as in the cited lemmas by
retracting relative cycles in $X_i$. This further implies that
$H_*(X_i;\mc{G})=H_*(\Sigma_{i,i};\mc{G})\cong
H_*(\Sigma_{i,i-1};\bar{\mc{G}})$. 

We can now continue by downward induction to show that, for $j\leq i$,
$H_*(\Sigma_{i,j-1},\Sigma_{i,j};\bar{\mc{G}})=0$. If $W_j=W_{j-1}$ this
is
again trivial. Otherwise, we need only note that once again
$\Sigma_{i,j}= U_i\cap W_{j}=U_i\cap (W_{j-1}-\bar
N_{j-1})$, while $\Sigma_{i,j-1}=U_i\cap W_{j-1}$. It now follows by the same methods as
the preceding paragraph that
$H_*(\Sigma_{i,j},\Sigma_{i,j-1};\bar{\mc{G}})=0$ for a similarly chosen
extension of the coefficient system and hence that
$H_*(\Sigma_{i,j};\bar{\mc{G}})\cong H_*(\Sigma_{i,j-1};\bar{
\mc{G}})$. But $\Sigma_{i,0}=U_i\cap W_0=U_i\cap S^n=U_i$. Therefore, 
$H_*(X_i;\mc{G})\cong H_*(U_i;\bar{\mc{G}})$. 

Now, with respect to our initial triangulation of $S^n$ or one of its
derived subdivisions, let $\bar N(\Sigma_{i-1})$ be a closed regular
neighborhood of
$\Sigma_{i-1}$ in $S^n$ with interior $N(\Sigma_{i-1})$. Note that $\bar
N(\Sigma_{i-1})$ and $N(\Sigma_{i-1})$ are stratum-preserving homotopy
equivalent as are $S^n-N(\Sigma_{i-1})$ and $S^n-\bar
N(\Sigma_{i-1})$ (for example, each member of the first  pair 
has a stratum-preserving deformation
retraction to any closed regular neighborhood of $\Sigma_{i-1}$ interior
to $N(\Sigma_{i-1})$, while each member of the latter pair has a 
stratum-preserving deformation retraction to any closed regular
neighborhood of $\Sigma_{i-1}$ which contains $\bar
N(\Sigma_{i-1})$ in its interior). Also, again by the
methods of the preceding
paragraphs, $H_*(U_i, U_i-U_i\cap \bar N(\Sigma_{i-1});\bar{\mc{G}})=0$
since
$(U_i, U_i-U_i\cap \bar N(\Sigma_{i-1}))=U_i\cap (S^n-\Sigma_{i-1},
S^n-\bar N(\Sigma_{i-1}))$. Therefore, $H_*(U_i;\bar{\mc{G}})\cong 
H_*(U_i-U_i\cap \bar N(\Sigma_{i-1});\bar{\mc{G}})\cong H_*(U_i-U_i\cap
N(\Sigma_{i-1});\bar{\mc{G}})$, the second isomorphism induced
by the
stratum-preserving homotopy equivalences noted above. But $U_i-U_i\cap
 N(\Sigma_{i-1})=\Sigma_i-(\Sigma_i\cap N(\Sigma_{i-1}))$ is a
closed subcomplex of $S^n$ and, in particular, a finite
complex. Thus, the homology module $ H_*(U_i-U_i\cap
N(\Sigma_{i-1});\bar{\mc{G}})$ can be calculated as the simplicial homology
with
local coefficients of a finite complex. Therefore, it is a finitely
generate module if the fiber $G$ of $\mc{G}$ is a finitely generated
module over a Noetherian ring. 

We have shown:

\begin{lemma}
If $\mc{G}$ is a local coefficient system on $X_i$ whose fiber $G$ is a finitely
generated module over a Noetherian ring, then $H_*(X_i;\mc{G})$ is a finitely generated module.
\end{lemma}

\begin{corollary}
$\bar E^2_{p,q}=H_p(X_{i,k}; \mc{IH}_q(L_{i,k};\vg))$ and
$E^2_{p,q}=H_p(X_{i,k}; \mc{IH}_q(cL;\vg))$ are finitely generated
$\Gamma$-torsion modules for all $(p,q)$.
\end{corollary}
\begin{proof}
It suffices to show this for the $\bar E^2_{p,q}$ since, as previously
noted,
each $E^2_{p,q}$ is equal to either $\bar E^2_{p,q}$ or zero.

Since $L_{i,k}$ is a knot pair (possibly non-locally-flat), $
IH_q(L_{i,k};\vg)$ is a finitely generated $\Gamma$-torsion module since, by
\cite[Propositions 2.2 and 2.4]{CS}, it is finitely generated as a $\Q$-vector space. It
then follows from the preceding lemma that
$H_p(X_{i,k}; \mc{IH}_q(L_{i,k};\vg))$ is finitely generated as a 
$\Gamma$-module and as a $\Q$ vector space and hence is a finitely generated
$\Gamma$-torsion module.

\end{proof} 

We can now say something about prime divisors of the polynomials of the modules
$ IH_j^{\bar p}(W_i;\vg)$.

\begin{lemma} 
Let $w_{ij}$ be the polynomial of the
$\Gamma$-module $IH_j^{\bar p}(W_i;\vg)$,
and let $\xi_{ikj}$ denote the $j$th intersection
Alexander polynomial of the link $L_{i,k}$ of the $k$th connected
component
$X_{i,k}$ of $X_i$. 
Suppose that $\gamma$ is a prime element of $\Gamma$. Then
$\gamma|w_{ij}$ only if one of the following holds:
\begin{enumerate}
\item $\gamma|w_{i+1,j}$,
\item $\gamma|\xi_{iks}$ for some some $k$ and for some $s$ such that $0\leq j-s
\leq i-1$ and $0\leq s < n-i-2$.
\end{enumerate}
Furthermore, if $\gamma\nmid w_{i+1,j}$ and, for each $k$,
$\gamma|\xi_{iks}$ only if $s<n-i-1-\bar p(n-i)$, then
$\gamma\nmid w_{ij}$.
\end{lemma}
\begin{proof}
For $W_{n-1}$, this holds vacuously because $W_n=W_{n-1}$ and $X_{n-1}$ is
empty. Also each $w_{n-1,j}$ is a well-defined polynomial since, by Lemma
\ref{L: IH(W)}, $IH_*^{\bar p}(W_n-W_n\cap \Sigma;\vg)\cong H_*(S^n-K;\vg)$, but $H_*(S^n-K;\vg)$ is the ordinary Alexander
module of the knot. 

The proof now proceeds by downward induction on $i$. Assuming that the statement
is true for $i+1$ and using the Mayer-Vietoris sequences \eqref{E: W sequence},
the lemma follows for $i$ as in the proofs of Theorems \ref{T: at
one} and \ref {T: divisors}. Although sheaf intersection cohomology is used
there, the same spectral sequence arguments concerning prime divisors holds
here, using
instead the spectral sequences of  Theorem \ref{P: Neigh IH} and the map
between sequences given there. It is only necessary to use homological
indexing throughout, instead of shifting to cohomological indexing and back
as was
done in the proofs of the theorems. That each $w_{ij}$ is well-defined as 
the polynomial of a finitely generated torsion $\Gamma$-module follows from
the fact that each of the other terms of the Mayer-Vietoris sequence is
a
finitely-generated torsion $\Gamma$-module either by induction or by the
proofs of Theorems \ref{T: at one} and \ref {T: divisors}.
\end{proof}

\begin{theorem}\label{T: divisors2}
Let $\xi_{iks}$ denote the $s$th intersection
Alexander polynomial of the link $L_{i,k}$ of the $k$th connected
component of 
$\Sigma_i-\Sigma_{i-1}$. A prime element $\gamma\in \Gamma$ divides the intersection
Alexander polynomial $I\lambda_j^{\bar p}$ only 
if $\gamma|\lambda_j$ or $\gamma|\xi_{iks}$ for some some set of indices $i$, $k$,
and $s$
such that $0\leq j-s \leq i-1$ and $0\leq s < n-i-2$. $\gamma\nmid
I\lambda_j^{\bar p}$ if, for each $i$ and $k$, $\gamma|\xi_{iks}$ only if $s<n-i-1-\bar
p(n-i)$.
\end{theorem}
\begin{proof}
Since $W_0=S^n$, this follows from the lemma by induction and the fact noted in its proof that
$IH_*^{\bar p}(W_n;\vg)\cong H_*(S^n-K;\vg)$.
\end{proof}

We can also say something about the maximum power to which a prime divisor of the
intersection Alexander polynomial can occur. Once again, suppose that $\gamma$ is a
prime
element of $\Gamma$. Let $\gamma_{ipq}$ be the maximum power to which $\gamma$ occurs as
a divisor of
the polynomial $e_{ipq}$ of 
$H_p(X_{i}; \mc{IH}_q(L;\vg))=\oplus_k H_p(X_{i,k}; \mc{IH}_q(L_{i,k};\vg))$,
i.e. $\gamma^{\gamma_{ipq}}|e_{ipq}$, but $\gamma^{\gamma_{ipq}+1}\nmid
e_{ipq}$. Note that $\gamma_{ipq}=0$ if $X_i=\emptyset$. Let
$\gamma_l$ denote the maximum power to which $\gamma$ occurs in the Alexander polynomial
$\lambda_l$ of the knot $K$, and let $\gamma_{ij}$ denote the maximum power to which
$\gamma$ occurs in $IH_j^{\bar p}(W_i;\vg)$. 

\begin{lemma}
The prime factor $\gamma$ cannot occur in the polynomial $w_{ij}$ to a power greater than
\begin{equation*}
N=\gamma_{i+1,j}+(\sum_{\overset{p+q=j}{
q=0, q<n-i-1-\bar p(n-i)}}\gamma_{ipq})+(\sum_{p+q=j-1}\gamma_{ipq}),
\end{equation*}
i.e. $\gamma^{N+1}\nmid w_{i,j}$.

\end{lemma}
\begin{proof}
This again follows from the fact that the polynomial $w_{ij}$ of $ IH_j^{\bar
p}(W_i;\vg)$ must divide the product of the polynomials of 
$IH_j^{\bar p}(W_{i+1};\vg)$, $ IH_j^{\bar p}(N'_i;\vg)$,  and
$IH_{j-1}^{\bar
p}(N'_i-X_i;\vg)$, so  the power of $\gamma$ occurring in $w_{ij}$ must
be bounded by the
sum of the powers to which it occurs in the other three polynomials. Hence,
the summand 
$\gamma_{i+1,j}$ of $N$ enters trivially, and it only remains to show that
the powers of
$\gamma$ in the polynomials of $ IH_j^{\bar p}(N'_i;\vg)$ and $IH_{j-1}^{\bar
p}(N'_i-X_i;\vg)$ are bounded by $\sum_{\overset{p+q=j}{
q=0, q<n-i-1-\bar p(n-i)}}\gamma_{ipq}$ and $\sum_{p+q=j-1}\gamma_{ipq}$,
respectively.

Let us first consider $IH_{*}^{\bar p}(N'_i-X_i;\vg)$. It is the direct sum over $k$ of
the
modules  $IH_{*}^{\bar p}(N'_{ik}-X_{i,k};\vg)$, which are the abutments of spectral
sequences with $E^2_{p,q}$ terms given by
$H_p(X_{i,k}; \mc{IH}_q(L_{i,k};\vg))$. Since each 
term of the spectral sequence $E^r_{p,q}$ is a quotient of a submodule of
$E^{r-1}_{p,q}$, the power of $\gamma$ occurring in the polynomial of
$E^r_{p,q}$, $r\geq 2$,  
must be
less than or equal to that occurring in the polynomial of $E^2_{p,q}$. But
we can also see from the
series of short exact sequences \eqref{E: filter} (suitably dualized for a homological
spectral sequence) that the polynomial of the $l$th grade of the abutment is the product
of the polynomials of the terms $E^{\infty}_{p,l-p}$. Therefore, the power of $\gamma$
in the polynomial of $IH_{l}^{\bar p}(N'_{ik}-X_{i,k};\vg)$ must be less
than or equal
to the sum of the
powers of $\gamma$ in $H_p(X_{i,k}; \mc{IH}_{l-p}(L_{i,k};\vg))$, where the sum is taken
over $p$.
Furthermore, since $IH_{l}^{\bar p}(N'_{ik}-X_{i,k};\vg)$ is
the direct sum over $k$ of the $IH_{l}^{\bar p}(N'_{ik}-X_{i,k};\vg)$, the
power of $\gamma$ in
the polynomial of the former is equal to the sum of the powers in the polynomials of the
latter. Thus the power of $\gamma$ in the polynomial of  $IH_{l}^{\bar
p}(N'_{i}-X_{i};\vg)$ is less than or equal to the sums of the powers of $\gamma$ in
$H_p(X_{i,k}; \mc{IH}_{l-p}(L_{i,k};\vg))$, the sum being taken over $p$ and $k$. But for
each fixed $p$, the sum over $k$ gives the power $\gamma_{ip,l-p}$ of
$\gamma$ in the polynomial of $H_p(X_{i}; \mc{IH}_{l-p}(L_{i};\vg))$,
so the entire sum is
$\sum_{p+q=l}\gamma_{ipq}$. Therefore, the desired bound for this term holds by
taking
$l=j-1$. 

The bound for the power of $\gamma$ in the polynomial of $ IH_j^{\bar p}(N'_i;\vg)$ is
determined in the same manner once we have, again, observed that the $E^2_{p,q}$ terms of
the
spectral sequence that compute it are the same as the $E^2_{p,q}$ terms of the spectral
sequence for $IH_{*}^{\bar p}(N'_{i}-X_{i};\vg)$ when $q=0$ or $q<n-i-1-\bar p(n-i)$
and are equal to $0$ otherwise. 
\end{proof}

\begin{theorem}\label{T: max power}
The prime $\gamma\in \Gamma$ cannot occur in the polynomial $I\lambda_j^{\bar p}$ to a
power greater than
\begin{equation*}
\gamma_j+\sum_{i=0}^{n-2}\left((\sum_{\overset{p+q=j}{
q=0, q<n-i-1-\bar p(n-i)}}\gamma_{ipq})+(\sum_{p+q=j-1}\gamma_{ipq})\right).
\end{equation*}
\end{theorem}
\begin{proof}
This follows by an induction on the preceding lemma and the facts that $W_0=S^n$ and
$IH_*^{\bar p}(W_n;\vg)=IH_*^{\bar p}(W_{n-1};\vg)\cong H_*(S^n-K;\vg)$.
\end{proof}

Lastly, we can relate the divisors of $I\lambda_j^{\bar p}$ to the divisors of
the ordinary Alexander polynomials $\zeta_{iks}$ of the link knots
$L_{ik}$. We have seen in Theorem \ref{T: divisors2} that a prime $\gamma\in
\Gamma$ can divide $I\lambda_j^{\bar p}$ only if it divides the ordinary 
Alexander polynomial $\lambda_j$ of $K$ or one the intersection Alexander
polynomials $\xi_{iks}$ of the link $L_{ik}$, for $s$ in a certain
range. But then, again by Theorem \ref{T: divisors2}, $\gamma$ can divide
$\xi_{iks}$ only if it divides the  ordinary Alexander polynomial of
the link knot or
the intersection Alexander polynomial of one of \emph{its} links. However, in the
stratified pseudomanifold $S^n$, the link of a stratum of a link is also a 
link of
the original pseudomanifold. To see this, recall that if $L_{ik}$ is the link of
the $k$th connected component of  the stratum  $\Sigma_i-\Sigma_{i-1}$, then 
for some
point $x\in\Sigma_i-\Sigma_{i-1}$ there exists a topological neighborhood of $x$
in $S^n$ which is PL-homeomorphic to $D^i\times cL_{ik}$ and such that the
filtration
of $L_{ik}$ and  that of $D^i\times cL_{ik}$ induced by the cone and
product
filtrations are the same as those  induced by the restriction of the
filtration
on $S^n$. Similarly, if the induced stratification on $L_{ik}$ is given by the
filtration $\{T_l\}$ and $\mf{L}$ is the link of the $j$th stratum of
$L_{ik}$, then for some
point $y\in T_j-T_{j-1}$ there exists a topological neighborhood of $y$
in $L_{ik}$ which is PL-homeomorphic to $D^j\times c\mf{L}$ and such that the
filtration
of $\mf{L}$ and  $D^j\times \mf{L}$  is the same as that induced by
restriction of the
filtration on $L_{ik}$. But now consider $D^i\times \R\times  D^j\times
c\mf{L}\cong D^{i+j+1}\times c\mf{L}$ as a PL-subspace of $D^i\times
\R\times L_{ik}\cong D^i\times (cL_{ik}-*)$. Since
the map of $D^i\times cL_{ik}$ to a neighborhood of $x$ is a PL-homeomorphism, its
restriction to $D^{i+j+1}\times c\mf{L}$ is also a PL-homeomorphism into
its image in $S^n$. In
particular, since the dimension of $D^{i+j+1}\times c\mf{L}$ is clearly $n$ by
its construction, this set is a neighborhood of the image of $y$ under the
homeomorphism. Furthermore, the filtrations are all compatible so that the
filtration of $D^{i+j+1}\times c\mf{L}$ as a product of a cone of a filtered
space must be the same as the restriction filtration induce by its inclusion in
$S^n$. Therefore,  $D^{i+j+1}\times c\mf{L}$ is a distinguished neighborhood of
the image of $y$, and $\mf{L}$ is in fact the link of one of the connected
components of the stratum $\Sigma_{i+j+1}-\Sigma_{i+j}$ in $S^n$. 

It now follows from Theorem \ref{T: divisors2} that $\gamma$ divides the intersection Alexander polynomial of 
the knot pair $L_{ik}$ only if it divides the ordinary Alexander polynomial of
that knot pair or one of the intersection Alexander polynomials of another link of a
stratum of $S^n$ with smaller codimension. Thus, by induction and
Theorem \ref{T: divisors2}, we have the following: 

\begin{theorem}\label{T: divisors3}
Let $\zeta_{iks}$ be the $s$th ordinary Alexander polynomial of the link
knot pair
$L_{ik}$. A prime element $\gamma\in \Gamma$ divides the $s$th 
intersection
Alexander polynomial $I\lambda_j^{\bar p}$ only
if $\gamma|\lambda_j$ or $\gamma|\zeta_{iks}$ for some set of indices $i$, $k$,
and $s$,
such that $0\leq j-s \leq i-1$ and $0\leq s < n-i-2$. 
\end{theorem}
\begin{proof}
It only remains to prove that this range of indices is the correct one. By
Theorem \ref{T: divisors2}, we know that $\gamma|I\lambda_j^{\bar p}$ only if
$\gamma|\lambda_j$ or $\gamma|\xi_{iks}$ for indices $i$, $k$, $s$ such that
$0\leq j-s \leq i-1$ and $0\leq s < n-i-2$ (call this Index
Condition 1 for $i$ and $s$). Similarly, applying Theorem \ref{T: divisors2} to
the
link knot $L_{ik}$ of dimension $n-i-1$, we know that $\gamma$ divides the
$s$th intersection Alexander
polynomial $\xi_{iks}$ of $L_{i,k}$ only if it divides the $s$th ordinary
Alexander polynomial of
$L_{ik}$, $\zeta_{iks}$, or the $r$th intersection Alexander polynomial of some
link of the
link, say $\mf{L}_{al}$, for some set of indices $a$, $l$, and $r$
satisfying $0\leq s-r\leq a-1$ and $0\leq r<n-i-1-a-2=n-(i+a+1)-2$ (call
this Index Condition 2 for $a$ and $r$). But by the discussion preceding the
statement of the
theorem, we know that $\mf{L}_{al}=L_{i+a+1,m}$ for some $m$. We will show that
if a collection of indices $i$, $s$, $a$, and $r$ satisfy Index Condition 1
for $i$ and $s$ and
Index Condition 2 for $a$ and $r$, then Index Condition 1 is satisfied for $i+a+1$ and
$r$,
i.e. with $i$ and $s$ replaced by $i+a+1$ and
$r$ (note that $j$ and $n$ are fixed throughout and the indexing of connected
components is irrelevant). This will imply that the relevant collection of
$r$th
intersection Alexander polynomials of the $\mf{L}_{al}$ will already have been
included among the collection $\zeta_{iks}$ satisfying Index Condition 1 for $i$
and $s$.  Hence, by an induction, we can conclude that $\gamma|I\lambda_j^{\bar
p}$ only
if it divides $\lambda_j$ or the ordinary Alexander polynomial $\zeta_{iks}$ for
$i$ and $s$ satisfying Index Condition 1.

We now prove the claim on the Index Conditions: The second part of
Index Condition 1 for $i+a+1$ and $r$ is exactly the second part of Index
Condition 2 for $a$ and $r$. From Index Condition 1 for $i$ and $s$, $0\leq
j-s\leq i-1$,
and by Index Condition 2 for $a$ and $r$, $0\leq s-r\leq a-1$. Adding these
inequalities gives  $0\leq j-r\leq a+i-2$, and certainly $a+i-2\leq a+i$.
Therefore $0\leq j-r\leq a+i$, which is Index Condition 1 for $i+a+1$ and $r$.
This completes the proof.
\end{proof}

\begin{remark}
Note that the results of this section seem to depend upon the particular choice of
stratification of the knot $K$. Hence, it is conceivable that specific choices of
stratification might yield more precise information. In particular, it might be
possible to
obtain extra information by making clever choices of stratification dependent
upon the specific prime $\gamma$ under discussion. 
\end{remark}

\bibliographystyle{amsplain}
\bibliography{/home/accts/friedman/docs/bib}

Several diagrams in this paper were typeset using the\TeX\, commutative
diagrams package by Paul Taylor.   

\end{document}